\documentclass{amsart}
\usepackage[left=2.65cm,right=2.65cm,top=2.7cm,bottom=2.7cm]{geometry}
\usepackage{amsmath,amssymb,amsthm,mathrsfs,mathtools,thmtools}   
\usepackage{hyperref}
\usepackage[capitalise,noabbrev]{cleveref}
\usepackage{dsfont}  
\usepackage{enumerate}
\usepackage{esint} 
\usepackage{comment}
\numberwithin{equation}{section}

\usepackage{xcolor}

\hypersetup{colorlinks=true, 	linkcolor=blue!75!black, 	citecolor=red!50!black, 	urlcolor=green!50!black, 	linktoc=all }

\def\NN{\mathbb{N}}
\def\RR{\mathbb{R}}
\def\ZZ{\mathbb{Z}}
\def\TT{\mathbb{T}}

\newtheorem{theorem}{Theorem}[section]
\newtheorem{proposition}[theorem]{Proposition}
\newtheorem{corollary}[theorem]{Corollary}

\newtheorem{lemma}[theorem]{Lemma}
\theoremstyle{definition}
\newtheorem*{definition*}{Definition}
\newtheorem{definition}[theorem]{Definition}

\theoremstyle{remark}
\newtheorem{remark}[theorem]{Remark}

\newcommand{\norm}[1]{\left\|#1\right\|}
\newcommand{\abs}[1]{\left\vert#1\right\vert}

\newcommand{\normxt}[2]{\left\|#2\right\|_{#1}}
\newcommand{\normxtjm}[2]{\left\|#2\right\|_{#1,\hspace{0.5pt}jm}}
\newcommand{\seminormxt}[2]{\left[#2\right]_{C^{#1}_{x,t}}}
\newcommand{\normx}[2]{\left\|#2\right\|_{C^{#1}_{x}}}

\DeclareMathOperator\Div{div}
\DeclareMathOperator\curl{curl}
\DeclareMathOperator\Id{\rm Id}

\DeclareMathOperator\dist{dist}
\DeclareMathOperator\supp{supp}
\DeclareMathOperator\tr{tr}

\newcommand\compt{\diamond}

\usepackage{units} 

\newcommand{\normatres}[1]{{\left\vert\kern-0.25ex\left\vert\kern-0.25ex\left\vert #1 
		\right\vert\kern-0.25ex\right\vert\kern-0.25ex\right\vert}}

\begin{document}

\title[H\"older continuous dissipative solutions to ideal MHD]{H\"older continuous dissipative solutions of  ideal MHD\\  with nonzero helicity}

\author{Alberto Enciso}
\address{Instituto de Ciencias Matem\'aticas, Consejo Superior de
	Investigaciones Cient\'\i ficas, 28049 Madrid, Spain}
\email{aenciso@icmat.es}

\author{Javier Pe\~nafiel-Tom\'as}
\address{Instituto de Ciencias Matem\'aticas, Consejo Superior de
	Investigaciones Cient\'\i ficas, 28049 Madrid, Spain}
\email{javier.penafiel@icmat.es}

\author{Daniel Peralta-Salas}
\address{Instituto de Ciencias Matem\'aticas, Consejo Superior de
	Investigaciones Cient\'\i ficas, 28049 Madrid, Spain}
\email{dperalta@icmat.es}

\newcommand{\cambios}[1]{\textcolor{Caca}{#1}}

%
%
\begin{abstract}
	We prove the existence of weak solutions  to the 3D ideal MHD equations, of class $C^\alpha$ with $\alpha=1/200$, for which the total energy and the cross helicity (i.e., the so-called Els\"asser energies) are not conserved. The solutions do not possess any symmetry properties and the magnetic helicity, which is necessarily conserved for H\"older continuous solutions, is nonzero. The construction, which works both on the torus~$\TT^3$ and on~$\RR^3$ with compact spatial support, is based on a novel convex integration scheme in which the magnetic helicity is preserved at each step. This is the first construction of continuous weak solutions at a regularity level where one conservation law (here, the magnetic helicity) is necessarily preserved while another (here, the total energy or cross helicity) is not, and where the preservation of the former is nontrivial in the sense that it does not follow from symmetry considerations.
\end{abstract}

\maketitle

\setcounter{tocdepth}{1}
\tableofcontents

\section{Introduction}

The 3D ideal magneto-hydrodynamics (MHD) system, which couples the Euler equations for an ideal fluid with Maxwell's equations of electromagnetism, models the macroscopic behavior of plasmas and other electrically conducting fluids. In terms of the {\em velocity}~$v(x,t)$ and the {\em magnetic field}~$B(x,t)$, which we can understand as time-dependent solenoidal vector fields in three dimensions, the MHD equations read as 
\begin{subequations}
\label{MHD}
\begin{align}
	& \partial_t v+v\cdot\nabla v-B\cdot\nabla B+\nabla p=0, \label{MHD a} \\
	& \partial_t B+v\cdot\nabla B-B\cdot\nabla v=0,  \label{MHD b}\\
	& \Div v=\Div B=0. \label{MHD c}
\end{align}
\end{subequations}
Here the scalar function $p(x,t)$ is the {\em pressure}\/, and can be written in terms of~$(v,B)$ by solving an elliptic equation, just as in the case of the Euler equations. For a preliminary discussion in this introduction, we can assume that the spatial variable takes values in $\TT^3:=(\RR/2\pi \ZZ)^3$. However, all along this article we will consider the somewhat harder case of compactly supported solutions on~$\RR^3$, which of course yields solutions on $\TT^3$ and on bounded domains.

\subsection{Conservation laws for the MHD system}

The MHD system has a number of conservation laws. The most basic ones are the mean velocity and magnetic field, which are preserved even for weak solutions of the equations. For this reason, we will always consider solutions satisfying
\[
\int_{\TT^3}v(x,t)\, dx=\int_{\TT^3}B(x,t)\, dx=0
\]
for all~$t$.

Much more interesting conservation laws are the {\em total energy}\/, the {\em magnetic helicity}\/, and the {\em cross helicity}\/, respectively defined as
\begin{align}\label{E.E_intro}
	\mathcal E(t)&:=\int_{\TT^3} \left(|v(x,t)|^2+|B(x,t)|^2\right) dx\,,\\
	\mathcal H(t)&:=\int_{\TT^3} (\curl^{-1}B)(x,t)\cdot B(x,t)\, dx\,,\\
	\mathcal H_\times(t)&:=\int_{\TT^3} B(x,t)\cdot v(x,t)\, dx\,.\label{E.Hx_intro}
\end{align}
Since~$B$ has zero mean, here one can take $\curl^{-1}B:=(-\Delta)^{-1}\curl B$, which is the only zero mean, divergence-free field whose curl is~$B$. More generally, replacing this field by any other vector potential of~$B$ does not change the value of the integral that defines the magnetic helicity. In the case that $B$ is a compactly supported vector field on $\mathbb R^3$, its $\curl^{-1}$ is defined using the Biot-Savart operator, see e.g.~\cite{PNAS16}.

The total energy, which is also the Hamiltonian of the system, leads to the consideration of solutions in the energy class $(v,B)\in L^\infty L^2$. The cross helicity is also defined at this regularity level, while the magnetic helicity remains well defined under the weaker condition $B\in L^\infty\dot H^{-\frac12}$. This is ultimately the reason for which magnetic helicity plays a different role than the total energy and the cross helicity (i.e., the so-called Els\"asser energies) in the context of non-smooth solutions.

\subsection{Non-conservative solutions of ideal MHD}

While smooth solutions to the MHD system do conserve the three quantities~\eqref{E.E_intro}--\eqref{E.Hx_intro}, this is not necessarily true in the case of rough solutions. In the closely related setting of the 3D incompressible Euler equations, 
\[
\partial_t u+u\cdot \nabla u+\nabla p=0\,,\qquad \Div u=0\,,
\]
which the MHD system reduces to when the magnetic field is identically~0, the situation is well understood by now. If $u$ is a weak Euler flow of class~$C^\alpha$, it is known that the energy $\mathcal{E}_{\mathrm E}(t):=\int_{\TT^3} |u(x,t)|^2\,dx$ is conserved~\cite{CET} if $\alpha>\frac13$, and that this is not necessarily the case if $\alpha<\frac13$, as shown in the proof of the flexible side of Onsager's conjecture by Isett~\cite{Isett2018} and by Buckmaster, De Lellis, Szek\'elyhidi and Vicol~\cite{BLSV} for dissipative solutions. 

The latter works hinge on Nash-type convex integration methods, which were first employed in the context of non-conservative Euler flows in the seminal work of De Lellis and Sz\'ekelyhidi~\cite{Continuous} and to which a number of authors subsequently made significant contributions~\cite{DeLellisSzekelyhidi12, BuckmasterDeLellisIsettSzekelyhidi15, Daneri17}. Previous examples of non-conservative weak solutions of the incompressible Euler equations had been obtained in~\cite{scheffer1993,shnirelman97,LS09}. Further results include the introduction of intermittent building blocks and their applications to the Navier--Stokes equations~\cite{BV, BuckmasterColomboVicol22}, flexibility results for 2D fluids~\cite{Giri}, and the analysis of other equations of fluid mechanics~\cite{SQG,CCF}. For details, we refer to the survey papers~\cite{LellisBAMS,VicolBAMS}.

An analogue of the Onsager conjecture for weak solutions to ideal MHD was proposed in~\cite{VicolBAMS}. The rigidity side of the problem is already well understood:  assuming $(v,B)\in C^\alpha$, Caflisch, Kappler and Steele have shown~\cite{ CKS} that the Els\"asser energies are conserved if $\alpha>\frac13$, and that the magnetic helicity is conserved under the weaker condition $\alpha>0$. Sharper estimates in $L^3$-based spaces are also available, including endpoint estimates. The latter mirror those for 3D Euler~\cite{CET, CheskidovConstantinFriedlanderShvydkoy08} in the case of the Els\"asser energies, and only require that $(v,B)\in L^3$ in the case of the magnetic helicity~\cite{KangLee07}.

The flexible side of the problem is far less clear cut. As emphasized in~\cite{VicolBAMS}, the key difficulty is that the fact that magnetic helicity is conserved by weak solutions under much milder regularity assumptions than the Els\"asser energies makes the construction of non-conservative solutions very challenging, so the construction of $C^\alpha$ weak solutions with $\alpha>0$ that preserve the magnetic helicity but not the Els\"asser energies requires the introduction of new methods. Morally, a similar difficulty arises in the SQG equation, where the kinetic energy is conserved for solutions with $C^{1/3}$ regularity but where the conservation of the Hamiltonian holds for continuous weak solutions. In this setting, Dai, Giri and Radu~\cite{DGR} and Looi and Isett~\cite{IsettSQG} have independently obtained discontinuous weak solutions for SQG with the sharp regularity which do not preserve the Hamiltonian. The construction of more regular solutions which necessarily preserve the Hamiltonian but not the kinetic energy remains wide open, since none of the existing convex integration schemes seem to yield solutions at a regularity level where one conservation law is necessarily preserved but not the other.

Concerning the MHD equations, the first nontrivial solutions for which the total energy is not conserved were constructed by Bronzi, Lopes Filho and Nussenzveig Lopes~\cite{BLN} using a class of symmetric weak solutions to MHD on~$\TT^3$ for which the problem boils down to the study of a $2\frac12$D Euler flow. Specifically, they exploited the observation that $x_3$-independent fields
\begin{equation}\label{E.intro_symm}
	v:=(v_1(x_1,x_2,t),v_2(x_1,x_2,t),0)\,,\qquad B:=(0,0,b(x_1,x_2,t))
\end{equation}
solve the MHD equations if and only if $u:=(v_1(x_1,x_2,t),v_2(x_1,x_2,t),b(x_1,x_2,t))$ is an $x_3$-independent Euler flow. Note that the magnetic helicity of symmetric solutions is automatically~0, and that this strategy cannot yield finite energy solutions on~$\RR^3$. Furthermore, the solutions constructed in~\cite{BLN} have zero cross helicity.

Faraco, Lindberg and Sz\'ekelyhidi~\cite{FLS} managed to adapt  the general setting of~\cite{LS09} to construct solutions $(v,B)\in L^\infty$ that do not preserve $\mathcal{E}$ and~$\mathcal{H}_\times$. Although nontrivial strict subsolutions had already been constructed in~\cite{FL}, the construction is highly nontrivial because the associated $\Lambda$-convex hull has empty interior. The solutions in~\cite{FLS} have zero magnetic helicity but they are not symmetric, so the analysis does not reduce to that of the Euler equations.  

The first non-conservative weak solutions with nonzero magnetic helicity were constructed by Beekie, Buckmaster and Vicol~\cite{BBV}. These solutions, of class $C^0H^\alpha$ for some small~$\alpha>0$, do not preserve {\em any}\/ of the quantities $\mathcal{E}$, $\mathcal{H}$, and~$\mathcal{H}_\times$. This is particularly interesting because, in view of Faraco and Lindberg's proof~\cite{FL} of Taylor's conjecture, these solutions cannot be obtained as weak ideal limits from the viscous and resistive MHD equations. The proof of this result uses an intermittent convex integration scheme adapted to the complicated structure of the MHD equations, which introduces serious difficulties. Later, Faraco, Lindberg and Sz\'ekelyhidi~\cite{FLS24} cleverly refined their approach in~\cite{FLS} to construct weak solutions $(v,B)\in L^\infty$ with nonzero magnetic helicity that do not preserve $\mathcal{E}$ or~$\mathcal{H}_\times$, as well as solutions $(v,B)$ in the sharp integrability class $L^\infty L^{3,\infty}$ that fail to preserve~$\mathcal E$, $\mathcal H$, and~$\mathcal H_\times$.

Nontrivial weak solutions to ideal MHD with the sharp H\"older regularity for which the Els\"asser energies are not conserved were recently constructed by Miao, Nie and Ye in the symmetric case~\cite{MNY}. More precisely, for any $\alpha<\frac13$, they constructed  $(v,B)\in C^\alpha$ of the form~\eqref{E.intro_symm} for which~$\mathcal{E}$ and~$\mathcal{H}_\times$ fail to be conserved. To prove this result, the authors formulate the problem in terms of solutions to the $2\frac12$D Euler equations as above and use the asymmetry of the errors to improve the Newton--Nash iteration scheme for 2D Euler developed by Giri and Radu~\cite{Giri}. This construction cleverly exploits that the magnetic helicity of a symmetric solution to ideal MHD is always zero to bypass the difficulties associated with the conservation of magnetic helicity for H\"older continuous solutions.

\subsection{Main result and strategy of the proof}

Our objective in this paper is to construct H\"older continuous solutions to ideal MHD with nonzero helicity that do not preserve the Els\"asser energies. This is the first construction of continuous weak solutions at a regularity level where one conservation law (here the magnetic helicity) is necessarily preserved while another (here the total energy or cross helicity) is not, and where the preserved quantity is nontrivial in the sense that its conservation is not a consequence of symmetry considerations.

More precisely, we prove the following:

\begin{theorem}
\label{main theorem}
Let $T>0$ and let $\alpha\coloneqq \frac{1}{200}$. Given a smooth divergence-free field $\bar{B}_0\in C^\infty_c(\RR^3,\RR^3)$ and $\varepsilon>0$, there exists a weak solution of ideal MHD $(v,B)\in C^\alpha(\RR^3\times[0,T])\cap C^0([0,T],L^2(\RR^3))$ such that: 
\begin{enumerate}
	\setlength\itemsep{5pt}
	\item the total energy is a strictly decreasing function of time,
	\item cross-helicity is not conserved,
	\item $\norm{B(\cdot,0)-\bar{B}_0}_{C^\alpha}<\varepsilon$.
\end{enumerate}
\end{theorem}
\begin{remark}
As mentioned above, magnetic helicity is conserved for H\"older continuous solutions. Therefore, choosing a field $\bar{B}_0$ with nontrivial helicity and taking $\varepsilon>0$ sufficiently small yields Hölder-continuous weak solutions of ideal MHD in which neither the total energy nor cross-helicity are conserved but with non-zero constant magnetic helicity. 
\end{remark}
\begin{remark}
The spatial support of the weak solution, $\supp\, (v,B)(\cdot,t)$, is contained in a bounded domain for all times $t\in[0,T]$. Hence, after a suitable rescaling, this result also provides infinitely many dissipative weak solutions of ideal MHD on $\TT^3$ and on bounded domains.
\end{remark}

Let us now give an informal idea of how the proof works. The details of the construction, modulo the technical results relegated to later sections, are given in Section~\ref{S.overview}. In view of the difficulties that convex integration schemes have faced in producing Hölder-continuous solutions for ideal MHD, our guiding principle is to employ a convex integration scheme that automatically preserves magnetic helicity. Since helicity is preserved under volume-preserving diffeomorphisms (see e.g.~\cite{PNAS16} and references therein), we thus start off the design of our scheme with the idea that, at each step, the magnetic field will be modified by pushing it forward with a (time dependent) volume-preserving diffeomorphism, thereby ensuring that the magnetic helicity is preserved throughout the process. Our general goal is to pick a diffeomorphism that involves high frequency oscillations, yet is  close to the identity in some suitable sense. The evolution equation for~$B$ shows that the choice of this diffeomorphism has a major impact on the structure of the velocity field~$v$. In our scheme, this diffeomorphism determines~$v$ through an ODE. 

From a technical point of view, the main challenges one finds in implementing this strategy---which are substantial---stem from the need to obtain good estimates for a velocity field that is defined essentially as the vector field generating certain rapidly oscillating volume-preserving diffeomorphism. Since we do not know how to construct a sufficiently rich family of explicit volume-preserving diffeomorphisms,  at each step our diffeomorphism~$\phi$ is obtained as the composition of an explicit main term~$\phi^0$ with a smaller correction diffeomorphism~$\phi^c$, which is needed to ensure that~$\phi$ is volume-preserving. The correction~$\phi^c$ is constructed by means of a fully nonlinear elliptic equation and is therefore very hard to analyze in sufficient detail. Our basic strategy is to choose a complicated, rather awkward class of main terms~$\phi^0$ for which one can show that~$\phi^c$ is suitably small. The price to pay is that, once we relate~$\phi$ to~$v$ through the ODE, a number of dangerous, highly nonstandard terms arise in the Reynolds stress, and they must be carefully estimated.

More precisely, at each iteration step, the diffeomorphism~$\phi$ induces a change in the velocity field which we decompose as the sum of an explicit, reasonably simple part, which one can eventually control for a suitable choice of~$\phi^0$, and another term that we call~$w_\phi$, which is highly nontrivial to estimate. Most of the choices made in this paper are motivated by the need to control this term, within the residual freedom granted by our general approach. The hardest estimates are those for $\mathcal R \partial_t w_\phi$, where the operator $\mathcal R$ should be understood as a sort of inverse of the divergence operator. Here we rely on a class of function spaces which allows us to describe functions whose time derivatives are comparable to their spatial derivatives, and which are well adapted to the structure of our system of equations. The estimates for $w_\phi$ take most of Sections~\ref{section estimates vector fields} and~\ref{section matrix}, which form the technical core of the paper.

Let us note in passing that volume-preserving diffeomorphisms are crucially used in another convex integration scheme we recently developed, motivated by the geometric problem of constructing stationary Euler flows that are topologically equivalent, in a weak sense, to a fixed divergence-free field~$v_0$ \cite{steady_Euler}. Diffeomorphisms enter these schemes in very different ways: whereas here they  essentially arise as non-autonomous flows of~$v$, in~\cite{steady_Euler} they are used to push forward~$v_0$ to define the velocity, which is pretty much what we do here with~$B$---{\em not}\/ with our main source of trouble, $v$. It is therefore not surprising that both the structure of the scheme and the estimates needed to close it are entirely different. However, we regard~\cite{steady_Euler} as a first proof of concept that volume-preserving diffeomorphisms can be effectively employed to construct convex integration schemes with novel features that are very hard to obtain using more traditional approaches.

\subsection{Organization of the paper} In~\cref{S.overview} we set the convex integration scheme that Theorem~\ref{main theorem} rests on, introducing a definition of subsolutions that captures the above discussion and stating the main estimates needed at each iteration step. The construction of subsolutions at each iteration step is presented in~\cref{section construction}. The necessary estimates for the diffeomorphism, for the velocity and for the Reynolds stress are respectively presented in Sections~\ref{section estimates diffeo}--\ref{section matrix}. Section~\ref{sec correction} introduces a final correction which allows us to control the spatial support of the subsolutions, completing the proof of \cref{prop steps}. Finally, in~\cref{section th main} we use \cref{prop steps} and a suitable choice of the initial subsolution to prove \cref{main theorem}. The paper concludes with three appendices containing some background material on Hölder and Besov spaces and several auxiliary estimates.



\section{Overview of our convex integration scheme}
\label{S.overview}

\subsection{Basic ideas}
We will now explain the basic ideas of our construction. Here we will focus on the heuristics and leave the details for later sections. A common factor in previous results concerning the construction of weak solutions of MHD is that they follow the standard ``hard analysis'' approach: at each step in the scheme, the fields $(v,B)$ are corrected simply by {adding} rapidly oscillating fields. 

Instead, our basic idea is to incorporate in our convex integration scheme the deep geometric meaning of Equation (\ref{MHD b}): by Alfv\'en's classical theorem, this equation encodes the fact that magnetic field lines are ``frozen in the fluid'' in the infinite conductivity limit. This is made more precise in the following lemma, whose proof is standard, see e.g.~\cite[Theorem 6.4.5]{AMR}. In the statement, $X\in C^\infty(\RR^3\times[0,T],\RR^3)$ is the (non-autonomous) flow of $v$, that is, the unique solution of
    \[\begin{cases}\partial_tX(x,t)=v(X(x,t),t), \\ X(x,0)=x,\end{cases}\]
    which satisfies that $X_t$ is a volume-preserving diffeomorphism of $\RR^3$ for all $t\in[0,T]$. Note that the conservation of magnetic helicity, which lies at the heart of the difficulties in constructing H\"older continuous solutions to ideal MHD, follows from this statement in the setting of smooth solutions, as helicity is preserved by volume-preserving diffeomorphisms (see e.g.~\cite{PNAS16} and references therein):

\begin{lemma}
\label{equivalencia ec B}
Two divergence-free fields $v,B\in C^\infty(\RR^3\times[0,T],\RR^3)$ satisfy (\ref{MHD b}) if and only if we can write
\[B(\cdot,t)=(X_t)_\ast\hspace{1pt} B(\cdot,0) \qquad \forall t\in[0,T],\]
where $X$ is the unique flow of $v$ and $Y_*w:=(DY w)\circ Y^{-1}$ denotes the pushforward of a vector field~$w$ by a diffeomorphism~$Y$.
\end{lemma}


Note that this property of the evolution of~$B$ mirrors that of the vorticity in the 3D incompressible Euler equations. The fact that this has never played any role in the construction of non-conservative weak Euler flows merely reflects that these flows are too rough for the vorticity to be useful. In ideal MHD, however, the situation is completely different.

Given a velocity $v$, Lemma~\ref{equivalencia ec B} provides an easy way to construct magnetic fields $B$ so that (\ref{MHD b}) holds.  However, it is not so easy to construct $v$ so that (\ref{MHD a}) holds, which is why we resort to convex integration. 

\cref{main theorem} is proved using an iterative scheme. At the beginning of the process we will have fields $(v_0, B_0)$ which are not a solution of (\ref{MHD}) but instead solve a ``relaxed'' system. Our goal is to perform successive high-frequency perturbations to these fields in order to ``reduce the error''. To make this precise, we define an appropriate notion of subsolutions. Here and in what follows, $\RR^{3\times3}_\text{sym}$ denotes the space of symmetric $3\times 3$ matrices and, given a map $f\in C^\infty(\RR^3\times[0,T])$, we set $f_t\coloneqq f(\cdot,t)\in C^\infty(\RR^3)$.
\begin{definition}
	\label{def subsolution}
	Given vector two fields $v,B\in C^\infty(\RR^3\times [0,T],\RR^3)$, a scalar field $p\in C^\infty(\RR^3\times [0,T],\RR)$ and a symmetric matrix $R\in C^\infty(\RR^3\times [0,T],\RR^{3\times3}_\text{sym})$, the 4-tuple $(v,B,p,R)$ is said to be a {\em subsolution}\/ if
	\begin{subequations}
		\label{subs}
		\begin{align}
			& \partial_t v+\Div(v\otimes v-B\otimes B)+\nabla p=\Div R, \label{subs a} \\ 
			& B(\cdot,t)=(X_t)_\ast\hspace{1pt}B(\cdot,0), \label{subs c}\\
			& \Div v=\Div B=0, \label{subs d}
		\end{align}
	\end{subequations}
As before, $X_t$ is the (non-autonomous) flow of $v$.	
\end{definition}

It follows from \cref{equivalencia ec B} that a subsolution will always satisfy (\ref{MHD b}) and (\ref{MHD c}). Furthermore, since $v$ and $B$ are divergence-free, it is easy to check that (\ref{subs a}) with $R=0$ implies (\ref{MHD a}). Hence, the matrix $R$, which is usually called {\it Reynolds stress}, measures the error of $(v,B,p)$. If it were identically zero, then $(v,B,p)$ would be a solution to (\ref{MHD}). Our goal is to construct a sequence of subsolutions $\{(v_q,B_q,p_q,R_q)\}_{q=0}^\infty$ such that $R_q$ converges uniformly to 0. Meanwhile, $v_q$ and $B_q$ will converge uniformly to some continuous fields $v$ and $B$, which will be a weak solution of (\ref{MHD}). 

Let us sketch how to construct the sequence of subsolutions. We will only informally explain the main ideas, as the actual construction will be considerably more involved. Given the subsolution $(v_q,B_q,p_q,R_q)$, we define the new velocity $v_{q+1}$ by prescribing its flow. Denoting by $X^q$ the flow of $v_q$, we define a new map $X^{q+1}\in C^\infty (\RR^3\times[0,T],\RR^3)$ as 
\[X^{q+1}_t\coloneqq \phi^{q+1}_t\circ X^q_t\circ \left(\phi^{q+1}_0\right)^{-1},\]
where $\phi^{q+1}\in C^\infty(\RR^3\times[0,T],\RR^3)$ is a map such that $\phi^{q+1}_t:\RR^3\to\RR^3$ is a volume-preserving diffeomorphism for each $t\in [0,T]$. Taking into account that $X^q_t$ is volume-preserving because it is generated by a divergence-free field, it follows from this definition that $X^{q+1}_t$ is volume-preserving for all times $t\in[0,T]$. Although we will explain how to choose $\phi^{q+1}$ later, at this point we can anticipate it will be roughly of the form
\begin{equation}
\phi^{q+1}(x,t)\approx x+\sum_{j,m} \frac{a_{jm}(x,t)}{\eta}\hspace{1pt}\zeta_j\sin[\theta_{jm}(x,t)],
\label{expresión phi intro}
\end{equation}
where $\zeta_j\in \RR^3$ with $j=1, \dots 6$ are fixed vectors, $a_{jm}$ are slowly varying scalar amplitudes, $\eta\gg1$ is a large parameter and the phases $\theta_{jm}$ encode oscillations with frequency $\lambda \gg \eta$. 

The new velocity is then defined so that it generates the flow $X^{q+1}$. Setting 
\begin{equation}
(v_{q+1})_t\coloneqq [\partial_t \phi^{q+1}_t+D\phi^{q+1}_t(v_q)_t]\circ (\phi^{q+1}_t)^{-1},
\label{def vq+1 intro}
\end{equation}
a direct computation shows
\begin{align}
	(v_{q+1})_t\circ X_t^{q+1}&=[\partial_t \phi^{q+1}_t+D\phi^{q+1}_t(v_q)_t]\circ X^q_t\circ (\phi^{q+1}_0)^{-1}\nonumber\\&=\left[(\partial_t\phi^{q+1}_t)\circ X^q_t+(D\phi^{q+1}_t\circ X^q_t)\,\partial_t X^q_t\right]\circ (\phi^{q+1}_0)^{-1}\label{vq+1 genera Xq+1}\\&=\partial_t(\phi^{q+1}_t\circ X^q_t)\circ (\phi^{q+1}_0)^{-1}=\partial_tX^{q+1}_t, \nonumber
\end{align}
where in the second equality we have used the fact that $X^q$ is the flow of $v_q$. Since $X^{q+1}_0=\Id$, by construction, we conclude that $X^{q+1}$ is the flow of $v_{q+1}$. Taking into account that $X^{q+1}_t$ is volume-preserving, we see that $v_{q+1}$ must be divergence-free. Again, although its precise form will be discussed later, we anticipate that the perturbation $w_{q+1}\coloneqq v_{q+1}-v_q$ will behave as
\begin{equation}
w_{q+1}\approx \sum_{j,m} a_{jm}\hspace{1pt}\zeta_j\cos(\theta_{jm}).
\label{approx wq+1}
\end{equation}

Next, we define the new magnetic field as
\[(B_{q+1})_t\coloneqq (\phi^{q+1}_t)_\ast\, (B_q)_t,\]
which guarantees that (\ref{subs c}) is satisfied. Indeed, using that (\ref{subs c}) holds for the previous subsolution, we have
\begin{align}
\begin{split}
(B_{q+1})_t&=(\phi^{q+1}_t)_\ast\,(X^q_t)_\ast\, (B_q)_0=(\phi^{q+1}_t)_\ast\,(X^q_t)_\ast\,[(\phi^{q+1}_0)^{-1}]_\ast\, (B_{q+1})_0 \\&= (X^{q+1}_t)_\ast\,(B_{q+1})_0,
\label{pushforward B intro}
\end{split}
\end{align}
as we wanted. Here we have used that the pushforward of a composition of diffeomorphisms is the composition of the pushforward associated to the individual diffeomorphism; this, in particular, implies that the pushforward commutes with taking the inverse. Since (\ref{subs c}) holds, it follows from \cref{equivalencia ec B} that (\ref{MHD b}) is automatically satisfied. We anticipate here that $\phi^{q+1}$ is chosen so that the perturbation $b_{q+1}\coloneqq B_{q+1}-B_q$ is very small.

Finally, $p_{q+1}$ and $R_{q+1}$ are  defined so that (\ref{subs a}) holds. Taking into account that (\ref{subs a}) holds for the previous subsolution, by hypothesis, it is easy to check that the new Reynolds stress must satisfy
\begin{align}
\nonumber
\Div R_{q+1}&=w_{q+1}\cdot\nabla v_q \\
\begin{split}
	\label{ec R intro}
	&\hspace{12pt}+\partial_t w_{q+1}+v_q\cdot\nabla w_{q+1}\\&\hspace{12pt}-\Div(b_{q+1}\otimes B_{q}+B_{q}\otimes b_{q+1}+b_{q+1}\otimes b_{q+1})
\end{split}
\\&\hspace{12pt}+\Div\left(w_{q+1}\otimes w_{q+1}+R_q+(p_{q+1}-p_q)\Id\right). \nonumber
\end{align}
Our goal is to construct a new subsolution such that $R_{q+1}$ is smaller than $R_q$. Fortunately, when solving this equation for $R_{q+1}$ we can expect a smoothing effect, so that high frequencies are attenuated (see \cref{appendix divergence}). Hence, there is no issue if the terms on the right-hand side are quite large: as long as their amplitude is much smaller than their frequency, their contribution to the Reynolds stress will be small. 

Taking into account that the perturbation $w_{q+1}$ will oscillate much faster than $v_q$, we see that the first term on the right-hand side is quite harmless, since the derivative falls on the low-frequency terms. Meanwhile, the second term is more problematic because the derivatives affect the high-frequency term. To prevent the amplitude from being unacceptably large, we must impose the cancelation
\begin{equation}
\abs{\partial_t\theta_{jm}+v_q\cdot\nabla\theta_{jm}}\ll \lambda,
\label{cancelaciones theta intro 1}
\end{equation}
where $\lambda$ is the frequency of the oscillations. This cancelation ensures that the amplitude of the second term on the right-hand side of (\ref{ec R intro}) is much smaller than its frequence, yielding a small contribution when we invert the divergence.

Regarding the third term, one can deduce from (\ref{expresión phi intro}) that imposing
\begin{equation}
	\abs{B_q\cdot\nabla\theta_{jm}}\ll \eta\ll \lambda
	\label{cancelaciones theta intro 2}
\end{equation}
makes the pertubation $b_{q+1}$ very small. So small, in fact, that the matrix in parenthesis on the third term of (\ref{ec R intro}) is negligible and can be readily absorbed into the new Reynolds stress. Since we only have to worry about the first equation in (\ref{MHD}) and we can essentially ignore the term $B\otimes B$, one could consider that we have reduced convex integration for MHD to convex integration for the Euler equations. The caveat is that the perturbation $w_{q+1}$ is given in terms of the complex expression (\ref{def vq+1 intro}) instead of being freely prescribed, which will be a source of major complications.

As a side remark, note that the necessity to impose both conditions (\ref{cancelaciones theta intro 1}) and (\ref{cancelaciones theta intro 2}) means that our construction has a chance to (and indeed does) work in 3D, but not in 2D. Indeed, we may write $\theta_{jm}(x,t)=\lambda\hspace{0.5pt} K_{jm}\cdot(x,t)$ for some fixed vector $K_{jm}\in \RR^4$ with $\abs{K_{jm}}\sim 1$. Then, (\ref{cancelaciones theta intro 1}) and (\ref{cancelaciones theta intro 2}) can be interpreted as the (almost) orthogonality of $K_{jm}$ with two vectors in $\RR^4$, which leaves the freedom to prescribe $K_{jm}$ to be orthogonal to a third vector $(\zeta_j,0)$. This orthogonality is needed to ensure that $\phi^{q+1}$ is volume-preserving, and flexibility in choosing $\zeta_j$ is necessary when dealing with the last term in (\ref{ec R intro}). Thus, we require three spatial dimensions.

Regarding the last term in (\ref{ec R intro}), as is standard in convex integration, we would like to use the low frequency terms in the product $w_{q+1}\otimes w_{q+1}$ to cancel the previous Reynolds stress. More specifically, we would like to have
\[w_{q+1}\otimes w_{q+1}\approx-[R_q+(p_{q+1}-p_q)\Id]+(\text{high frequency terms}),\]
where the high-frequency terms will be under control due to (\ref{cancelaciones theta intro 1}). To achieve this, we use \cref{geometric lemma} with a suitable choice of $p_{q+1}$ to write the previous Reynolds stress as a sum of simple terms:
\begin{equation}
(p_q-p_{q+1})\Id-R_q=\sum_{j=1}^6\gamma_{qj}^2\hspace{1pt}\zeta_j\otimes\zeta_j,
\label{descomposición geométrica R intro}
\end{equation}
where $\gamma_{qj}\in C^\infty_c(\RR^3\times[0,T])$ are certain amplitudes and $\zeta_j\in \RR^3$ are fixed unitary vectors. With a proper choice of the parameters, a perturbation like (\ref{approx wq+1}) can be used to cancel the previous Reynolds stress $R_q$, up to some small errors.

The diffeomorphism $\phi^{q+1}$ can be engineered to produce the desired perturbation of the velocity to first order, but there will be an extra term $w_\phi$. This term, which arises from the fact that the new velocity is defined through a diffeomorphism, is problematic and genuinely non-standard. Controlling $w_\phi$ is the main challenge of this construction because it is given implicitly by a long, complicated formula involving components that are themselves also not explicit. It is particularly challenging to estimate the term $\mathcal{R}\partial_t w_\phi$, where the operator $\mathcal{R}$ is a right-inverse of the divergence defined in \cref{appendix divergence}. We can only expect good estimates for the operator $\mathcal{R}\partial_t$ when applied to functions that oscillate in space-time in a particular manner. Verifying these estimates for $w_\phi$ is highly nontrivial due to the non-explicit nature of many of its components.

This lack of explicit formulas for many elements in our construction stems from the fact that we cannot construct by hand a volume-preserving diffeomorphism with the appropriate behavior. As a result, $\phi^{q+1}$ is actually obtained as the composition of two  diffeomorphisms: a main diffeomorphism $\phi^0$ that produces the desired perturbations, to first order, and a small correction $\phi^c$ that ensures that $\phi^{q+1}$ is divergence-free. This correction is constructed by means of a fully nonlinear elliptic equation, and it is therefore very hard to obtain any information other than rough estimates on the size of $\phi^c_t-\Id$. Our approach to limiting the effect of this little-known element is to make a very precise choice of $\phi^{0}_t$ so that this diffeomorphism is extremely close to being volume-preserving. As a result, $\phi^c_t$ can be chosen to be very close to the identity in the $C^2$-norm.

An important matter that we have omitted here for simplicity is the support of the fields. We would like $v_q$ and $B_q$ to be compactly supported. Therefore, the Reynolds stress must also be compactly supported, since we will have to perturb the velocity wherever the Reynolds stress is nonzero in order to obtain a weak solution in the limit. As we will see, $R_q$ can only chosen to be compactly supported if the angular momentum of $v_q$ is constant. In principle, this is not guaranteed by our construction, so we have to add a further modification. It is easier to add this perturbation in a region where $B_q$ vanishes. As a result, the subsolution will always be compactly supported, but its support grows with $q$ (albeit in a controllable way).

In summary, constructing the new fields in each step through a diffeomorphism ensures that the second equation of the MHD system (\ref{MHD}) is always satisfied, so we only have to work on the first equation. In addition, by imposing the cancelation (\ref{cancelaciones theta intro 2}) we can essentially ignore the term $B\otimes B$. Hence, we have reduced the problem to convex integration for the Euler equations. The caveat is that the perturbation to the velocity is defined through a diffeomorphism instead of being freely prescribed, which is a source of major complications. Although the diffeomorphism can be engineered to produce the desired perturbation to the velocity to first order, there are extra terms that are complex and not given by an explicit formula. Estimating these terms is a lengthy and challenging task, and obtaining suitable bounds requires a very careful choice of the diffeomorphism.

\subsection{Inductive hypotheses} 
In this subsection, we describe the iterative process in detail. We refer to Appendix~\ref{appendix a} for the definition of the Hölder norms that we use, in which we take into account regularity in space and time.

The perturbation to the velocity and to the magnetic field will oscillate at a frequency controlled by a parameter $\lambda_q$ and the amplitude of the oscillations will be controlled by a parameter $\delta_{q}$, which are given by
\begin{align}
	\lambda_q&\coloneqq a^{b^q}, \label{def lambdaq}\\
	\delta_q&\coloneqq \lambda_q^{-2\beta}, \label{def deltaq}
\end{align}
where $b\coloneqq 1200$ and $a>1$ is a very large parameter that will be chosen later. Meanwhile, $\beta>0$ will determine the Hölder regularity of the final fields, and it is assumed to be sufficiently small. Due to technical reasons that will be clear in \cref{section th main}, we need to introduce this extra parameter $\beta>\alpha$ instead of working directly with the $\alpha$ given in the statement of \cref{main theorem}.

The velocity $v_q$ and the magnetic field $B_q$ will be compactly supported but we must allow for some growth of the support after each step. Let $\bar{r}\geq 1$ such that the supports of $v_0$, $B_0$ and $R_0$ are contained in $B(0,\bar{r}/2)\times[0,T]$.  We define the following sets:
\begin{equation}
\label{def Omegaq}
\Omega_q\coloneqq B(0,(1-2^{-q})\bar{r}\hspace{1pt}).
\end{equation}
We also introduce some intermediate sets that will be useful later:
\begin{equation}
	\label{def Omegaqi}
	\Omega_{q,i}\coloneqq B\hspace{-1pt}\left(0,\left(1-2^{-(q+i/3)}\right)\hspace{-1pt}\bar{r}\right) \hspace{30pt} i\in\{0,1,2\}.
\end{equation}

For any $q\geq 0$, the complete list of hypotheses satisfied by the subsolution $(v_q,B_q,p_q,R_q)$ is as follows: 
\begin{align}
\allowdisplaybreaks
(v_q,B_q,p_q,R_q)&=(0,0,0,0) \qquad \text{outside }\Omega_q\times[0,T], \label{inductive support q}\\
\normxt{0}{B_q}&\leq 1-2^{-q}, \label{inductive tamaño Bq}\\
\normxt{1}{B_q}&\leq \delta_{q}^{1/2}\lambda_q, \label{inductive cota Bq C1}\\
\normxt{0}{v_q}&\leq 1-2^{-q}, \label{inductive tamaño vq} \\
\normxt{N}{v_q}&\leq \delta_{q}^{1/2}\lambda_q^N, \qquad N=1,2,3,4,\label{inductive cotas vq CN}\\
\normxt{0}{R_q}&\leq \delta_{q+1}\lambda_{q}^{-2\tau}, \label{inductive cota Rq C0} \\
\normxt{1}{R_q}&\leq \delta_{q+1}\lambda_q^{1-2\tau}, \label{inductive cota Rq C1} \\
\frac{1}{4}\delta_{q+1}\lambda_q^{-\tau}&\leq e(t)-\int_{\RR^3}\left(\abs{v_q}^2+\abs{B_q}^2\right)dx\leq \frac{3}{4}\delta_{q+1}\lambda_q^{-\tau}, \label{inductive q energía} 
\end{align}
where $\tau>0$ is a very small parameter that will be chosen later.

Given a subsolution satisfying the previous inductive hypotheses, the following proposition shows how to construct the next term in the sequence, which satisfies these hypotheses with $q$ replaced by $q+1$. 
\begin{proposition} \label{prop steps}
There exists $\beta>\frac{1}{200}$ such that given a sufficiently small positive number $\tau>0$ and a number $a>1$ sufficiently large (depending on $\tau$), the following holds:

Let $(v_q,B_q,p_q,R_q)\in C^\infty_c(\RR^3\times[0,T])$ be a subsolution satisfying (\ref{inductive support q})-(\ref{inductive q energía}). Then, there exists a subsolution $(v_{q+1},B_{q+1},p_{q+1},R_{q+1})\in C^\infty_c(\RR^3\times[0,T])$  satisfying (\ref{inductive support q})-(\ref{inductive q energía}) with $q$ replaced by $q+1$. Furthermore, we have the estimates 
\begin{align}
	\normxt{0}{v_{q+1}-v_q}+\lambda_{q+1}^{-1}\normxt{1}{v_{q+1}-v_q}&\leq \delta_{q+1}^{1/2}, \label{cambio v q}\\[3pt] \normxt{0}{B_{q+1}-B_q}+\lambda_{q+1}^{-1}\normxt{1}{B_{q+1}-B_q}&\leq \delta_{q+1}^{1/2}, \label{cambio B q} \\[3pt]
	\abs{\int_{\RR^3}v_{q+1}\cdot B_{q+1}-\int_{\RR^3}v_q\cdot B_q}&\leq \delta_{q+2}. \label{cambio helicidad q}
\end{align}
\end{proposition}
Most of this paper (Sections~\ref{section construction}-\ref{sec correction}) is devoted to the proof of this proposition, which is the core of our construction. In \cref{section construction} we construct the new subsolution $(v_{q+1},B_{q+1},p_{q+1},R_{q+1})$, which is obtained from the previous one through a volume-preserving diffeomorphism $\phi^{q+1}$. In \cref{section estimates diffeo} we derive suitable bounds on $\phi^{q+1}$. These are necessary for estimating the new velocity $v_{q+1}$ and the new magnetic field $B_{q+1}$, which is done in \cref{section estimates vector fields}. In \cref{section matrix} we check that the perturbations to both vector fields produce the desired correction to the Reynolds stress. In \cref{sec correction} we introduce a further correction that ensures that we can fully control the support of the new subsolution. This will complete the proof of \cref{prop steps}. Finally, in \cref{section th main} we use \cref{prop steps} with a suitable initial subsolution $(v_0,B_0,p_0,R_0)$ to construct a sequence that converges to the desired weak solution, thereby proving \cref{main theorem}.

\section{Construction of the new subsolution}\label{section construction}
In this section we construct the new subsolution $(v_{q+1},B_{q+1},p_{q+1},R_{q+1})$, which is the first step in the proof of \cref{prop steps}. The new subsolution is obtained from the previous one through a volume-preserving diffeomorphism $\phi^{q+1}$. This section is devoted to the definition of these objects. We begin with some preliminaries, which include some notation that we will use and certain parameters that are important in our construction. In \cref{subsec simplification} we write the Reynolds stress $R_q$ in a more suitable form. In \cref{subsec def diffeos} we construct the diffeomorphism $\phi^{q+1}$, which allows us to define the new velocity and magnetic field in \cref{subsec def vectors}. Having defined the new fields, in \cref{subsec def Reynolds} we check how much they fail to solve \eqref{MHD}, that is, we construct the new Reynolds stress. Finally, in \cref{subsec final correction} we introduce a further correction that allows us to keep control on the support of the new Reynolds stress, i.e., \eqref{inductive support q}.

\subsection{Preliminaries}\label{section preliminaries}
There is a series of parameters that describes the size of the elements of our construction. The diffeomorphism $\phi^0$ and the perturbations to the velocity and to the magnetic field will oscillate with frequency $\lambda_{q+1}$. The deviation of $\phi^0$ from the identity (more precisely, from $\pi$) will be measured by a parameter $\eta$. Finally, we introduce a parameter $\mu^{-1}$ that will control the size of certain cutoffs that we will use. The hierarchy of these parameters is as follows:
\begin{equation}
	\label{hierarchy}
	\lambda_{q}\ll\mu\ll \eta \ll \lambda_{q+1}
\end{equation}
They are given by
\begin{align}
	\mu&\coloneqq\frac{\delta_{q+2}}{\delta_{q+1}^{1/2}}\lambda_{q+1}^{1/3-12\tau}, \label{def mu} \\ \eta&\coloneqq \frac{\delta_{q+2}}{\delta_{q+1}^{1/2}}\lambda_{q+1}^{1-12\tau}. \label{def eta}
\end{align}
As we will see, these definitions guarantee certain relationships between the parameters that lead to suitable bounds for the new subsolution. Except for $\lambda_q\ll \mu$, which will be checked at the beginning of \cref{section estimates diffeo}, the rest of hierarchy (\ref{hierarchy}) is clear in view of these definitions.

Next, we fix some notation that we will use. First of all, we refer again to Appendix~\ref{appendix a} for the definition of the norms that we will use. We warn the reader that our notation differs from the standard notation in convex integration: we use space-time Hölder norms, as opposed to the supremum in time of the spatial Hölder norms.

We will work with maps defined in $\RR^3\times[0,T]$ that have particular properties when we fix a time $t\in [0,T]$. For instance, $\phi^{q+1}_t:\RR^3\to\RR^3$ will be a diffeomorphism. To avoid writing too many subscripts and to simplify our exposition, we make the convention to just say that $\phi^{q+1}$ is a diffeomorphism, meaning that it is a diffeomorphism of $\RR^3$ for each fixed time $t\in [0,T]$. In addition, it will be convenient to introduce the  definitions
\begin{align}
	\pi(x,t)&\coloneqq x, \label{def pi} \\
	(f\compt g)(x,t)&\coloneqq f(g(x,t),t) \label{def circt}
\end{align}
for the spatial projection and for ``spatial composition'' of maps $f:\RR^3\times[0,T]\to\RR^3$ and $g:\RR^3\times[0,T]\to\RR^3$, for any $m\geq1$.
Note that the operation $\compt$ is associative, and that
\begin{equation}
	f\compt \pi=f,\qquad \pi\compt g=g \label{f pi}
\end{equation}
with $f,g$ as above. For each $t\in[0,T]$, we will also set
\[
f_t:=f(\cdot,t).
\]

We will write $A\lesssim B$ to denote $A\leq C\hspace{0.5pt}B$ for some constant $C>0$ that is independent of the iteration index $q$ or the parameter $a$. However, the implicit constants are allowed to depend on $\beta$ and $\tau$. In addition, when estimating a $C^N$-norm, they may also depend on $N$. This is not an issue because we will only use a finite number of derivatives.

The goal of this notation is to focus on the dependence on the iteration parameters $\delta_q$, $\lambda_q$, $\mu$ and $\eta$. The implicit constants will be irrelevant in most inequalities because we will have an extra factor that can be made arbitrarily small by increasing $a$. Therefore, once we have chosen suitable values of $\beta$ and $\tau$, choosing $a$ sufficiently large will compensate for all the implicit constants, provided that they do not depend on $a$ and that there is a finite number of them. This is ensured by excluding dependence on $q$ and by using only a finite number of derivatives.

\subsection{Simplification of the Reynolds stress}\label{subsec simplification}
The first step is to write the Rey\-nolds stress $R_q$ in a more suitable form. We want to express it as a sum of rank-1 matrices that have disjoint support. These simpler matrices are precisely of the form that can be corrected with our perturbation to the velocity. Their having disjoint support allows us to correct all of these components at the same time without influencing one another. Otherwise, we would need to correct each component in a different step, due to technical reasons that will be discussed in \cref{remark problemas si soporte no disjunto}. This would dramatically decrease the regularity of the final solution because the frequencies would have to grow much faster.

The separation of the supports is essentially a separation in space using the Mikado profiles introduced in \cite{Daneri17}. Unlike in \cite{BLSV,Isett2018}, here Mikados are used to control the amplitude of the perturbations, not as the oscillation itself. In~\cite{Giri} the authors introduce a Newton iteration that produces decoupling in time. This allowed them to prove Onsager's conjecture in 2D and has been highly successful in other contexts~\cite{DGR}. However, adapting this technique to our problem would be highly nontrivial, especially since we want (\ref{MHD b}) to hold exactly. Nevertheless, at the currently achievable regularity it is not necessary, as our decoupling in space is more than enough.

First of all, by absorbing $-\frac{1}{3}\tr{R_q}$ into the pressure, we may assume
\begin{equation}
\tr R_q=0.
\label{traza de Rq}
\end{equation}
This can be done while preserving \eqref{inductive support q}, \eqref{inductive cota Rq C0} and \eqref{inductive cota Rq C1}. Next, we recall a geometric lemma that is standard in convex integration. Its proof can be found in \cref{app.C}.
\begin{lemma}\label{geometric lemma}
	There exist $c_0>0$, unitary vectors $\zeta_1, \dots, \zeta_6\in \RR^3$ and smooth maps $\Gamma_1, \dots, \Gamma_6$ defined in the ball $\overline{B}(\Id,c_0)\subset\RR^{3\times 3}_{\text{sym}}$ such that
	\[S=\sum_{j=1}^6\Gamma_j(S)^2\zeta_j\otimes\zeta_j \qquad \forall S\in \overline{B}(\Id,c_0).\]
\end{lemma}
We now introduce a time-dependent amplitude $\rho_q$ that will be used when applying the previous lemma. Its role is to artificially increase the size of the coefficients of the decomposition, which will allow us to prescribe the energy. First, we fix a smooth cutoff $\overline{\sigma}_q\in C^\infty_c(\Omega_{q,1},[0,1])$ that is identically 1 on $\Omega_q$. This cutoff will be used to restrict the perturbations to a neighborhood of the support of the previous subsolution.
\begin{lemma}
	\label{lema rhoq}
	If $a>1$ is sufficiently large, the function $\rho_q\in C^\infty([0,T])$ given by
	\begin{equation}
		\rho_q(t)\coloneqq\frac{1}{3}\norm{\overline{\sigma}_q}_{L^2(\RR^3)}^{-2}\left[e(t)-\frac{1}{2}\delta_{q+2}\lambda_{q+1}^{-\tau}-\int_{\TT^3}\left(\abs{v_q}^2+\abs{B_q}^2\right)dx\right]
		\label{def rhoq}
	\end{equation}
	satisfies
	\begin{equation*}
		\delta_{q+1}\lambda_q^{-\tau}\lesssim \norm{\rho_q}_{C^0_t}\lesssim\delta_{q+1}\lambda_q^{-\tau}, \hspace{50pt}\norm{\rho_q}_{C^1_t}\lesssim\delta_{q+1}\lambda_q^{1-\tau}.
	\end{equation*}
\end{lemma}
\begin{proof}
	First of all, by construction of $\overline{\sigma}_q$ and the sets $\Omega_q$ and $\Omega_{q,1}$, we see that $\overline{\sigma}_q$ takes values in $[0,1]$, is supported in $B(0,\bar{r})$ and equals 1 on $B(0,\bar{r}/2)$. Therefore, its $L^2$ norm is bounded above and below by a constant independent of $q$.	Hence, it follows from \eqref{inductive q energía} that
	\[\abs{\rho(t)}\geq \frac{1}{3}\norm{\overline{\sigma}_q}_{L^2(\RR^3)}^{-2}\left(\frac{1}{4}\delta_{q+1}\lambda_q^{-\tau}-\frac{1}{2}\delta_{q+2}\lambda_{q+1}^{-\tau}\right)\geq \frac{1}{15}\norm{\overline{\sigma}_q}_{L^2(\RR^3)}^{-2}\delta_{q+1}\lambda_{q}^{-\tau}\gtrsim\delta_{q+1}\lambda_{q}^{-\tau},\]
	where we have taken $\tau>0$ sufficiently small and $a>1$ large.  Similarly, by \eqref{inductive q energía}:
	\[\abs{\rho(t)}\leq \frac{1}{4}\norm{\overline{\sigma}_q}_{L^2(\RR^3)}^{-2}\delta_{q+1}\lambda_q^{-\tau}\lesssim \delta_{q+1}\lambda_q^{-\tau}.\]
	Regarding the derivative, it follows from \eqref{MHD} and standard vector calculus identities that
	\[\frac{d}{dt}\int_{\TT^3}\left(\abs{v_q}^2+\abs{B_q}^2\right)dx=2\int_{\TT^3}(v_q\cdot\partial_t v_q+B_q\cdot\partial_tB_q)\,dx=2\int_{\TT^3}v_q\cdot\Div R_q\,dx.\]
	Using \eqref{inductive cota Rq C1}, we conclude that the time derivative of the third term in (\ref{def rhoq}) is bounded (up to a constant) by $\delta_{q+1}\lambda_q^{1-2\tau}$. Since the size of $\partial_t e$ will be negligible with respect to this quantity for $a>1$ sufficiently large, the claimed bound for the $C^1$ norm of $\rho_q$ follows.
\end{proof}
It follows from the previous lemma and \eqref{inductive cota Rq C0} that $$\normxt{0}{\rho_q^{-1}R_q}\lesssim \lambda_q^{-\tau},$$ which for large $a>1$ will be smaller than the radius $c_0>0$ that appears in \cref{geometric lemma}. Therefore, we can write:
\begin{equation}
	\label{descomposición geométrica R_q}
	\overline{\sigma}_q^2\rho_q\Id-R_q=\sum_{j=1}^6\overline{\sigma}_q^2\rho_q\,\Gamma_j\hspace{-2pt}\left(\Id-\rho_q^{-1}R_q\right)^2\zeta_j\otimes\zeta_j\equiv \sum_{j=1}^6\gamma_{qj}^2\hspace{1pt}\zeta_j\otimes\zeta_j.
\end{equation}
where we have used that the cutoff $\overline{\sigma}_q$ is identically 1 on $\Omega_q$ and, thus, on the support of $R_q$.
Note that we only have control on the $C^0$ and $C^1$ norms of these coefficients because only \eqref{inductive cota Rq C0} an \eqref{inductive cota Rq C1} are available. However, in our construction we will need control on all $C^N$ norms, so the next step is to mollify these coefficients in space and time: 
\begin{lemma}
	\label{lema aprox gammas}
	For any $j\in \{1,\dots, 6\}$ there exists $\overline{\gamma}_{qj}\in C^\infty_c(\Omega_{q,2}\times[0,T])$ such that
	\begin{align}
		\normxt{0}{\gamma_{qj}^2-\overline{\gamma}_{qj}^2}&\lesssim \delta_{q+2}\lambda_{q+1}^{-11\tau}, \label{diferencia gamatilde} \\
		\normxt{0}{\overline{\gamma}_{qj}}&\lesssim\delta_{q+1}^{1/2}\lambda_q^{-\tau/2}, \label{cota gamatilde C0}\\
		\normxt{N}{\overline{\gamma}_{qj}}&\lesssim \delta_{q+1}^{1/2}\left(\delta_{q+1}\delta_{q+2}^{-1}\lambda_q\lambda_{q+1}^{11\tau}\right)^N \leq  \delta_{q+1}^{1/2}\mu^{N} \hspace{30pt} \forall N\geq 0. \label{cotas gamatilde}
	\end{align}
\end{lemma}
\begin{proof}
Since $2^q\lambda_q^{-1}$ can be made arbitrarily small by taking $a>1$ sufficiently large (independent of $q$), the cutoff $\overline{\sigma}_q$ can be chosen to satisfy $\normxt{1}{\overline{\sigma}_q}\lesssim \lambda_q$. See \cref{lemma cutoff}. Then, it follows from \eqref{inductive cota Rq C1} and \cref{lema rhoq} that
\begin{equation}
\normxt{0}{\gamma_{qj}}+\lambda_q^{-1}\normxt{1}{\gamma_{qj}}\lesssim \delta_{q+1}^{1/2}\lambda_q^{-\tau/2}.
\label{aux cotas gammaqj}
\end{equation}
We consider a smooth extension of $\gamma_{qj}$ to the whole $\RR^3\times\RR$, which we still denote as $\gamma_{qj}$. This extension may be assumed to satisfy \eqref{aux cotas gammaqj}. Next, we fix a standard convolution kernel $\varphi\in C^\infty_c(\RR^4)$ with $\int \varphi=1$ and whose support is contained in the unit ball. We set $\varphi_{\bar{\varepsilon}}\coloneqq \bar{\varepsilon}^{\hspace{1pt}4}\varphi(\bar{\varepsilon} \,\cdot\,)$ with $\bar{\varepsilon}\coloneqq (\delta_{q+1}\lambda_q)^{-1}\delta_{q+2}\lambda_{q+1}^{-11\tau}$ and we define:
\[\overline{\gamma}_{qj}\coloneqq \gamma_{qj}\ast \varphi_{\bar{\varepsilon}},\]
where the convolution is in space and time. The bounds (\ref{aux cotas gammaqj}) and standard estimates for mollifiers yield \eqref{diferencia gamatilde}, \eqref{cota gamatilde C0} and the first inequality in \eqref{cotas gamatilde}. We are particularly interested in keeping the factor $\lambda_q^{-\tau/2}$ in the bound for the $C^0$ norm because it will help us get rid of the implicit constants later on. On the other hand, in \eqref{relación entre lambdas} we will prove that $\bar{\varepsilon}^{\hspace{1pt}-1}$ is much smaller than $\mu$, which yields the second inequality in \eqref{cotas gamatilde}. Finally, since $2^q\bar{\varepsilon}\leq 2^q\lambda_q^{-1}$, we may assume that the support of $\overline{\gamma}_{qj}(\cdot,t)$ is contained in $\Omega_{q,2}$ for all times $t\in[0,T]$.
\end{proof}

Having decomposed the Reynolds stress as a sum of simpler matrices, now our goal is to make the different components have disjoint support. We will use the Mikado profiles introduced in~\cite{Daneri17}. Since a rescaling of our vectors $\zeta_1, \dots, \zeta_6\in \RR^3$ is in $\ZZ^3$, the proof is the same as in~\cite{Daneri17}.
\begin{lemma}\label{lemma mikado}
	Let $\zeta_1, \dots, \zeta_6\in \RR^3$ be as in the previous lemma. There exist $c_1>0$ and smooth functions $\sigma_1,\dots, \sigma_6\in C^\infty(\TT^3)$ such that the following holds for $j\in\{1, \dots, 6\}$:
	
	\vspace{-2pt}
	\begin{enumerate}
		\setlength\itemsep{3pt}
		\item $\zeta_j\cdot \nabla\sigma_j=0$,
		\item $\dist(\supp\sigma_j\hspace{0.5pt},\, \supp \sigma_{j'})>2c_1\;\;$ for $j\neq j'$,
		\item $\fint_{\TT^3}\sigma_j^2(x,t)\hspace{0.5pt}dx=1$.
	\end{enumerate}
\end{lemma}

Since $\zeta_j\cdot\nabla\sigma_j=0$, taking the divergence of (\ref{descomposición geométrica R_q}) leads to
\begin{equation}
	\label{descomposición final R_q}
	\begin{aligned}
		\nabla(\overline{\sigma}_q^2\rho_q)-\Div R_q&=\Div\hspace{-2pt}\left(\sum_{j=1}^6[\gamma_{qj}^2-\overline{\gamma}_{qj}^2]\,\zeta_j\otimes \zeta_j\right)+\Div\hspace{-2pt}\left(\sum_{j=1}^6\sigma_j(\mu \,\cdot\,)^2\hspace{1.5pt}\overline{\gamma}_{qj}^2\hspace{0.5pt}\zeta_j\otimes \zeta_j\right)\\&\hspace{12pt}+\sum_{j=1}^6[1-\sigma_j(\mu \,\cdot\,)^2]\Div\hspace{-1pt}\left(\overline{\gamma}_{qj}^2\hspace{0.5pt}\zeta_j\otimes \zeta_j\right).\end{aligned}
\end{equation}
In \cref{section matrix} we will see that the first matrix in parenthesis is sufficiently small to be absorbed into the new Reynolds stress $R_{q+1}$. Similarly, we will see that we can find a suitably small matrix whose divergence equals the third term. The key for this is the fact that the mean of $\sigma_j^2$ is 1, so $1-\sigma_j(\mu\,\cdot\,)^2$ is a high-frequency term. Therefore, we only have to correct the matrix in the second parenthesis, which is a sum of rank-1 matrices with disjoint support, as we wanted.  

\subsection{Definition of the diffeomorphism}\label{subsec def diffeos}
The diffeomorphism $\phi^{q+1}$ will be of the form
\[\phi^{q+1}=\phi^c\compt\phi^0,\]
where $\phi^0$ is the main correction term, which will yield the main perturbation $w_0$ to the velocity, while $\phi^c$ is very close to the identity and its only role is to ensure that $\phi^{q+1}$ is volume preserving. This guarantees that $v_{q+1}$ and $B_{q+1}$ are divergence-free.

In practice, we do not define the diffeomorphism $\phi^0$ directly. Instead, it is more convenient to define its inverse. To this end, we construct an auxiliary map $\psi^0$, show that it is invertible and define $\phi^0$ as its inverse, i.e., $\phi^0_t:=(\psi^0_t)^{-1}$. This approach is motivated by the fact that obtaining an explicit formula for the inverse of a diffeomorphism is rarely possible; typically, we can only write down one of the two explicitly. We have chosen to define $\psi^0$ explicitly because it simplifies certain estimates. Nevertheless, both $\phi^0$ and its inverse appear in the definition of the pushforward, so we will need good control on both maps.

First, we define certain cutoffs that we will need. We fix a smooth cutoff function $\chi\in C^\infty_c\left(\left(-\frac{3}{4},\frac{3}{4}\right)^4\right)$ such that 
\[\sum_{m\in \ZZ^4}\chi(y-m)^2=1\]
and for $m\in \ZZ^4$ we define 
\begin{equation}
	\chi_m\coloneqq \chi(\mu \,\cdot\,-\hspace{1.5pt}m).
	\label{def chim}
\end{equation}
Since our fields are compactly supported, we only need a finite number of these cutoffs. We define
\[\Lambda\coloneqq \left\{m\in\ZZ^4:\,\supp \chi_m\cap \big([B(0,\overline{r})\times[0,T]\big)\neq \varnothing\right\},\]
where $\overline{r}>0$ is given by \eqref{def Omegaq}. By \eqref{inductive support q}, the support of all the objects that we will consider will be contained in $B(0,\overline{r})\times[0,T]$, so we do not need to consider values $m\notin\Lambda$. Note that the number of elements in $\Lambda$ grows like $\mu^4$, so it will be quite large for large values of $q$. However, any point $(x,t)\in \RR^3\times[0,T]$ is in the support of at most 16 of the cutoffs. This will be very important because it will ensure that our estimates are independent of the number of elements in $\Lambda$.

We define the amplitude of the perturbations as
\begin{equation}
	\label{def am}
	a_{jm}\coloneqq \sqrt{2}\,\chi_m\,\sigma_j(\mu\,\cdot\,)\,\overline{\gamma}_{qj}.
\end{equation}
Hence, up to a numerical constant, the amplitude is given by $\overline{\gamma}_{qj}$, but we localize to a small set by means of cutoffs. Indeed, the purpose of $\chi_m$ is to restrict to a very small region in which $v_q$ and $B_q$ can be considered to be approximately constant. For each $m\in \Lambda$ we define the constant vectors
\begin{align}
	v^q_{jm}&\coloneqq v_q(\mu^{-1}m), \label{def vJm}\\ \qquad B^q_{jm}&\coloneqq B_q(\mu^{-1}m).
\end{align}
We then choose a unitary vector $k_{jm}\in \mathbb{S}^2$ such that
\[k_{jm}\,\bot\, \zeta_j\qquad \text{and}\qquad k_{jm}\,\bot\, B^q_{jm},\]
where $\zeta_j$ are given by \cref{geometric lemma}. This vector determines the direction of oscillation. 

As a final ingredient, we need to introduce certain coefficients that will prevent unwanted interference later on. We define
\begin{equation}\label{eq.L}
\mathcal{L}\coloneqq\left\{\sqrt{p_1}\hspace{1pt},\hspace{1pt} \dots\hspace{1pt},\hspace{1pt} \sqrt{p_{16}}\hspace{1pt}\right\},
\end{equation}
where $p_n$ is the $n$-th prime. Since the elements of $\mathcal{L}$ are linearly independent over $\mathbb{Q}$, any nontrivial linear combination of elements of $\mathcal{L}$ with integer coefficients must be nonzero. It is easy to see that for each $m\in\Lambda$ we can choose $\ell_m\in \mathcal{L}$ such that
\[\supp\chi_m\cap\supp\chi_{m'}\neq \varnothing, \quad m\neq m'\qquad \Rightarrow\qquad \ell_m\neq \ell_{m'}.\]
Then, we define the phase
\begin{equation}
	\theta_{jm}(x,t)\coloneqq \ell_m\lambda_{q+1}k_{jm}\cdot(x-v^q_{jm}t)+\ell_m\eta t.
	\label{def thetam}
\end{equation}
We emphasize some properties of these phases that will be very important in our construction:
\begin{align}
	&\abs{\nabla\theta_{jm}}=\ell_m\lambda_{q+1}, \\
	& \zeta_j\cdot\nabla \theta_{jm}=0, \\
	&B^q_{jm}\cdot\nabla \theta_{jm}=0, \label{cancelation B}\\
	&\partial_t\theta_{jm}+v^q_{jm}\cdot\nabla \theta_{jm}=\ell_m\eta. \label{derivada material thetam}
\end{align}

After all of these previous definitions, we are finally ready to construct the auxiliary diffeomorphism $\psi^0\in C^\infty(\RR^3\times[0,T],\RR^3)$:
\begin{equation}
\psi^0\coloneqq \pi- \sum_{j=1}^6\sum_{m\in \Lambda}\frac{a_{jm}}{\ell_m\eta}\,\zeta_j\hspace{0.5pt}\sin(\theta_{jm})+\sum_{j=1}^6\sum_{m\in \Lambda}\frac{\nabla a_{jm}\times(\zeta_j\times k_{jm})}{\ell_m^2\eta\lambda_{q+1}}\hspace{0.5pt}\cos(\theta_{jm}).
\label{def psi}
\end{equation}
Since the support of the coefficients $a_{jm}$ is contained in $\Omega_{q,2}\times [0,T]$ because so is the support of $\overline{\gamma}_{qj}$, we see that $\psi^0\equiv \pi$ in a neighborhood of $(\RR^3\backslash\Omega_{q,2})\times[0,T]$.

\begin{remark}
\label{remark soporte ajm disjuntos}
For a fixed value of $j\in\{1,\dots, 6\}$, the support of each $a_{jm}$ intersects the support of at most 15 other coefficients $a_{jm'}$, by construction of the cutoffs $\chi_m$. However, the support of $a_{jm}$ does not intersect the support of a coefficient with a different value of $j$ because the supports of $\sigma_j(\mu\,\cdot\,)$ are disjoint for different values of $j$. This will prevent undesired terms in the next Reynolds stress $R_{q+1}$.
\end{remark}

\begin{remark}
	\label{fase invariante}
	Since $k_{jm}$ is orthogonal to $\zeta_j$ for any $m\in \Lambda$, we see that
	\[\sup_{(x,t)\in\supp a_{jm}}\abs{\theta_{jm}\compt\psi^0-\theta_{jm}}\lesssim \frac{1}{\eta}\max_{n\in \Lambda}
	\normxt{1}{a_{jn}}
	.\] 
	Note that we only have to consider the terms with the same value of $j$, by the previous remark. We will see in \cref{lemma estimates am} that the right-hand side is much smaller than 1. Thus, the effect of the diffeomorphism on the phases $\theta_{jm}$ is very small. This is essential for obtaining a suitable formula for the inverse of $\psi^0$, which is necessary to identify the main correction $w_0$ to the velocity and to derive the appropriate estimates. 
	\end{remark}
	\begin{remark}
	\label{remark problemas si soporte no disjunto}
	If the supports of the terms associated to different values of $j$ were not disjoint, we would have no orthogonality between $k_{jm}$ and $\zeta_{j'}$, yielding no control on the effect of the diffeomorphism on the phases. The final consequence would be that we would have no information on the form of the perturbations, so the construction would break down. The alternative would be to correct each component $\zeta_j\otimes \zeta_j$ in different small steps, which would dramatically decrease the regularity of the final solution.
\end{remark}

It is not clear whether $\psi^0$ is a diffeomorphism (for each $t\in[0,T]$), as we claim. However, we will see that its Jacobian is very close to 1. This is quite remarkable, because $D\psi^0$ is quite large. The reason behind this is the orthogonality that we are imposing. The Jacobian being nonzero means that $\psi^0$ is a local diffeomorphism. To prove that it is a global diffeomorphism, we will apply Banach's fixed point theorem to the map
\begin{align}
\begin{split}
	\label{def T contractiva}
	T:C^0(\RR^3\times[0,T],\RR^3)\,&\to\; C^0(\RR^3\times[0,T],\RR^3) \\  f\hspace{25pt}\,&\mapsto\; \pi+\sum_{j=1}^6\sum_{m\in\Lambda}\frac{a_{jm}\compt f }{\ell_m\eta}\,\zeta_j\sin(\theta_{jm}\compt f)\\&\hspace{24pt}-\sum_{j=1}^6\sum_{m\in\Lambda}\frac{(\nabla a_{jm}\compt f)\times(\zeta_j\times k_{jm})}{\ell_m^2\eta\lambda_{q+1}}\hspace{0.5pt}\cos(\theta_{jm}\compt f).
\end{split}
\end{align}
We will see that the fixed point $\phi^0\in C^0(\RR^3\times[0,T],\RR^3)$ satisfies $\phi^0=(\psi^0)^{-1}$. Since the local inverse given by the inverse function theorem is smooth, so is $\phi^0$. It is not possible to obtain an explicit formula for $\phi^0$, but the identity $T\phi^0=\phi^0$ will be enough to derive all the estimates that we will need.

We will see that $\phi^0$ is not volume-preserving, so we need to introduce a correction $\phi^c$ in order to obtain the final diffeomorphism $\phi^{q+1}$. Let 
\begin{equation}
	\label{def fc}
	f_c\coloneqq \det(D\psi^0)-1.
\end{equation}
It is easy to check that
\[\det(D \phi^c)=1+f_c \quad \Rightarrow \quad \det[D(\phi^c\circ\phi^0)]=1.\]
Note that if we tried to express $f_c$ in terms of $\phi^0$ instead of $\psi^0$, we would obtain a more complicated expression. This is another advantage of having an explicit formula for the inverse of $\phi^0$ instead of a formula for $\phi^0$ itself.

As we will see, $D\psi^0$ can be quite large. However, its determinant will turn out to very close to 1, so $f_c$ will be very small. As a result, the correction $\phi^c$ can be chosen to be very close to the identity. This is crucial in our construction: since we have little information on $\phi^c$, the only way to ensure that it does not have unwanted effects is to make certain that it is sufficiently close to the identity. For this, a very precise choice of $\psi^0$ is essential.

Since $\psi^0$ is the identity in a neighborhood of $\RR^3\backslash\Omega_{q,2}$ for all $t\in[0,T]$, the support of $f_c$ is contained in $\Omega_{q,2}\times[0,T]$. Hence, we do not have to correct the diffeomorphism in the whole space: we will be able to choose the correction $\phi^c$ so that $\phi^c$ is the identity on $\RR^3\backslash\Omega_{q,2}$. This is desirable because we only wish to perturb the subsolution in this set, in view of (\ref{inductive support q}).  

To construct a diffeomorphism $\phi^c$ such that $\det(D \phi^c)=1+f_c$, we will use the following lemma, which is just a simplified version of the results in \cite{DM}: 
\begin{proposition}
	\label{lema dacorogna}
	Let $\Omega\subset \RR^3$ be a bounded domain with smooth boundary. Let $f\in C^\infty_c(\Omega\times[0,T],\RR)$ such that
	$\int_{\Omega}f=0$ for all $t\in[0,T]$. Let $u\in C^\infty_c(\Omega\times[0,T],\RR^3)$ such that
	\[\Div u=f.\]
	In addition, let us assume that
	\begin{align}
		\normxt{0}{f}&\leq 1/2, \label{assumption dacorogna 1} \\
		\normxt{1}{u}+\normxt{0}{u}\normxt{1}{f} &\leq 1. \label{assumption dacorogna 2}
	\end{align}
	Then, there exists a map $\phi\in C^\infty(\RR^3\times[0,T],\RR^3)$ such that $\phi_t:\RR^3\to\RR^3$ is a diffeomorphism for all $t\in[0,T]$ and \[\det(D \phi_t)=1+f.\] Furthermore, $\phi\equiv \pi$ in $(\RR^3\backslash\Omega)\times[0,T]$ and it satisfies
	\begin{align}
		\normxt{0}{\phi-\pi}&\leq C(\Omega)\normxt{0}{u}, \label{estimate phi dacorogna C0}\\
		\normxt{N}{\phi-\pi}&\leq C(\Omega,N)\left( \normxt{N}{u}+\normxt{0}{u}\normxt{N}{f}\right) \qquad \forall N\geq 1  \label{estimate phi dacorogna CN}.
	\end{align}
\end{proposition}
Although the construction in the previous lemma follows a standard approach, we postpone the proof until \cref{section estimates diffeo}, where we establish the necessary bounds for the diffeomorphism.

Applying the previous lemma to the function $f_c$, the set $\Omega_{q,2}$ and a suitable vector field $u_c$, we will obtain a diffeomorphism $\phi^c$ such that $\det(D\phi^c)=1+f_c$ and  $\phi^c\equiv\pi$ on $(\RR^3\backslash\Omega_{q,2})\times[0,T]$. We define the final map as $\phi^{q+1}\coloneqq \phi^c\compt \phi^0$. Hence, for each fixed time $t\in[0,T]$, we have a volume-preserving diffeomorphism $\phi^{q+1}_t= \phi^c_t\circ \phi^0_t$ of $\RR^3$ that, by construction, is the identity outside of $\Omega_{q,2}$.

The correction $\phi^c$ will be crucial to ensure that the velocity and the magnetic field remain divergence-free, yet its influence on the Reynolds stress will be negligible. Consequently, the primary contribution to the velocity perturbation will come from $\phi^0$, which is specifically designed to introduce the required correction. Meanwhile, the perturbation to the magnetic field due to the action of $\phi^{q+1}$ will be so small that it can be essentially ignored.

\subsection{The new velocity and magnetic field}\label{subsec def vectors}
Given the flow $X^q$ of the previous velocity $v_q$, we define a new flow $X^{q+1}\in C^\infty(\RR^3\times[0,T],\RR^3)$ as follows:
\[X^{q+1}\coloneqq \phi^{q+1}\compt X^q\compt\left(\phi^{q+1}_0\right)^{-1}.\]
 Since $\phi^{q+1}$ is volume-preserving by construction and $X^q$ is volume-preserving due to being generated by a divergence-free field, we conclude that $X^{q+1}$ is volume-preserving for all times $t\in[0,T]$, as required. 

Next, we define the new magnetic field as
\[B_{q+1}\coloneqq (\phi^{q+1}_t)_\ast\, B_{q}.\]
Arguing as in (\ref{pushforward B intro}), one can check that this definition implies
\[B_{q+1}(\cdot,t)=(X^{q+1}_t)_\ast\,B_{q+1}(\cdot,0).\]
Taking into account that $v_{q+1}$ will be defined so that it generates $X^{q+1}$, we see that $B_{q+1}$ will satisfy (\ref{subs c}). Furthermore, since $\phi^{q+1}_t$ is volume-preserving and $B_q$ is divergence-free, so is $B_{q+1}$. In addition, $B_{q+1}=B_q=0$ on $(\RR^3\backslash\Omega_{q,2})\times[0,T]$ because $\phi^{q+1}=\pi$ on that set. On the other hand, we will see that the perturbation 
\[b_{q+1}\coloneqq B_{q+1}-B_q\]
is very small, due to the cancellation (\ref{cancelation B}). In fact, it is so small that it can be incorporated into $R_{q+1}$. As a result, we do not need a detailed study of $b_{q+1}$, just to estimate its size.

The new velocity field is defined as:
\begin{equation}
v_{q+1}\coloneqq [\partial_t \phi^{q+1}_t+D\phi^{q+1}_t(v_q)]\circ (\phi^{q+1}_t)^{-1},
\label{def vq+1}
\end{equation}
which ensures that the map $X^{q+1}$ is the flow of $v_{q+1}$. Indeed, the same reasoning as in (\ref{vq+1 genera Xq+1}) applies. Since $X^{q+1}$ is volume-preserving, $v_{q+1}$ must be divergence-free. In addition, $v_{q+1}=v_q=0$ on $(\RR^3\backslash\Omega_{q,2})\times[0,T]$ because $\phi^{q+1}=\pi$ on that set. 

It will be convenient to decompose the perturbation $w_{q+1}\coloneqq v_{q+1}-v_q$ as a sum $w_{q+1}=w_0+w_c+w_\phi$, where $w_0$ is the main correction term, while the other two terms are much smaller. We define
\begin{align}
	w_0(x,t)&\coloneqq \sum_{j=1}^6\sum_{m\in \Lambda}a_{jm}(x,t)\,\zeta_j\hspace{0.5pt}\cos[\theta_{jm}(x,t)],
	\label{def w0} \\ w_c(x,t)&\coloneqq \sum_{j=1}^6\sum_{m\in \Lambda} \frac{\nabla a_{jm}(x,t)\times(\zeta_j\times k_{jm})}{\ell_m\lambda_{q+1}}\sin[\theta_{jm}(x,t)] \label{def wc} \\
	w_\phi&\coloneqq v_{q+1}-v_q-w_0-w_c. 
\end{align}
We will see later on that the perturbation $w_0$ provides the desired correction to $R_J$, while $w_c$ and $w_\phi$ are negligible with regards to the Reynolds stress.

The term $w_c$ can be viewed as the small correction that is usually present in convex integration schemes in order to ensure that the perturbation is divergence-free. Indeed, one can check that
\begin{equation}
w_0+w_c=\curl\left(\sum_{j=1}^6\sum_{m\in \Lambda}\frac{a_{jm}}{\ell_m\lambda_{q+1}}(\zeta_j\times k_{jm})\sin(\theta_{jm})\right).
\label{fórmula w0+wc}
\end{equation}
We have chosen to isolate this term mostly for convenience; however,  the remaining term $w_\phi$ is just as large, if not larger, than $w_c$. 

In a usual convex integration scheme, one would simply define the perturbation as (\ref{fórmula w0+wc}). Nevertheless, in our setting we must make sure that (\ref{MHD b}) is satisfied, so we must define the new velocity as (\ref{def vq+1}). This leads to a new perturbation term $w_\phi$.

Therefore, we can think of $w_\phi$ as an artifact caused by defining the new velocity via diffeomorphisms, which is a unique (and inconvenient) feature of our construction. As we will see, $w_\phi$ has a very complex expression, which makes handling this term very difficult. In fact, estimating $w_\phi$ and its effects on the Reynolds stress is the most delicate part of this paper.

As a final thought note that, although $X^{q+1}$ is conceptually very important, in reality we will never work with it. Indeed, $v_{q+1}$ is defined as the generator of $X^{q+1}$, while the evolution of $B_{q+1}$ is determined by $X^{q+1}$ through (\ref{subs c}), which is the key to ensure that (\ref{MHD b}) holds. However, in practice both fields are constructed merely in terms of the previous fields and the map $\phi^{q+1}$.

\subsection{The new Reynolds stress}\label{subsec def Reynolds}
To construct a new subsolution we still need to define the new pressure $p_{q+1}$ and the new Reynolds stress $R_{q+1}$. We set $p_{q+1}=p_q-\overline{\sigma}_q^2\rho_q$. Substituting this, along with our ansatz for $v_{q+1}$ and $B_{q+1}$ into the definition of subsolution leads to the following equation for $R_{q+1}$:
\begin{equation}
\label{ec Rq+1} 
\begin{split}
\Div R_{q+1}&=\partial_tw_0+v_q\cdot\nabla w_0+(\Div w_0)v_q+w_0\cdot \nabla v_q+\partial_t w_c \\ &\hspace{12pt}+ \partial_t w_\phi+\Div(-\overline{\sigma}_q^2\rho_q\Id+R_q+w_0\otimes w_0)  \\&\hspace{12pt}+\Div\left[(w_c+w_\phi)\otimes v_{q+1}+v_{q+1}\otimes(w_c+w_\phi)-(w_c+w_\phi)\otimes(w_c+w_\phi)\right] \\ &\hspace{12pt}-\Div(b_{q+1}\otimes B_{q}+B_{q}\otimes b_{q+1}+b_{q+1}\otimes b_{q+1}).
\end{split}
\end{equation}
To solve this equation for $R_{q+1}$, we need a right inverse operator for the divergence. In Appendix~\ref{appendix divergence} we define an operator $\mathcal{R}$ that maps a field $f\in C^\infty_c(\RR^3,\RR^3)$ to a smooth symmetric matrix such that $\Div(\mathcal{R}f)=f$. Once we have introduced this operator, we can define
\begin{align}
S_1&\coloneqq \mathcal{R}(\partial_tw_0+v_q\cdot\nabla w_0+(\Div w_0)v_q+w_0\cdot \nabla v_q+\partial_t w_c) \equiv \mathcal{R}(z_1), \label{def S1} \\
S_2&\coloneqq \mathcal{R}(\partial_t w_\phi), \label{def S2}\\
S_3&\coloneqq \sum_{j=1}^6(\overline{\gamma}_{qj}^2-\gamma_{qj}^2)\,\zeta_j\otimes \zeta_j, \label{def S3}\\
S_4&\coloneqq \mathcal{R}\Div\hspace{-1.5pt}\left(w_0\otimes w_0-\sum_{j=1}^6\sigma_j(\mu\,\cdot\,)^2\hspace{1.5pt}\overline{\gamma}_{qj}^2\hspace{1.5pt}\zeta_j\otimes \zeta_j\right) \equiv \mathcal{R}(z_2), \label{def S4} \\
S_5&\coloneqq -\mathcal{R}\hspace{-1.5pt}\left(\sum_{j=1}^6[1-\sigma_j(\mu \,\cdot\,)^2]\Div\hspace{-1pt}\left(\overline{\gamma}_{qj}^2\hspace{0.5pt}\zeta_j\otimes \zeta_j\right)\right)\equiv -\mathcal{R}(z_3), \label{def S5}\\
S_6&\coloneqq (w_c+w_\phi)\otimes v_{q+1}+v_{q+1}\otimes(w_c+w_\phi)-(w_c+w_\phi)\otimes(w_c+w_\phi)\\&\hspace{13pt}-b_{q+1}\otimes B_{q}-B_{q}\otimes b_{q+1}-b_{q+1}\otimes b_{q+1}. \label{def S6}
\end{align}
Using \eqref{descomposición final R_q}, we see that \eqref{ec Rq+1} will hold if we set
\begin{equation}
R_{q+1}\coloneqq S_1+S_2+S_3+S_4+S_5+S_6.\label{def Rq+1}
\end{equation}

We claim that with our definitions for the new velocity $v_{q+1}$ and magnetic field $B_{q+1}$, the new Reynolds stress will be suitably small, meaning that \eqref{inductive cota Rq C0} and \eqref{inductive cota Rq C1} will hold for $q+1$. Sections~\ref{section estimates diffeo}-\ref{section matrix} are devoted to proving the necessary bounds. We anticipate that the most problematic term is $S_2$. Indeed, the other terms are similar to the expressions that one encounters when doing convex integration for the Euler equations, and dealing with them is by now fairly standard. However, estimating $S_2$ is extremely challenging due to the complex structure of $w_\phi$. This is one of the most delicate parts of our construction, along with estimating $w_\phi$ itself.  

\subsection{One final correction}\label{subsec final correction}
We will see that the constructed subsolution satisfies all of the desired estimates (\ref{inductive tamaño Bq})-(\ref{cambio helicidad q}). Furthermore, we have seen that the support of $(v_{q+1},B_{q+1},p_{q+1})$ is contained in $\Omega_{q,2}\times[0,T]$ because we are only perturbing inside $\Omega_{q,2}$. However, the matrices $S_1$, $S_2$, $S_4$ and $S_5$ are defined through a non-local operator, so they are not compactly supported, in general. Hence, the new Reynolds stress $R_{q+1}$ will not be supported on $\Omega_{q+1}\times[0,T]$, as required.

We would like to correct or modify $R_{q+1}$ to make it compactly supported, but there is an obstacle: angular momentum. Indeed, since $v_q$, $B_q$, $p_q$ and $R_q$ are compactly suported, it follows from (\ref{subs a}) and (\ref{identidad de Green matrices}) that
\begin{align}
\frac{d}{dt}\int_{\RR^3}\xi\cdot v_q&=\int_{\RR^3}\xi\cdot\Div\left(R_q+B_q\otimes B_q-v_q\otimes v_q-p_q\Id\right)\notag\\
&=\int_{\RR^3}\nabla_{\text{sym}}\xi: \left(R_q+B_q\otimes B_q-v_q\otimes v_q-p_q\Id\right)=0,
\label{momento angular constante}
\end{align}
for every time $t\in[0,T]$ and every vector field $\xi$ in the kernel of the operator $\nabla_{\text{sym}}$, which is defined in Appendix~\ref{appendix divergence}. By the same reasoning, the new velocity must also satisfy this condition if (\ref{inductive support q}) holds. However, the vector field that we have constructed need not satisfy this condition. Therefore, our goal is to define a new correction $w_L$ so that the new velocity
\begin{equation}
\tilde{v}_{q+1}\coloneqq v_{q+1}+w_L
\label{velocidad final}
\end{equation}
satisfies
\begin{equation}
\int_{\RR^3}\xi(x)\cdot[\tilde{v}_{q+1}(x,t)-v_q(x,t)]\,dx=0 \qquad \forall \xi\in \ker\nabla_{\text{sym}},\,\forall t\in[0,T].
\label{momento angular cero}
\end{equation}
Hence, it will follow from (\ref{momento angular constante}) that an analogous condition holds for $\tilde{v}_{q+1}$.

The kernel of the operator $\nabla_{\text{sym}}$ are the so-called Killing fields, and they are the generators of isometries of $\RR^3$. A basis $\{e_1,e_2,e_3,\xi_1,\xi_2,\xi_3\}$  of this vector space is given in Appendix~\ref{appendix divergence}; it consists of three constant vectors (the elements of the Cartesian basis, $e_i$), which are the generators of translations, and the generators of rotations, $\xi_i:=  e_i\times x$. Since the velocity is a compactly supported divergence-free field, it is mean-free, so we only have to worry about the product with the vectors $\xi_i$ for $i=1,2,3$. We define the angular momentum $L\in C^\infty([0,T],\RR^3)$ of the perturbation $w_0+w_c$ as
\begin{equation}
L(t)\coloneqq\int_{\RR^3}x\times (w_0+w_c)(x,t)\,dx,
\label{def L}
\end{equation}
where the integration is understood to be component-wise. By the defition of the vectors $\xi_i$, we have
\[L_i=\int_{\RR^3}e_i\cdot[x\times (w_0+w_c)]=\int_{\RR^3}(w_0+w_c)\cdot(e_i\times x)=\int_{\RR^3}(w_0+w_c)\cdot\xi_i.\]
The angular momentum of $w_\phi$ could be defined in a similar manner, but it would be harder to estimate. Fortunately, this field is so small that we can follow a different strategy.

To describe the support of the perturbation, we fix a point $x_q$ such that
\begin{equation}
B\left(x_q,\lambda_{q+1}^{-\tau}\bar{r}\}\right)\subset\joinrel\subset \Omega_{q+1}\backslash\overline{\Omega}_{q,2},
\label{perturbación en corona}
\end{equation}
which is possible because
\[\dist(\partial\Omega_{q+1}, \Omega_{q,2})=\left(2^{-q-2/3}-2^{-q-1}\right)\bar{r}\geq \frac{1}{8}2^{-q}\hspace{1pt}\bar{r} > 4\lambda_{q+1}^{-\tau}\bar{r}\]
for sufficiently large $a>1$ (depending on $\tau$). We fix a cutoff $\widetilde{\chi}\in C^\infty_c(B(0,1))$ such that $\int \widetilde{\chi}=1$ and we define
\[\chi_L(x)\coloneqq \lambda_{q+1}^{3\tau}\,\widetilde{\chi}(\lambda_{q+1}^\tau(x-x_q)).\]
We can finally define the correction $w_L$ as
\begin{equation}
\label{def wL}
w_L(x,t)\coloneqq -\frac{1}{2}\curl\left[\chi_L(x)L(t)\right]-\lambda_{q+1}^{4\tau}w_\phi(\lambda_{q+1}^\tau(x-x_q),t).
\end{equation}
Since the support of $\chi$ is contained in the unit ball and the support of $w_\phi$ is contained in $B(0,\bar{r})\times[0,T]$, it follows from (\ref{perturbación en corona}) that the suppor of $w_L$ is contained in $(\Omega_{q+1}\backslash\overline{\Omega}_{q,2})\times[0,T]$, as we wanted. Remember that we are assuming $\bar{r}\geq 1$.

Let us check that this perturbation produces the desired correction to the angular momentum. By the definition of the vector product and the curl, for any $A\in C^\infty_c(\RR^3,\RR^3)$ we have
\begin{align}
\begin{split}
\int_{\RR^3}(x\times \curl A)_i&=\int_{\RR^3}\varepsilon_{ijk}x_j\varepsilon_{klm}\partial_lA_m=-\int_{\RR^3}\varepsilon_{ijk}(\partial_lx_j)\varepsilon_{klm}A_m\\&=-\int_{\RR^3}\varepsilon_{ijk}\varepsilon_{kjm}A_m=\int_{\RR^3}\varepsilon_{ijk}\varepsilon_{mjk}A_m=2\int_{\RR^3}A_i,
\end{split}
\label{identidad momento angular}
\end{align}
where $\varepsilon_{ijk}$ is the Levi-Civita symbol and where summation over repeated indices is understood. Hence, we see that the first term in (\ref{def wL}) cancels the angular momentum of $w_0+w_c$, whereas the second term cancels the angular momentum of $w_\phi$. We conclude that if we define a new velocity as in (\ref{velocidad final}), it satisfies (\ref{momento angular cero}), as we wanted.

Since we have corrected the velocity, the rest of the elements of the construction must be modified accordingly in order to obtain a subsolution. First of all, taking into account that we have only changed the velocity outside of the support of $B_{q+1}$, the flow remains unchanged on the support of $B_{q+1}$, which means that (\ref{subs c}) is still satisfied. Thus, we do not need to modify the magnetic field. Regarding the Reynolds stress, we define
\begin{align}
S_7&\coloneqq\mathcal{R}(\partial_t w_L), \label{def S7}\\ S_8&\coloneqq w_L\otimes v_{q+1}+v_{q+1}\otimes w_L+w_L\otimes w_L. \label{def S8}
\end{align}
By (\ref{subs a}), we then have
\[\partial_t \tilde{v}_{q+1}+\Div(\tilde{v}_{q+1}\otimes \tilde{v}_{q+1}-B_{q+1}\otimes B_{q+1})+\nabla p_{q+1}=\Div (R_{q+1}+S_7+S_8),\]
which can be written as 
\[\partial_t \tilde{v}_{q+1}+\Div\left(\tilde{v}_{q+1}\otimes \tilde{v}_{q+1}-B_{q+1}\otimes B_{q+1}+p_{q+1}\Id-S_3-S_6-S_8\right)=\Div(S_1+S_2+S_4+S_5+S_7).\]
Note that the support of the matrix in parenthesis on the left-hand side is contained in $\Omega_{q+1}\times[0,T]$. By (\ref{momento angular cero}) and (\ref{identidad de Green matrices}), we conclude for all times $t\in[0,T]$ we have
\begin{equation}
\int_{\Omega_{q+1}}\xi(x)\cdot \Div(S_1+S_2+S_4+S_5+S_7)(x,t)\;dx=0 \qquad \forall \xi\in\ker\nabla_{\text{sym}}.
\label{momento angular suma matrices}
\end{equation}
By \cref{invertir divergencia matrices} and \cref{remark modificar matriz}, there exists a matrix $S_9\in C^\infty(\RR^3\times[0,T],\RR^{3\times 3}_\text{sym})$ such that $\Div S_9=0$ and
\[S_1+S_2+S_4+S_5+S_7+S_9\]
is supported on $\Omega_{q+1}\times[0,T]$. Hence, the matrix
\begin{equation}
\widetilde{R}_{q+1}\coloneqq R_{q+1}+S_7+S_8+S_9
\label{def tilde Rq+1}
\end{equation}
vanishes outside of $\Omega_{q+1}\times[0,T]$. Since $S_9$ is divergence-free, by construction of $S_7$ and $S_8$ we conclude that $(\tilde{v}_{q+1},B_{q+1},p_{q+1},\widetilde{R}_{q+1})$ is a subsolution satisfying (\ref{inductive support q}). In Sections~\ref{section estimates diffeo}-\ref{sec correction} we will derive the necessary estimates to ensure that it also satisfies the bounds \eqref{inductive tamaño Bq}-\eqref{cambio helicidad q}, thereby concluding \cref{prop steps}.

\section{Estimates on the diffeomorphism} \label{section estimates diffeo}
We begin by establishing some inequalities relating the various parameters, which will be used to ensure that some error terms encountered throughout the proof are small enough. The key inequalities are (\ref{param 1}) and (\ref{param 2}). The former motivates our choice of $\mu$, while the latter determines the relationship between $\lambda_q$ and $\lambda_{q+1}$, forcing us to choose a large value of $b$ and a small value of $\beta$. The other inequalities are just consequences of these choices.
\begin{lemma}
\label{lema relaciones entre parámetros}
Let $\tau>0$ be sufficiently small. There exists $\beta>\frac{1}{200}$ such that
\begin{align}
\frac{\delta_{q+1}^{3/2}\mu^3\lambda_{q+1}}{\eta^2}&\leq \delta_{q+2}\lambda_{q+1}^{-11\tau}, \label{param 1}	 \\
\frac{\delta_q^{1/2}\delta_{q+1}^{1/2}\lambda_q\lambda_{q+1}^{1+\tau}}{\mu\eta}&\leq \delta_{q+2}\lambda_{q+1}^{-11\tau}, \label{param 2} \\
\frac{\delta_{q+1}\mu^3}{\eta\lambda_{q+1}}\leq \frac{\delta_{q+1}\mu^2}{\eta}&\leq \delta_{q+2}\lambda_{q+1}^{-11\tau}, \label{param 3}\\
\frac{\delta_{q+1}^{3/2}\mu\lambda_{q+1}^2}{\eta^3}&\leq \delta_{q+2}\lambda_{q+1}^{-11\tau}, \label{param 4} \\
\frac{\delta_{q+1}\mu^2\lambda_{q+1}^{1+\tau}}{\eta^2}&\leq \delta_{q+2}\lambda_{q+1}^{-11\tau}, \label{param 5}\\
\frac{\delta_q^{1/2}\lambda_q}{\mu}\leq\frac{\delta_{q+1}^2\lambda_q\lambda_{q+1}^{11\tau}}{\delta_{q+2}\mu^{1-\tau}}&\leq \delta_{q+2}\lambda_{q+1}^{-11\tau}, \label{relación entre lambdas}\\
\frac{\delta_q^{1/2}\delta_{q+1}^{1/2}\lambda_q^2\lambda_{q+1}^2}{\eta^3}&\leq \delta_{q+2}\lambda_{q+1}^{-11\tau},\label{param 6}\\
\frac{\delta_q^{1/2}\delta_{q+1}\lambda_q\mu^2}{\eta^2}&\leq \delta_{q+2}\lambda_{q+1}^{-11\tau}, \label{param 7}\\
\frac{\delta_{q+1}^{1/2}\mu\lambda_{q+1}^\tau}{\eta} &\leq \delta_{q+2}\lambda_{q+1}^{-11\tau}. \label{param 9}
\end{align}

\end{lemma}
\begin{proof}
It is straighforward to check that (\ref{def mu}) and (\ref{def eta}) imply (\ref{param 1}) and (\ref{param 3}). Regarding (\ref{param 2}), we compute
\begin{align}
\log_a\left(\frac{\delta_q^{1/2}\delta_{q+1}^{1/2}\lambda_q\lambda_{q+1}^{1+12\tau}}{\delta_{q+2}\mu\eta}\right)&=\log_a\left(\frac{\delta_q^{1/2}\delta_{q+1}^{3/2}\lambda_q}{\delta_{q+2}^3\lambda_{q+1}^{1/3-36\tau}}\right) \nonumber \\&= 6\beta b^{q+2}-\beta b^q-3\beta b^{q+1}+b^{q}-\left(\frac{1}{3}-36\tau\right)b^{q+1} \label{aux log param} \\&= -b^q\left[\left(\frac{1}{3}-36\tau\right)b-1-\beta(6b^2-3b-1)\right], \nonumber
\end{align}
where we have substituted the definitions (\ref{def deltaq}), (\ref{def lambdaq}), (\ref{def mu}) and (\ref{def eta}). Evaluating the expression
\[\frac{\frac{1}{3}b-1}{6b^2-3b-1}\]	
at $b=6$, we see that if $\tau>0$ is chosen sufficiently small, there exists $\beta>\frac{1}{200}$ such that (\ref{aux log param}) is negative for any $q\geq 0$, so (\ref{param 2}) holds for such $\beta$. Concerning (\ref{param 4}) and (\ref{param 5}), substituting (\ref{def deltaq}), (\ref{def lambdaq}), (\ref{def mu}) and (\ref{def eta}) leads to 
\begin{align*}
	\frac{\delta_{q+1}^{3/2}\mu\lambda_{q+1}^{2+11\tau}}{\delta_{q+2}\eta^3}&=\frac{\delta_{q+1}^{5/2}}{\delta_{q+2}^3}\lambda_{q+1}^{-2/3+35\tau}\leq \lambda_{q+2}^{6\beta}\lambda_{q+1}^{-1/2}\leq \lambda_{q+1}^{6b\beta-1/2}, \\ \frac{\delta_{q+1}\mu^2\lambda_{q+1}^{1+12\tau}}{\delta_{q+2}\eta^2}&=	\frac{\delta_{q+1}}{\delta_{q+2}}\lambda_{q+1}^{-1/3+12\tau}\leq \lambda_{q+2}^{2\beta}\lambda_{q+1}^{-1/4}\leq\lambda_{q+1}^{2b\beta-1/4}, 
\end{align*}
where we have assumed that $210\tau\leq 1$. One can check that the chosen value of $\beta$ ensures that the exponents are negative, so (\ref{param 4}) and (\ref{param 5}) hold, too. Next, using \eqref{def lambdaq} and \eqref{def deltaq}, we have
\[\frac{\delta_q^{1/2}\delta_{q+2}}{\delta_{q+1}^2\mu^\tau\lambda_{q+1}^{11\tau}}\leq \frac{\delta_{q+2}}{\delta_{q+1}^2}=\lambda_{q+1}^{2\beta(2-b)}\leq 1,\]
so the first inequality in (\ref{relación entre lambdas}) holds. Regarding the second inequality, we have
\[\frac{\delta_{q+1}^2\lambda_q \lambda_{q+1}^{22\tau}}{\delta_{q+2}^2\mu^{1-\tau}}\leq\frac{\delta_{q+1}^{5/2}\lambda_q}{\delta_{q+2}^{3}\lambda_{q+1}^{1/3-35\tau}}=\lambda_q^{(-5b\beta+1+6b^2\beta-b/3)+35b\tau},\]
where we have used \eqref{def mu}. Substituting $b=6$ and assuming $199<\beta^{-1}<200$, one can check that the term in parenthesis is negative, so the exponent will be negative for $\tau>0$ sufficiently small. This yields the second inequality in (\ref{relación entre lambdas}). Concerning the remaining inequalities, dividing the left-hand side of (\ref{param 2}) by the left-hand side of (\ref{param 6}), we obtain
\[\frac{\eta^2}{\lambda_q\mu\lambda_{q+1}^{1-\tau}}\geq \frac{\eta^2}{\mu^2\lambda_{q+1}^{1-\tau}}=\lambda_{q+1}^{1/3+\tau}\geq 1,\]
so (\ref{param 2}) implies (\ref{param 6}). Similarly, dividing the left-hand side of (\ref{param 2}) by the left-hand side of (\ref{param 7}) yields
\[\frac{\eta\lambda_{q+1}^{1+\tau}}{\delta_{q+1}^{1/2}\mu^3}=\frac{\delta_{q+1}^{1/2}}{\delta_{q+2}^2}\lambda_{q+1}^{1+25\tau}\geq 1,\]
so (\ref{param 2}) implies (\ref{param 7}), too. Finally, for \eqref{param 9} we compute
\[\frac{\delta_{q+1}^{1/2}\mu\lambda_{q+1}^{12\tau}}{\delta_{q+2}\eta}=\frac{\delta_{q+1}^{1/2}}{\delta_{q+2}}\lambda_{q+1}^{-2/3+12\tau}=\lambda_{q+1}^{2b\beta-2/3-\beta+12\tau}.\]
With our choice of $b$ and $\beta$ the exponent will be negative for sufficiently small $\tau>0$, so \eqref{param 9} holds. This completes the proof.
\end{proof}

From now on, we assume that we have chosen $\tau>0$ small enough and that the parameter $\beta>0$ is the one given by the previous lemma, so that (\ref{param 1})-(\ref{param 3}) hold. This will ensure that several terms that appear in the following sections are sufficiently small. In particular, it follows from (\ref{relación entre lambdas}) that $\lambda_q\ll \mu$, as claimed.

Next, we derive some basic estimates for the amplitudes that appear in the definition of $\psi^0$:
\begin{lemma} \label{lemma estimates am}
	If $a>1$ is sufficiently large (depending on $\tau$), for any $j\in\{1,\dots, 6\}$ and any $m\in \Lambda$ we have
	\begin{align}
		\normxt{0}{a_{jm}}&\leq 2^{-6}\delta_{q+1}^{1/2}, \label{estimate am C0}\\
		\label{estimates am}
		\normxt{N}{a_{jm}}&\lesssim \delta_{q+1}^{1/2}\mu^N\qquad \forall N\geq 0.
	\end{align}
\end{lemma}
\begin{proof}
It follows from (\ref{def chim}) that 
\[\normxt{N}{\sigma_j(\mu\,\cdot\,)}+\normxt{N}{\chi_m}\lesssim \mu^{N} \qquad \forall N\geq 0.\]
We apply (\ref{estimate product}) to (\ref{def am}) using these bounds and \eqref{cotas gamatilde}, which yields \eqref{estimates am}. To obtain the improved bound for the $C^0$ norm, we use \eqref{cota gamatilde C0} and spend the extra $\lambda_q^{-\tau/2}$ to get rid of the constants by taking $a>1$ sufficiently large.
\end{proof}

These bounds on the coefficients directly translate into bounds on the diffeomorphism:
\begin{lemma}
For any $N\geq 0$ we have
\begin{equation}
	\label{estimates psi}
	\normxt{N}{\psi^0-\pi}\lesssim \delta_{q+1}^{1/2}\eta^{-1}\lambda_{q+1}^N,
\end{equation}
where $\pi$ was defined in (\ref{def pi}).
\end{lemma}
\begin{proof}
By \cref{remark soporte ajm disjuntos}, at most 16 of the coefficients $a_{jm}$ are nonzero at any given point. Thus, the estimate for the $C^0$-norm is a direct consequence of (\ref{hierarchy}) and (\ref{estimates am}). Regarding the $C^N$-norm for $N\geq 1$, we see that the second and higher order derivatives of $\theta_{jm}$ vanish because it is a polynomial of degree 1. In addition, $\normxt{1}{\theta_{jm}}\lesssim \lambda_{q+1}$. Hence, applying (\ref{estimate product}) with estimates (\ref{estimates am}) yields
\begin{align*}
\normxt{N}{\psi^0-\pi}&\lesssim \eta^{-1}\max_{m\in\Lambda}\left(\normxt{N}{a_{jm}}+\normxt{0}{a_{jm}}\lambda_{q+1}^N\right)\\&\hspace{12pt}+\eta^{-1}\lambda_{q+1}^{-1}\max_{m\in\Lambda}\left(\normxt{N+1}{a_{jm}}+\normxt{1}{a_{jm}}\lambda_{q+1}^N\right) \\&\lesssim \delta_{q+1}^{1/2}\eta^{-1}\lambda_{q+1}^N,
\end{align*}	
where we have used again that at most 16 of the coefficients $a_{jm}$ are nonzero in a neighborhood of any given point and (\ref{hierarchy}). 
\end{proof}

Next, we will study the Jacobian of $\psi^0$. As we have seen, $D\psi^0$ can be quite large. Fortunately, $\det(D\psi^0)$ will turn out to be very close to 1 due to cancellations in our construction. This is essential to ensure that the correction $\phi^c$ needed to obtain a volume-preserving map is small. If this were not the case, the perturbation to the velocity would change substantially and it would not produce the desired correction in the error matrix.
\begin{lemma}
	\label{lema fc}
	The function $f_c\coloneqq [\det(D\psi^0)-1]\in C^\infty_c(\Omega_{q,2}\times[0,T])$ satisfies
	\begin{align}
		&\int_{\Omega_{q,2}}f_c(x,t)\,dx=0 \qquad \forall t\in[0,T], \label{fc mean free} \\
		&\normxt{N}{f_c}\leq C(N)\, \frac{\delta_{q+1}\mu^2}{\eta^2}\lambda_{q+1}^{N} \qquad \forall N\geq 0. \label{estimates fc CN}
	\end{align}
	In addition, there exists a field $u_c\in C^\infty_c(\Omega_{q,2}\times[0,T],\RR^3)$ such that $\Div u_c=f_c$  satisfying the estimates
	\begin{equation}
		\normxt{N}{u_c}\leq C(N,\tau)\, \frac{\delta_{q+1}\mu^2}{\eta^2}\max\{\mu^N,\lambda_{q+1}^{N-1+\tau}\} \qquad \forall N\geq 0. \label{estimates uc CN}
	\end{equation}
\end{lemma}
\begin{proof}
We begin by obtaining an expression for $f_c$. First, it follows from (\ref{def psi}) that
\begin{align*}
	D\psi^0=\Id&-M_0+M_c\equiv \Id+M,
\end{align*}
where
\begin{align*}
	M_0&\coloneqq\sum_{m\in \Lambda}\left[\frac{\lambda_{q+1}}{\eta}\,a_{jm}\,\zeta_j\otimes k_{jm}\cos(\theta_{jm})+\frac{1}{\ell_m\eta}\,\zeta_j\otimes \nabla a_{jm}\sin(\theta_{jm})\right]\\ &\hspace{12pt}+\sum_{m\in \Lambda}\frac{1}{\ell_m\eta}\,[\nabla a_{jm}\times(\zeta_j\times k_{jm})]\otimes k_{jm}\sin(\theta_{jm}), \\[3pt]
	M_c&\coloneqq \sum_{m\in\Lambda}\frac{D[\nabla a_{jm}\times(\zeta_j\times k_{jm})]}{\ell_m^2\eta\lambda_{q+1}}\,\cos(\theta_{jm}).
\end{align*}
As a consequence of the Cayley-Hamilton theorem, we have the following formula for the determinant of $3\times 3$ matrices~\cite{Gant}:
\[\det(A)=\frac{1}{6}\tr(A)^3-\frac{1}{2}\tr(A)\tr(A^2)+\frac{1}{3}\tr(A^3).\]
Since we are working with symmetric matrices, it is easy to check that the formula holds by substituting $\tr(A^n)=\sum_{j=1}^3\lambda_{A,j}^n$, where $\lambda_{A,j}$ are the eigenvalues of $A$. In our case, we have:
\begin{align}
\begin{split}
f_c=\det\left(\Id+M\right)-1&= \tr(M)+\frac{1}{2}\tr(M)^2-\frac{1}{2}\tr(M^2)\\&\hspace{12pt}-\frac{1}{2}\tr(M)\tr(M^2)+\frac{1}{6}\tr(M)^3+\frac{1}{3}\tr(M^3).
\label{fc en función de las trazas general}
\end{split}
\end{align}
However, $\psi^0-\pi$ is divergence-free because it can be written as
\[\psi^0-\pi=\sum_{j,m}\curl\left[\frac{a_{jm}}{\ell_m^2\eta\lambda_{q+1}}(\zeta_j\times k_{jm})\cos(\theta_{jm})\right].\]
Here we have used that $k_{jm}\times(\zeta_j\times k_{jm})=\zeta_j$ because $k_{jm}$ is unitary and orthogonal to $\zeta_j$. We conclude that $M$ is trace-free, so \eqref{fc en función de las trazas general} simplifies to
\begin{equation}
\label{fc en función de las trazas}
f_c=-\frac{1}{2}\tr(M^2)+\frac{1}{3}\tr(M^3).
\end{equation}
Therefore, everything reduces to estimating the trace of $M^2$ and $M^3$. First, it follows from (\ref{estimates am}) that
\begin{equation}
\label{estimate Mc}
\normxt{0}{M_c}\lesssim \frac{\delta_{q+1}^{1/2}\mu^2}{\eta\lambda_{q+1}}.
\end{equation}
Concerning $M_0$, taking into account that $\zeta_j$ is orthogonal to $k_{jm}$ for any $m\in \Lambda$, we compute
\begin{align}
M_0^2&=\sum_{j=1}^6\sum_{n,m\in\Lambda}\bigg[\;\frac{\lambda_{q+1}}{\ell_m\eta^2}\,a_{jn}\, k_{jn}\cdot[\nabla a_{jm}\times(\zeta_j\times k_{jm})]\,\zeta_j\otimes k_{jm} \cos(\theta_{jn})\sin(\theta_{jm}) \nonumber \\ &\hspace{66pt}+\frac{\lambda_{q+1}}{\ell_m\eta^2}\,a_{jm}(\zeta_j\cdot\nabla a_{jn})\,\zeta_j\otimes k_{jm}\sin(\theta_{jn})\cos(\theta_{jm}) \nonumber\\ &\hspace{66pt}+\frac{1}{\ell_m\ell_n\eta^2}\,(\zeta_j\cdot\nabla a_{jn})\,\zeta_j\otimes \nabla a_{jm}\sin(\theta_{jn})\sin(\theta_{jm}) \label{expresión M0^2} \\ &\hspace{66pt}+\frac{1}{\ell_m\ell_n\eta^2}\,\nabla a_{jn}\cdot[\nabla a_{jm}\times (\zeta_j\times k_{jm})]\, \zeta_j\otimes k_{jm} \sin(\theta_{jn})\sin(\theta_{jm}) \nonumber\\ &\hspace{66pt}+\frac{1}{\ell_m\ell_n\eta^2}\,k_{jn}\cdot[\nabla a_{jm}\times(\zeta_j\times k_{jm})]\,[\nabla a_{jn}\times(\zeta_j\times k_{jn})]\otimes k_{jm} \sin(\theta_{jn})\sin(\theta_{jm})\,\bigg]. \nonumber
\end{align}
Note that we only have to sum over a single value of $j$ due to \cref{remark soporte ajm disjuntos}. Again, since $\zeta_j$ and $k_{jm}$ are orthogonal, we have $\tr(\zeta_j\otimes k_{jm})=0$. Therefore, only the third and fifth terms contribute to the trace, so
\[\normxt{0}{\tr(M_0^2)}\lesssim \frac{\delta_{q+1}\mu^2}{\eta^2},\]
where we have used (\ref{estimates am}). We conclude:
\begin{equation}
	\begin{split}
		\normxt{0}{\tr(M^2)}&\lesssim \normxt{0}{\tr(M_0^2)}+\normxt{0}{M_0}\normxt{0}{M_c}+\normxt{0}{M_c}^2 \\ &\lesssim \frac{\delta_{q+1}\mu^2}{\eta^2}+\left(\frac{\delta_{q+1}^{1/2}\lambda_{q+1}}{\eta}\right)\left(\frac{\delta_{q+1}^{1/2}\mu^2}{\eta\lambda_{q+1}}\right) \lesssim \frac{\delta_{q+1}\mu^2}{\eta^2}.
	\end{split}
	\label{traza M^2}
\end{equation}

Let us now study $M^3$. Since $\tr(AB)=\tr(BA)$, when estimating $\tr{M^3}$ we only have to account for four kind of terms: $M_0^3$, $M_0^2M_c$, $M_0M_c^2$ and $M_c^3$. To estimate $\tr(M_0^3)$, we consider the product $M_0\cdot M^2_0$:
\begin{itemize}
	\item Due to the orthogonality between $\zeta_j$ and the vectors $k_{jm}$, the first term in $M_0$ only yields nonzero terms when multiplied by the last term in $M_0^2$. However, they will be proportional to $\zeta_j\otimes k_{jm}$, which does not contribute to the trace.
	\item Since the terms $\zeta_j\otimes k_{jm}$ do not contribute to the trace, the second term in $M_0$ only contributes to the trace when multiplied by the third term in $M_0^2$.
	\item Again, due to the orthogonality between $\zeta_j$ and the vectors $k_{jm}$, the last term in $M_0$ only contributes to the trace when multiplied by the last term in $M_0^2$.
\end{itemize}
Hence, we conclude that \[\normxt{0}{\tr(M_0^3)}\lesssim \left(\frac{\delta_{q+1}^{1/2}\mu}{\eta}\right)^3.\]
Meanwhile, using the bound \eqref{estimates am} in \eqref{expresión M0^2} yields:
\[\normxt{0}{M_0^2M_c}\leq \normxt{0}{M_0^2}\normxt{0}{M_c}\lesssim \left(\frac{\delta_{q+1}\mu\lambda_{q+1}}{\eta^2}\right)\frac{\delta_{q+1}^{1/2}\mu^2}{\eta\lambda_{q+1}} = \left(\frac{\delta_{q+1}^{1/2}\mu}{\eta}\right)^3\]
where we have also used \eqref{estimate Mc}. The terms $M_0M_c^2$ and $M_c^3$ satisfy even better bounds, so we conclude
\[\normxt{0}{\tr(M^3)}\lesssim \left(\frac{\delta_{q+1}^{1/2}\mu}{\eta}\right)^3.\]
Substituting this and \eqref{traza M^2} into (\ref{fc en función de las trazas}), we finally obtain
\[\normxt{0}{f_c}\lesssim \frac{\delta_{q+1}\mu^2}{\eta^2}.\]
In fact, since we get at most an extra $\lambda_{q+1}$ factor every time we differentiate any of the terms in the product, we see that (\ref{estimates fc CN}) holds. In particular, we can ensure that $\normxt{0}{f_c}\leq 1/2$ if we take $a>1$ sufficiently large.

Let us now prove that $f_c(\cdot,t)$ is mean-free for any time $t\in[0,T]$. Since $\psi^0$ is a diffeomorphism of $\RR^3$ that is the identity in a neighborhood of $\RR^3\backslash\Omega_{q,2}$, it is also a diffeomorphism of $\Omega_{q,2}$. By the change of variables formula we have
\[\int_{\Omega_{q,2}}1\,dx=\int_{\Omega_{q,2}}\abs{\det(D\psi^0)}\,dx=\int_{\Omega_{q,2}}[1+f_c(x,t)]\,dx,\]
where we have used that $\normxt{0}{f_c}\leq 1/2$. We conclude (\ref{fc mean free}). 

By \cref{invertir div vectores}, there exists $u_c\in C^\infty_c(\Omega_{q,2}\times[0,T],\RR^3)$ such that $\Div u_c=f_c$ and satisfying the estimates (\ref{cotas invertir div vectores}). Therefore, it follows from the interpolation of (\ref{estimates fc CN}) that for any $M\geq0$ and $N\geq 1$ we have
\begin{align}
\begin{split}
\label{aux uc N>0}
\normx{N+\tau}{\partial_t^Mu_c(\cdot,t)}&\lesssim \normx{N-1+\tau}{\partial_t^Mf_c(\cdot,t)}\lesssim \normxt{M+N-1+\tau}{f_c}\\&\lesssim \frac{\delta_{q+1}\mu^2\lambda_{q+1}^{M+N-1+\tau}}{\eta^2}.
\end{split}
\end{align}
In order to cover the case $N=0$, we must first derive suitable bounds on the Besov norms of $f_c$. To do so, we must separate into high and low frequencies. Taking into account the form of each of the terms in the product, we see that we can write
\begin{equation}
\label{descomposición fc}
f_c=g_0^0+\sum_{j,l} g_{jl}^0\,\sin(\tilde{\theta}_{jl})
\end{equation}
where $\tilde{\theta}_{jl}$ is a sum of up to three phases $\theta_{jm}$ with $m\in\Lambda$ plus possibly a constant phase $\frac{\pi}{2}$, and
\[\normxt{N}{g_0^0}+\max_{j,l} \normxt{N}{g_{j,l}^0}\lesssim \frac{\delta_{q+1}\mu^{2+N}}{\eta^2}.\]
Since the supports of at most $16$ of the coefficients $a_{jm}$ have nonempty intersection, we see that there must exists a numerical constant $J_\ast$ such that the supports of at most $J_\ast$ of the functions $g^0_{jl}$ have nonempty intersection. The value of $J_\ast$ depends on the combinatorics involved in computing the determinant.

The low-frequency term $g_0$ appears due to destructive interference between two terms that are multiplied by a trigonometric function with the same argument $\theta_{jm}$. This also leads to a high-frequency term with $\tilde{\theta}_{jl}=2\theta_{jm}$. On the other hand, the coefficients $\ell_m$ were constructed so that $\ell_m\neq \ell_n$ whenever the supports of two different amplitudes $a_{jm}$, $a_{jn}$ have nonempty intersection. As a result, there cannot be destructive interference: 
\[|\nabla\tilde{\theta}_{jl}|\sim \lambda_{q+1}.\]
Similarly, the sum of three coefficients $\ell_m$ is nonzero, by construction. Thus, the product of three oscillatory terms also leads to a high-frequency term satisfying the above estimate.

It is easy to see that $\partial_t^Mf_c$ has an analogous decomposition to (\ref{descomposición fc}), but the corresponding amplitudes $g_0^M$, $g_{jl}^M$ will satisfy
\begin{align*}
\normxt{N}{g_0^M}&\lesssim  \frac{\delta_{q+1}\mu^{2+N+M}}{\eta^2}, \\
\max_{j,l} \normxt{N}{g_{jl}^M}&\lesssim \frac{\delta_{q+1}\mu^{2+N}\lambda_{q+1}^M}{\eta^2}.
\end{align*}
By \cref{stationary phase lemma}, for any integer $m\geq 0$ we then have
\begin{align*}
\norm{\partial_t^Mf_c(\cdot,t)}_{B^{-1+\tau}_{\infty,\infty}}&\lesssim \normxt{0}{g_0^M}+\max_{j,l}\left(\lambda_{q+1}^{-1+\tau}\normx{0}{g_{jl}^M(\cdot,t)}+\lambda_{q+1}^{-m}\normx{m}{g_{jl}^M}\right) \\
&\lesssim \frac{\delta_{q+1}\mu^{2+M}}{\eta^2}+\frac{\delta_{q+1}\mu^2\lambda_{q+1}^{M-1+\tau}}{\eta^2}+\frac{\delta_{q+1}\mu^{2+m}\lambda_{q+1}^{M-m}}{\eta^2},
\end{align*}
where we have estimated the Besov norm of the low-frequency term by its $C^0$-norm. Choosing $m=2$, it follows from (\ref{def mu}) that $\mu^2\leq \lambda_{q+1}$, so
\[\normx{\tau}{\partial_t^Mu_c(\cdot,t)}\lesssim \norm{\partial_t^Mf_c(\cdot,t)}_{B^{-1+\tau}_{\infty,\infty}} \lesssim \frac{\delta_{q+1}\mu^2}{\eta^2}\max\{\mu^M,\lambda_{q+1}^{M-1+\tau}\}.\]
Combining this with (\ref{aux uc N>0}), we obtain
\[\normxt{N}{u_c}\lesssim \sum_{0\leq M_1+M_2\leq N} \max_{t\in [0,T]}\normx{M_2}{\partial_t^{M_1}u_c(\cdot,t)}\lesssim \frac{\delta_{q+1}\mu^2}{\eta^2}\max\{\mu^N,\lambda_{q+1}^{N-1+\tau}\},\]
which completes the proof.
\end{proof}

Due to the previous lemma and the inverse function theorem, each point has a sufficiently small neighborhood where $\psi^0_t$ has a smooth inverse, meaning $\psi^0_t$ is a local diffeomorphism. Since $\psi^0$ depends smoothly on $t$, so do the local inverses. In the next lemma we will verify that $\psi^0_t$ is, in fact, a global diffeomorphism. 

First, we introduce some notation. Many of the functions that appear are be multiplied by one of the coefficients $a_{jm}$ or by its derivatives. Therefore, only their behavior on the support of $a_{jm}$ is relevant (or on a slightly enlarged set, when we compose with a suitable diffeomorphism). For this reason, it is convenient in some cases to work with the following localized norms:
\begin{definition}
For $j\in \{1,\dots,6\}$, $m\in \Lambda$ and $N\geq 0$, we define
\begin{align*}
G_{jm}&\coloneqq \supp a_{jm}+B(0,c_1\mu^{-1}), \\ \normxtjm{N}{\,\cdot\,}&\coloneqq \norm{\,\cdot\,}_{C^N(G_{jm})},
\end{align*}
where $c_1$ is given by Lemma~\ref{lemma mikado}.
\end{definition}
Note that by \cref{remark soporte ajm disjuntos}, we have
\begin{equation}
G_{jm}\cap G_{j'm'}=\varnothing \qquad  \forall j\neq j'.
\label{Gjm disjuntos}
\end{equation}

We are now ready to prove that $\psi^0$ is invertible:
\begin{lemma} \label{identidad inversa}
	There exists a map $\phi^0\in C^\infty(\RR^3\times[0,T],\RR^3)$ such that $\phi^0_t$ is the inverse of $\psi^0_t$. Furthermore, it satisfies 
	\[\phi^0=\pi+\sum_{j, m}\left[\frac{a_{jm}\compt\phi^0}{\ell_m\eta}\hspace{0.5pt}\zeta_j\,\sin(\theta_{jm}\compt\phi^0)-\frac{(\nabla a_{jm}\compt \phi^0)\times(\zeta_j\times k_{jm})}{\ell_m^2\eta\lambda_{q+1}}\hspace{0.5pt}\cos(\theta_{jm}\compt \phi^0)\right],\]
where $\pi$ was defined in (\ref{def pi}).
\end{lemma}
\begin{proof}
Let us define a distance on $C^0(\RR^3\times[0,T],\RR^3)$ through the following expression:
\[d(f,g)\coloneqq \normxt{0}{f-g}+\mu^{-1}\max_{j, m}\normxtjm{0}{\theta_{jm}\compt f-\theta_{jm}\compt g}\]
and consider the space
\[X\coloneqq \{f\in C^0(\RR^3\times[0,T],\RR^3): d(f,\pi)\leq c_1\mu^{-1}\}.\]
It is easy to check that $d(\,\cdot\,,\,\cdot\,)$ is, indeed, a distance and that $X$ endowed with this distance is a complete metric space. Next, let us consider the map $T: C^0(\RR^3\times[0,T],\RR^3)\to C^0(\RR^3\times[0,T],\RR^3)$ defined in (\ref{def T contractiva}). We will show that the image of the restriction $T|_X$ is contained in $X$ and that $T|_X:X\to X$ is a contractive map.

First of all, it follows from the definition of $T$ that
\[T\pi=\pi+\sum_{j, m}\left[\frac{a_{jm}}{\ell_m\eta}\hspace{0.5pt}\zeta_j\,\sin(\theta_{jm})-\frac{(\nabla a_{jm})\times(\zeta_j\times k_{jm})}{\ell_m^2\eta\lambda_{q+1}}\hspace{0.5pt}\cos(\theta_{jm})\right].\]
Hence, by \eqref{estimates am} we have $\normxt{0}{T\pi-\pi}\lesssim \delta_{q+1}^{1/2}\eta^{-1}$. Next, we fix $j\in\{1,\dots 6\}$ and $m\in\Lambda$. It follows from definition \eqref{def thetam} that:
\[\theta_{jm}\compt (T\pi)-\theta_{jm}=\sum_{j', n}\left[\frac{\ell_m\lambda_{q+1} }{\ell_n\eta}\hspace{0.5pt}a_{j'n}\hspace{0.5pt}(k_{jm}\cdot\zeta_{j'})\,\sin(\theta_{j'n})-\frac{\ell_m}{\ell_n^2\eta}\hspace{0.5pt}[k_{jm}\cdot(\nabla a_{j'n})\times(\zeta_j\times k_{j'n})]\hspace{0.5pt}\cos(\theta_{j'n})\right].\]
If we evaluate on $(x,t)\in G_{jm}$, only the terms with $j'=j$ will contribute, due to \eqref{Gjm disjuntos}. Furthermore, the first term in brackets will vanish because $\zeta_j$ and $k_{jm}$ are orthogonal. Using \eqref{estimates am}, we obtain
\[\normxtjm{0}{\theta_{jm}\compt (T\pi)-\theta_{jm}}\lesssim \delta_{q+1}^{1/2}\mu\eta^{-1}.\]
Taking the maximum in $j\in\{1,\dots 6\}$ and $m\in\Lambda$, we obtain
\[d(T\pi, \pi)\leq C \delta_{q+1}^{1/2}\eta^{-1}\leq \frac{1}{2}c_1\mu^{-1}\]
for some constant $C>0$ and $a>1$ sufficiently large. In particular, $T\pi\in X$.

We will now estimate the distance between $Tf$ and $Tg$ for $f,g\in X$. This will allow us to prove that the image of $T|_X$ is contained in $X$ and that it is a contractive map. First of all, note that for any smooth maps $f,g, h_1, h_2$, one has
\begin{align}
\nonumber
\|(h_1\compt &f)[h_2\compt(\theta_{jm}\compt f)]-(h_1\compt g)[h_2\compt(\theta_{jm}\compt g)]\|_{C^0_{x,t}} \lesssim \\
\begin{split}
\label{aux invertir psi0}
&\lesssim \normxt{0}{(h_1\compt f)[h_2\compt(\theta_{jm}\compt f)]-(h_1\compt g)[h_2\compt(\theta_{jm}\compt f)]}  \\
&\hspace{12pt} +\normxt{0}{(h_1\compt g)[h_2\compt(\theta_{jm}\compt f)]-(h_1\compt g)[h_2\compt(\theta_{jm}\compt g)]}
\end{split}  \\
&  \lesssim \normxt{0}{h_2}\normxt{1}{h_1}\normxt{0}{f-g}+\normxt{0}{h_1}\normxt{1}{h_2}\normxt{0}{\theta_{jm}\compt f-\theta_{jm}\compt g}, \nonumber 
\end{align}
where the last norm can be replaced by $\normxtjm{0}{\,\cdot\,}$ if the support of $h_1$ is contained in $\supp a_{jm}$ and $f,g\in X$. Indeed, in that case, the support of $h_1\compt f$ and $h_1\compt g$ is contained in $G_{jm}$.

Next, it follows from the definition of $T$ that
\begin{align*}
Tf-Tg&=\sum_{j,m}\left[\frac{a_{jm}\compt f}{\ell_m\eta}\zeta_j\hspace{0.5pt}\sin(\theta_{jm}\compt f)-\frac{a_{jm}\compt g}{\ell_m\eta}\zeta_j\hspace{0.5pt}\sin(\theta_{jm}\compt g)\right] \\&\hspace{12pt}-\sum_{j,m}\left[\frac{(\nabla a_{jm}\compt f)\times(\zeta_j\times k_{jm})}{\ell_m^2\eta\lambda_{q+1}}\hspace{0.5pt}\cos(\theta_{jm}\compt f)+\sum_{j,m}\frac{(\nabla a_{jm}\compt g)\times(\zeta_j\times k_{jm})}{\ell_m^2\eta\lambda_{q+1}}\hspace{0.5pt}\cos(\theta_{jm}\compt g)\right].
\end{align*}
Hence, by (\ref{aux invertir psi0}), we have
\begin{align*}
	\normxt{0}{Tf-Tg}&\lesssim \max_{j,m}\bigg(\eta^{-1}\normxt{1}{a_{jm}}\normxt{0}{f-g}+\eta^{-1}\normxt{0}{a_{jm}}\normxtjm{0}{\theta_{jm}\compt f-\theta_{jm}\compt g}\\&\hspace{40pt}+\eta^{-1}\lambda_{q+1}^{-1}\normxt{2}{a_{jm}}\normxt{0}{f-g}+\eta^{-1}\lambda_{q+1}^{-1}\normxt{1}{a_{jm}}\normxtjm{0}{\theta_{jm}\compt f-\theta_{jm}\compt g}\bigg) \\ &\lesssim \frac{\delta_{q+1}^{1/2}\mu}{\eta}\,d(f,g),
\end{align*}
where we have used (\ref{estimates am}). Next, since $\zeta_j$ is orthogonal to $k_{jn}$ for any $n\in \Lambda$, we may write
\begin{align*}
\theta_{jm}\compt (Tf)\big|_{G_{jm}}-\theta_{jm}\compt (Tg)\big|_{G_{jm}}=&-\sum_{n\in \Lambda}\frac{l_m k_m\cdot[(\nabla a_{jn}\compt f)\times(\zeta_j\times k_{jn})]}{\ell_n^2\eta}\cos(\theta_{jn}\compt f) \\ &+\sum_{n\in \Lambda}\frac{l_m k_m\cdot[(\nabla a_{jn}\compt g)\times(\zeta_j\times k_{jn})]}{\ell_n^2\eta}\cos(\theta_{jn}\compt g)
\end{align*}
by the same reasoning as with $\theta_{jm}\compt (T\pi)$. Thus, using (\ref{estimates am}) and (\ref{aux invertir psi0}) leads to
\[\normxtjm{0}{\theta_{jm}\compt (Tf)-\theta_{jm}\compt (Tg)}\lesssim \frac{\delta_{q+1}^{1/2}\mu^2}{\eta}\,d(f,g).\]
We conclude that
\[d(Tf,Tg)\lesssim \frac{\delta_{q+1}^{1/2}\mu}{\eta}\,d(f,g)\leq \delta_{q+1}^{1/2}\,d(f,g) \qquad \forall f,g\in X.\]
Therefore, for sufficiently large $a>1$:
\[d(Tf,\pi)\leq d(Tf,T\pi)+d(T\pi,\pi)\leq \delta_{q+1}^{1/2}d(f,\pi) +\frac{1}{2}c_1\mu^{-1} \leq c_1\mu^{-1}\]
for all $f\in X$. Thus, the image of $T|_X$ is contained in $X$. Furthermore, we see that $T|_X:X\to X$ is a contractive mapping if $a$ is sufficiently large. Hence, by Banach's fixed point theorem, there exists a unique $\phi^0\in X$ such that $T\phi^0=\phi^0$. We will now check that $\phi^0_t$ is the inverse of $\psi^0_t$. Note that $\psi^0$ has the same expression as $T\pi$ up to the sign of some terms. Therefore, by the same argument, $\psi^0\in X$. Using (\ref{f pi}), definition (\ref{def psi}) and the fact that the operation $\compt$ is associative, we compute
\begin{align*}
	\phi^0\compt \psi^0-\pi&=\sum_{j,m}\left[\frac{a_{jm}\compt (\phi^0\compt \psi^0)}{\ell_m\eta}\zeta_j\hspace{0.5pt}\sin[\theta_{jm}\compt (\phi^0\compt \psi^0)]-\frac{a_{jm}}{\ell_m\eta}\zeta_j\hspace{0.5pt}\sin(\theta_{jm})\right] \\&\hspace{12pt}-\sum_{j,m}\frac{[\nabla a_{jm}\compt (\phi^0\compt \psi^0)]\times(\zeta_j\times k_{jm})}{\ell_m^2\eta\lambda_{q+1}}\hspace{0.5pt}\cos[\theta_{jm}\compt (\phi^0\compt \psi^0)]
	\\&\hspace{12pt}+\sum_{j,m}\frac{\nabla a_{jm}\times(\zeta_j\times k_{jm})}{\ell_m^2\eta\lambda_{q+1}}\hspace{0.5pt}\cos(\theta_{jm}).
\end{align*}
Taking into account that $\zeta_j$ is orthogonal to $k_{jn}$ for any $n\in \Lambda$, we also have
\begin{align*}
	\theta_{jm}\compt (\phi^0\compt \psi^0)\big|_{G_{jm}}-\theta_{jm}\big|_{G_{jm}}=&-\sum_{n\in \Lambda}\frac{l_m k_m\cdot[(\nabla a_{jn}\compt (\phi^0\compt \psi^0))\times(\zeta_j\times k_{jn})]}{\ell_n^2\eta}\cos[\theta_{jn}\compt (\phi^0\compt \psi^0)] \\ &+\sum_{n\in \Lambda}\frac{l_m k_m\cdot[\nabla a_{jn}\times(\zeta_j\times k_{jn})]}{\ell_n^2\eta}\cos(\theta_{jn}).
\end{align*}
Therefore, arguing as before leads to
\[d(\phi^0\compt \psi^0, \pi)\lesssim \frac{\delta_{q+1}^{1/2}\mu}{\eta}d(\phi^0\compt \psi^0,\pi)\leq \delta_{q+1}^{1/2}\hspace{1pt}d(\phi^0\compt \psi^0,\pi).\]
Thus, we will have $d(\phi^0\compt \psi^0,\pi)\leq \frac{1}{2}d(\phi^0\compt \psi^0,\pi)$ for sufficiently large $a>1$, so $\phi^0\compt \psi^0=\pi$. Meanwhile, direct computation of the composition shows that $\psi^0\compt\phi^0=\pi$. We conclude that $\psi^0_t$ is invertible and $\phi^0_t$ is its inverse. Since the local inverses were smooth, $\phi^0$ is automatically smooth.
\end{proof}

The expression derived in the previous lemma will allow us to obtain good estimates for $\phi^0$. This task is still nontrivial because $\phi^0$ appears in the right-hand side of the expression. The most problematic term is $\theta_{jm}\compt \phi^0$, whose derivatives may increase due to the composition with $\phi^0$. Fortunately, we will see that certain cancellations ensure that the change is not too large.
\begin{lemma}
For any $N\geq0$, the diffeomorphism $\phi^0$ satisfies
\begin{equation}
\label{estimates phi0}
\normxt{N}{\phi^0-\pi}+\frac{1}{\mu}\max_{j,m}\normxtjm{N}{\theta_{jm}\compt\phi^0-\theta_{jm}}\lesssim \frac{\delta_{q+1}^{1/2}\lambda_{q+1}^N}{\eta}.
\end{equation}
\end{lemma}
\begin{proof}
The bound for the $C^0$-norm of the difference $\phi^0-\pi$ follows from \cref{identidad inversa} using (\ref{hierarchy}) and (\ref{estimates am}). Similarly, the bound for the $C^0$-norm of the restriction of $\theta_{jm}\compt\phi^0-\theta_{jm}$ to $G_{jm}$ follows from
\begin{equation}
\theta_{jm}\compt\phi^0\big|_{G_{jm}}-\theta_{jm}\big|_{G_{jm}}=-\sum_{n\in \Lambda}\frac{\ell_m k_{jm}\cdot[(\nabla a_{jn}\compt\phi^0)\times(\zeta_j\times k_{jn})]}{\ell_n^2\eta}\cos(\theta_{jn}\compt\phi^0).
\label{diferencia thetam con y sin phi0}
\end{equation}
As in the previous lemma, we have used that only the terms with the same value of $j$ contribute to the sum, due to \eqref{Gjm disjuntos}, and that that $\zeta_j$ is orthogonal to $k_{jn}$ for any $n\in \Lambda$. 

Let us now establish (\ref{estimates phi0}) for $N\geq1$. First, by (\ref{estimate product}), we have
\begin{align}
\normxt{N}{\phi^0-\pi}&\lesssim \frac{1}{\eta}\max_{j,m}\left(\normxt{N}{a_{jm}\compt\phi^0}+\normxt{0}{a_{jm}}\normxtjm{N}{\sin(\theta_{jm}\compt\phi^0)} \right) \nonumber \\&\hspace{12pt}+\frac{1}{\eta\lambda_{q+1}}\max_{j,m}\left(\normxt{N}{\nabla a_{jm}\compt\phi^0}+\normxt{1}{a_{jm}}\normxtjm{N}{\cos(\theta_{jm}\compt\phi^0)}\right). \label{aux estimates phi0 1}
\end{align}
Here we have used that $a_{jm}\compt \phi^0$ and $\nabla a_{jm}\compt \phi^0$ are supported on $G_{jm}$ if $a>1$ is sufficiently large, which follows from the bound for the $C^0$-norm of the difference $\phi^0-\pi$. As a result, the behavior of $\theta_{jm}\compt\phi^0$ only matters on $G_{jm}$. Meanwhile, applying (\ref{estimate product}) to (\ref{diferencia thetam con y sin phi0}) yields:
\begin{equation}
\normxtjm{N}{\theta_{jm}\compt\phi^0-\theta_{jm}}\lesssim \frac{1}{\eta}\max_{n\in \Lambda}\left(\normxt{N}{\nabla a_{jn}\compt\phi^0}+\normxt{1}{a_{jn}}\norm{\cos(\theta_{jm}\compt\phi^0)}_{N,\hspace{0.5pt}jn}\right). \label{aux estimates phi0 2}
\end{equation}
To estimate each term, we use (\ref{estimate composition}):
\begin{align*}
\seminormxt{N}{a_{jm}\compt\phi^0}&\lesssim \seminormxt{1}{a_{jm}}\left(1+\normxt{N-1}{\partial_t\phi^0}+\normxt{N-1}{D\phi^0}\right)+\normxt{N}{a_{jm}}\seminormxt{1}{\phi^0}^N \\
&\lesssim \normxt{1}{a_{jm}}\left(1+\normxt{N}{\phi^0-\pi}\right)+\normxt{N}{a_{jm}}\left(1+\normxt{1}{\phi^0-\pi}^N\right), \\[5pt]
\left[\sin(\theta_{jm}\compt\phi^0)\right]_{C^N_{x,t}(G_{jm})}&\lesssim \norm{\partial_t(\theta_{jm}\compt\phi^0)}_{C^{N-1}_{x,t}(G_{jm})}+\norm{D(\theta_{jm}\compt\phi^0)}_{C^{N-1}_{x,t}(G_{jm})}+\left[\theta_{jm}\compt\phi^0\right]_{C^1_{x,t}(G_{jm})}^N \\
&\lesssim \lambda_{q+1}+\norm{\theta_{jm}\compt\phi^0-\theta_{jm}}_{C^N_{x,t}(G_{jm})}+\left[\theta_{jm}\compt\phi^0\right]_{C^1_{x,t}(G_{jm})}^N,
\end{align*}
where we have used that $\seminormxt{1}{\theta_{jm}}\lesssim \lambda_{q+1}$. The terms $\nabla a_{jm}\compt\phi^0$ and $\cos(\theta_{jm}\compt \phi^0)$ satisfy analogous bounds. Let us focus on $N=1$. Substituting these estimates into our previous expressions and using (\ref{estimates am}), we obtain
\begin{align*}
\normxt{1}{\phi^0-\pi}&\lesssim\frac{\delta_{q+1}^{1/2}\mu}{\eta}\normxt{1}{\phi^0-\pi}+\frac{\delta_{q+1}^{1/2}\lambda_{q+1}}{\eta}+\frac{\delta_{q+1}^{1/2}}{\eta}\max_{j,m}\normxtjm{1}{\theta_{jm}\compt\phi^0-\theta_{jm}}, \\[3pt] \normxtjm{1}{\theta_{jm}\compt\phi^0-\theta_{jm}}&\lesssim\frac{\delta_{q+1}^{1/2}\mu^2}{\eta}\normxt{1}{\phi^0-\pi}+\frac{\delta_{q+1}^{1/2}\mu\lambda_{q+1}}{\eta}+\frac{\delta_{q+1}^{1/2}\mu}{\eta}\max_{j,n}\norm{\theta_{jn}\compt\phi^0-\theta_{jn}}_{1,\hspace{0.5pt}jn}.
\end{align*}
Taking the maximum in $j\in\{1,\dots, 6\}$ and $m\in \Lambda$ and combining both bounds, we obtain
\begin{align*}
	&\normxt{1}{\phi^0-\pi}+\frac{1}{\mu}\max_{j,m}\normxtjm{1}{\theta_{jm}\compt\phi^0-\theta_{jm}}\lesssim\frac{\delta_{q+1}^{1/2}\lambda_{q+1}}{\eta}+\\&\hspace{195pt}+\frac{\delta_{q+1}^{1/2}\mu}{\eta}\left(\normxt{1}{\phi^0-\pi}+\frac{1}{\mu}\max_{j,m}\normxtjm{1}{\theta_{jm}\compt\phi^0-\theta_{jm}}\right).
\end{align*}
Since the coefficient multiplying the term in parenthesis can be made arbitrarily small by increasing $a>1$, due to (\ref{hierarchy}), we conclude (\ref{estimates phi0}) for $N=1$. In particular, we have
\[\left[\theta_{jm}\compt\phi^0\right]_{C^1_{x,t}(G_{jm})}\lesssim \lambda_{q+1}.\]
We can substitute this into our previous estimates, obtaining
\begin{align*}
\seminormxt{N}{a_{jm}\compt\phi^0}&\lesssim \delta_{q+1}^{1/2}\mu\normxt{N}{\phi^0-\pi}+\delta_{q+1}^{1/2}\left(\frac{\delta_{q+1}^{1/2}\mu\lambda_{q+1}}{\eta}\right)^N, \\
\left[\sin(\theta_{jm}\compt \phi^0)\right]_{C^N_{x,t}(G_{jm})}&\lesssim \lambda_{q+1}^N+\normxtjm{N}{\theta_{jm}\compt\phi^0-\theta_{jm}} 
\end{align*}
and analogous bounds for $\nabla a_{jm}\compt\phi^0$ and $\sin(\theta_{jm}\compt \phi^0)$. Substituting this into (\ref{aux estimates phi0 1}) and (\ref{aux estimates phi0 2}) leads to
\begin{align*}
	\normxt{N}{\phi^0-\pi}&\lesssim\frac{\delta_{q+1}^{1/2}\lambda_{q+1}^N}{\eta}+\frac{\delta_{q+1}^{1/2}\mu}{\eta}\left(\normxt{N}{\phi^0-\pi}+\mu^{-1}\max_{j,m}\normxtjm{N}{\theta_{jm}\compt\phi^0-\theta_{jm}}\right), \\ \normxtjm{N}{\theta_{jm}\compt\phi^0-\theta_{jm}}&\lesssim\frac{\delta_{q+1}^{1/2}\mu\lambda_{q+1}}{\eta}+\frac{\delta_{q+1}^{1/2}\mu^2}{\eta}\left(\normxt{N}{\phi^0-\pi}+\mu^{-1}\max_{j,n}\norm{\theta_{jn}\compt\phi^0-\theta_{jn}}_{N,\hspace{0.5pt}jn}\right),
\end{align*}
where we have used (\ref{hierarchy}) to simplify some terms. Taking the maximum in $j\in\{1,\dots, 6\}$ and $m\in\Lambda$ and combining both bounds, we have
\begin{align*}
&\normxt{N}{\phi^0-\pi}+\frac{1}{\mu}\max_{j,m}\normxtjm{N}{\theta_{jm}\compt\phi^0-\theta_{jm}}\lesssim\frac{\delta_{q+1}^{1/2}\lambda_{q+1}^N}{\eta}+\\&\hspace{90pt}+\frac{\delta_{q+1}^{1/2}\mu}{\eta}\left(\normxt{N}{\phi^0-\pi}+\frac{1}{\mu}\max_{j,m}\normxtjm{N}{\theta_{jm}\compt\phi^0-\theta_{jm}}\right).
\end{align*}
Since the coefficient multiplying the term in parenthesis can be made arbitrarily small by increasing $a>1$, we conclude (\ref{estimates phi0}) for $N\geq1$.
\end{proof}

The diffeomorphism $\phi^0_t$ will not be volume-preserving, in general. Fortunately, it is very nearly so, which means that it merely needs a very small correction. We will now prove the general strategy for constructing a correcting diffeomorphism and we will later apply it to our case.
\begin{proof}[Proof of \cref{lema dacorogna}]
Let us define a map $\Phi\in C^\infty(\RR^3\times[0,T]\times[0,1],\RR^3)$ as the unique solution of the following initial value problem:
\begin{equation}
\label{edo phi}
\begin{cases}
	\frac{\partial}{\partial s}\Phi(x,t,s)=\frac{u(\Phi(x,t,s),t)}{1+(1-s)f(\Phi(x,t,s),t)}, \\ \Phi(x,t,0)=x.
\end{cases}
\end{equation}
Due to (\ref{assumption dacorogna 1}), the right-hand side of the ODE is bounded by $2\normxt{0}{u}$ for all $s\in[0,1]$. Hence, the solution for a fixed $(x,t)\in \RR^3\times[0,T]$ is, indeed, defined in the considered interval. We see that $\Phi(\cdot,s)=\pi$ outside of the support of $u$ and, hence, on $(\RR^3\backslash \Omega)\times[0,T]$. In addition, it can be proved (see \cite[proof of Lemma 3]{DM}) that 
\[\det[D_x\Phi(x,t,1)]=1+f(x,t),\]
so it suffices to define $\phi(x,t)=\Phi(x,t,1)$. Therefore, all we have to do is to prove the appropriate bounds for $\phi$. Note that the ODE can be written as
\[\begin{cases}
		\frac{\partial}{\partial s}[\Phi(x,t,s)-x]=\frac{u(\Phi(x,t,s),t)}{1+(1-s)f(\Phi(x,t,s),t)}, \\ \Phi(x,t,0)-x=0,
	\end{cases}\]
so estimate (\ref{estimate phi dacorogna C0}) follows immediately. Concerning the higher order norms, it follows from (\ref{estimate product}) and (\ref{estimate quotient}) that
\[\frac{d}{ds}\normxt{N}{\Phi(\cdot,s)-\pi}\lesssim \normxt{N}{u\compt\Phi(\cdot,s)}+\normxt{0}{u}\normxt{N}{f\compt\Phi(\cdot,s)},\]
where we have used that $\normxt{0}{f}\leq 1/2$. Given $N\geq 1$, we use (\ref{estimate composition}) to obtain
\begin{align*}
\seminormxt{N}{u\compt\Phi(\cdot,s)}&\lesssim \seminormxt{1}{u}\left(1+\normxt{N-1}{\partial_t\Phi(\cdot,s)}+\normxt{N-1}{D\Phi(\cdot,s)}\right)\\ &\hspace{12pt}+\normxt{N}{u}\left(1+\seminormxt{1}{\Phi(\cdot,s)}^N\right) \\ &\lesssim \normxt{1}{u}\left(1+\normxt{N}{\Phi(\cdot,s)-\pi}\right)+\normxt{N}{u}\left(1+\normxt{1}{\Phi(\cdot,s)-\pi}^N\right).
\end{align*} 
Similarly,
\[\seminormxt{N}{f\compt\Phi(\cdot,s)}\lesssim \normxt{1}{f}\left(1+\normxt{N}{\Phi(\cdot,s)-\pi}\right)+\normxt{N}{f}\left(1+\normxt{1}{\Phi(\cdot,s)-\pi}^N\right).\]
Therefore, we conclude
\begin{align*}
\frac{d}{ds}\normxt{N}{\Phi(\cdot,s)-\pi}&\lesssim\left(\normxt{N}{u}+\normxt{0}{u}\normxt{N}{f}\right)\left(1+\normxt{1}{\Phi(\cdot,s)-\pi}^N\right) \\&\hspace{12pt}+\left(\normxt{1}{u}+\normxt{0}{u}\normxt{1}{f}\right)\left(1+\normxt{N}{\Phi(\cdot,s)-\pi}\right).
\end{align*}
In the particular case $N=1$, we have
\begin{align*}
	\frac{d}{ds}\normxt{1}{\Phi(\cdot,s)-\pi}&\lesssim \left(\normxt{1}{u}+\normxt{0}{u}\normxt{1}{f}\right)+\left(\normxt{1}{u}+\normxt{0}{u}\normxt{1}{f}\right)\normxt{1}{\Phi(\cdot,s)-\pi} \\
	&\lesssim \left(\normxt{1}{u}+\normxt{0}{u}\normxt{1}{f}\right)+\normxt{1}{\Phi(\cdot,s)-\pi},
\end{align*}
where we have used (\ref{assumption dacorogna 2}). Taking into account that $\Phi(\cdot,0)=\pi$ and applying Grönwall's inequality, we conclude that for any $s\in[0,1]$ we have
\[\normxt{1}{\Phi(\cdot,s)-\pi}\lesssim  \normxt{1}{u}+\normxt{0}{u}\normxt{1}{f}.\]
Substituting $s=1$ yields (\ref{estimate phi dacorogna CN}) for $N=1$. In particular, we have
\begin{equation}\label{eq.otra}
	\normxt{1}{\Phi(\cdot,s)} \lesssim 1 \qquad \forall s\in[0,1],  
\end{equation}
due to (\ref{assumption dacorogna 2}). This, along with (\ref{assumption dacorogna 2}), allows us to simplify the inequality for the general case $N\geq 1$:
\[\frac{d}{ds}\normxt{N}{\Phi(\cdot,s)-\pi}\lesssim\left(\normxt{N}{u}+\normxt{0}{u}\normxt{N}{f}\right)+\normxt{N}{\Phi(\cdot,s)-\pi}.\]
Hence, applying Grönwall's inequality yields (\ref{estimate phi dacorogna CN}) for $N\geq 1$.
\end{proof}

We are now ready to construct and estimate the correction $\phi^c$:
\begin{lemma}
There exists a map $\phi^c\in C^\infty(\RR^3\times[0,T],\RR^3)$ such that $\phi^c_t:\RR^3\to\RR^3$ is a diffeomorphism for any $t\in[0,T]$ and such that
\[\det[D(\phi^c_t\circ\phi^0_t)]=1.\]
In addition, $\phi^c \equiv \pi$ on $(\RR^3\backslash\Omega_{q}^{(2)})\times[0,T]$ and it satisfies the following estimates:
\begin{equation}
	\normxt{N}{\phi^c-\pi}\lesssim \delta_{q+1}\mu^2\eta^{-2}\max\{\mu^N,\lambda_{q+1}^{N-1+\tau}\} \qquad \forall N\geq 0. \label{estimates phic}
\end{equation}
Let us denote by $\psi^c\in C^\infty(\RR^3\times[0,T],\RR^3)$ the map such that $\psi^c_t$ is the inverse of $\phi^c_t$ for any  time $t\in[0,T]$. Then, it satisfies
\begin{equation}
	\normxt{N}{\psi^c-\pi}\lesssim \delta_{q+1}\mu^2\eta^{-2}\max\{\mu^N,\lambda_{q+1}^{N-1+\tau}\} \qquad \forall N\geq 0. \label{estimates psic}
\end{equation}
In particular, by (\ref{param 3}) and (\ref{param 5}) we have
\[\normxt{2}{\phi^c-\pi}+\normxt{2}{\psi^c-\pi}\lesssim \delta_{q+2}\lambda_{q+1}^{-11\tau},\]
which can be made arbitrarily small by increasing $a>1$.
\end{lemma}
\begin{proof}
By \cref{lema fc}, the density $f_c$ and the field $u_c$ satisfy the hypotheses of \cref{lema dacorogna} because
\[\normxt{0}{f_c}+ \normxt{1}{u_c}+\normxt{0}{u}\normxt{1}{f}\lesssim\frac{\delta_{q+1}\mu^3}{\eta^2}+\frac{\delta_{q+1}^2\mu^4\lambda_{q+1}}{\eta^4},\]
so they can be made arbitrarily small by taking $a>1$ sufficiently large, due to (\ref{def mu}) and (\ref{def eta}). Therefore, there exists a map $\phi^c\in C^\infty(\RR^3\times[0,T],\RR^3)$ such that $\phi^c_t$ is a diffeomorphism for all $t\in[0,T]$ and
\[\det(D\phi^c_t)=1+f_c=\det(D\psi^0_t).\]
Since $\psi^0_t$ is the inverse of $\phi^0_t$, it is easy to see that
\[\det[D(\phi^c_t\circ \phi^0_t)]=1.\]

Regarding the estimates, substituting bounds (\ref{estimates fc CN}) and (\ref{estimates uc CN}) into (\ref{estimate phi dacorogna C0}) and (\ref{estimate phi dacorogna CN}) leads to
\[\normxt{N}{\phi^c-\pi}\lesssim \frac{\delta_{q+1}\mu^2}{\eta^2}\left(\mu^N+\lambda_{q+1}^{N-1+\tau}+\frac{\delta_{q+1}\mu^2\lambda_{q+1}^N}{\eta^2}\right),\]
which yields (\ref{estimates phic}) because of (\ref{param 5}). Hence, it suffices to estimate the map $\psi^c$. Since $\phi^c_t$ is a diffeomorphism of $\RR^3$, we have
\[\normxt{0}{\psi^c-\pi}=\normxt{0}{\psi^c\compt \phi^c-\pi\compt\phi^c}=\normxt{0}{\pi-\phi^c},\]
so (\ref{estimates psic}) for $N=0$ follows from (\ref{estimates phic}) with $N=0$. To estimate the higher-order norms in a more systematic manner, let us define diffeomorphisms $\Phi, \Psi$ of $\RR^3\times[0,T]$ given by
\[ \Phi(x,t)\coloneqq (\phi^c(x,t),t), \qquad \Psi(x,t)\coloneqq (\psi^c(x,t),t).\]
We see that $\Psi$ is the inverse of $\Phi$. It also follows from (\ref{estimates phic}) that
\begin{align}
	\normxt{1}{\Phi-\Id}&\lesssim \delta_{q+1}\mu^3\eta^{-2}, \label{aux estimates Phi C1} \\
	\normxt{N}{\Phi-\Id}&\lesssim \delta_{q+1}\mu^2\eta^{-2}\lambda_{q+1}^{N-1+\tau} \qquad \forall N\geq 2 \label{aux estimates Phi}
\end{align}
where we have used $\mu^2\leq \lambda_{q+1}$, which is clear from definiton (\ref{def mu}). Regarding the $C^1$-norm, differentiating the identity $\Phi\circ\Psi=\Id$ leads to
\[D_{x,t}\Psi=(D_{x,t}\Phi)^{-1}\circ \Psi,\]
where we have written $D_{x,t}$ to denote the differential on $\RR^3\times[0,T]$. We thus have
\begin{align*}
\normxt{0}{D_{x,t}\Psi-\Id}&=\normxt{0}{(D_{x,t}\Phi)^{-1}-\Id}\lesssim \normxt{0}{(D_{x,t}\Phi)^{-1}}\normxt{0}{\Id-D_{x,t}\Phi} \\ &\lesssim \normxt{0}{(D_{x,t}\Phi)^{-1}-\Id}\normxt{0}{\Id-D_{x,t}\Phi}+\normxt{0}{\Id-D_{x,t}\Phi}.
\end{align*}
If follows from (\ref{param 3}) and (\ref{aux estimates Phi C1}) that $\normxt{0}{\Id-D_{x,t}\Phi}$ can be made smaller than $1/2$ by taking $a>1$ sufficiently large. Therefore,
\[\normxt{0}{D_{x,t}\Psi-\Id}=\normxt{0}{(D_{x,t}\Phi)^{-1}-\Id}\lesssim \normxt{0}{\Id-D_{x,t}\Phi}\lesssim \delta_{q+1}\mu^3\eta^{-2},\]
which, combined with (\ref{estimates psic}) with $N=0$, yields the desired bound for $N=1$. In particular, by (\ref{param 3}), we have
\begin{equation}
\normxt{1}{\Phi}+\normxt{1}{\Psi}\lesssim 1. \label{aux segunda derivada psic 1}
\end{equation}
Let us fix $N\geq 2$. Due to the previous estimate, it follows from (\ref{estimate inverse}) that
\begin{equation}
\normxt{N}{\Psi}\lesssim \normxt{N}{\Phi}\lesssim 1+\normxt{N}{\Phi-\Id}. \label{aux segunda derivada psic 2}
\end{equation}
Hence, it follows from (\ref{estimate composition}) that
\begin{align*}
\normxt{N}{\Psi-\Id}&\lesssim \normxt{N}{(\Id-\Phi)\circ\Psi}\lesssim \normxt{0}{\Id-\Phi}+\normxt{N}{\Id-\Phi}\normxt{1}{\Psi}^N+\normxt{1}{\Id-\Phi}\normxt{N}{\Psi} \\ &\lesssim \normxt{N}{\Id-\Phi}\stackrel{(\ref{aux estimates Phi})}{\lesssim} \delta_{q+1}\mu^2\eta^{-2}\lambda_{q+1}^{N-1+\tau},
\end{align*}
where we have used (\ref{aux segunda derivada psic 1}) and (\ref{aux segunda derivada psic 2}). From this, we deduce (\ref{estimates psic}) for $N\geq 2$.
\end{proof}

We define the final map as $\phi^{q+1}\coloneqq \phi^c\compt \phi^0$. By the previous lemma, for each fixed time $t\in[0,T]$, we have a volume-preserving diffeomorphism $\phi^{q+1}_t= \phi^c_t\circ \phi^0_t$ of $\RR^3$ that is the identity outside of $\Omega_{q,2}$. As we have seen, $\phi^c$ is very close to the identity in $C^2$, so it will produce small corrections to the vector fields. However, it is essential to ensure that they are divergence-free.

\section{Estimates on the vector fields} \label{section estimates vector fields}
In this section we will estimate the perturbation to the magnetic field $b_{q+1}=B_{q+1}-B_q$ and to the velocity $w_{q+1}=v_{q+1}-v_q$. We begin by proving that $b_{q+1}$ is very small. In fact, we will see in \cref{section matrix} that it is so small that it can be essentially ignored when dealing with the Reynolds stress.
\begin{lemma}
The perturbation $b_{q+1}$ satisfies
\begin{equation}
	\label{estimates b}
	\normxt{0}{b_{q+1}}+\lambda_{q+1}^{-1}\normxt{1}{b_{q+1}}\lesssim \delta_{q+2}\lambda_{q+1}^{-11\tau}.
\end{equation}
\end{lemma}
\begin{proof}
Let $B^0_{q+1}$ be the pushforward of $B_q$ by $\phi^0_t$, that is, 
\[B^0_{q+1}\coloneqq (D\phi^0\,B_q)\compt\psi^0_t,\]
where we have used that $\psi^0_t$ is the inverse of $\phi^0_t$. Since $\phi^c$ is very close to the identity in $C^1$, we expect this to be the dominating term in $B_{q+1}$. By standard properties of the vector product, we have
\begin{equation}
(\nabla a_{jm}\compt\phi^0)\times (\zeta_j\times k_{jm})=(\zeta_j\otimes k_{jm}-k_{jm}\otimes \zeta_j)(\nabla a_{jm}\compt\phi^0)
\label{desarrollar producto nabla am}
\end{equation}
which leads to
\[D[(\nabla a_{jm}\compt\phi^0)\times (\zeta_j\times k_{jm})]\,B_q=(\zeta_j\otimes k_{jm}-k_{jm}\otimes \zeta_j)\hspace{0.5pt}(D^2 a_{jm}\compt \phi^0)\hspace{0.5pt}D\phi^0B_q.\]
Therefore, by \cref{identidad inversa}, we have
\begin{align*}
	D\phi^0B_{q}=B_{q}&+\sum_{j,m}\frac{\lambda_{q+1}a_{jm}\compt \phi^0}{\eta}[k_{jm}\cdot(B_{q}-B^q_{jm})]\,\zeta_j\cos(\theta_{jm}\compt\phi^0) \\&+\sum_{j,m}\frac{a_{jm}\compt \phi^0}{\ell_m\eta}B_q\cdot\nabla(\theta_{jm}\compt \phi^0-\theta_{jm})\,\zeta_j \cos(\theta_{jm}\compt\phi^0) \\&+\sum_{j,m}\frac{1}{\ell_m\eta}[(\nabla a_{jm}\compt \phi^0)D\phi^0B_q]\,\zeta_j\sin(\theta_{jm}\compt\phi^0)\\&-\sum_{j,m}\frac{(\zeta_j\otimes k_{jm}-k_{jm}\otimes \zeta_j)\hspace{0.5pt}(D^2 a_{jm}\compt \phi^0)\hspace{0.5pt}D\phi^0B_q}{\ell_m^2\eta\lambda_{q+1}}\cos(\theta_{jm}\compt \phi^0)\\&+\sum_{j,m} \frac{(\nabla a_{jm}\compt \phi^0)\times(\zeta_j\times k_{jm})}{\ell_m^2\eta\lambda_{q+1}}[B_q\cdot\nabla(\theta_{jm}\compt\phi^0)]\sin(\theta_{jm}\compt\phi^0),
\end{align*}
where we have used that $k_{jm}\cdot B^q_{jm}=0$ for all $j\in\{1,\dots 6\}$ and $m\in \Lambda$. Since $\psi^0_t$ is the inverse of $\phi^0_t$, we obtain
\begin{align}
	B^0_{q+1}=B_{q}\compt \psi^0&+\sum_{j,m}\frac{\lambda_{q+1}a_{jm}}{\eta}[k_{jm}\cdot(B_{q}\compt\psi^0-B^q_{jm})]\,\zeta_j\cos(\theta_{jm})\nonumber\\&+\sum_{j,m}\frac{a_{jm}}{\ell_m\eta}\left[(B_j\cdot\nabla(\theta_{jm}\compt\phi^0-\theta_{jm}))\compt\psi^0\right]\zeta_j\cos(\theta_{jm})  \nonumber\\&+\sum_{j,m}\frac{1}{\ell_m\eta}(\nabla a_{jm}\cdot  B^0_{q+1})\,\zeta_j\sin(\theta_{jm}) \label{expresión tilde Bq+1}\\&-\sum_{j,m}\frac{(\zeta_j\otimes k_{jm}-k_{jm}\otimes \zeta_j)D^2 a_{jm}B^0_{q+1}}{\ell_m^2\eta\lambda_{q+1}}\cos(\theta_{jm}) \nonumber\\&+\sum_{j,m}\frac{\nabla a_{jm}\times (\zeta_j\times k_{jm})}{\ell_m^2\eta\lambda_{q+1}}\left[(B_q\cdot\nabla(\theta_{jm}\compt\phi^0))\compt\psi^0\right]\sin(\theta_{jm}). \nonumber
\end{align}
Given $j\in\{1,\dots 6\}$ and $m\in\Lambda$, the support of $a_{jm}$ with  is contained in a cube centered at $\mu^{-1}m$ and whose sidelength is comparable to $\mu^{-1}$. Hence,
\begin{equation}
\begin{split}
\normxtjm{0}{B_q\compt \psi^0-B^q_{Jm}}&\leq \normxtjm{0}{B_q\compt \psi^0-B_q}+\normxtjm{0}{B_q(\,\cdot\,)-B_q(\mu^{-1}m)}\\&\lesssim \normxt{1}{B_q}\normxt{0}{\psi^0-\pi}+\mu^{-1}\normxt{1}{B_q}\lesssim \mu^{-1}\normxt{1}{B_q},
\end{split}
\label{diferencia aproximación B}
\end{equation}
where we have used (\ref{hierarchy}) and (\ref{estimates psi}), as well as the definition of $B^q_{jm}$. Meanwhile, by (\ref{inductive tamaño Bq}) and (\ref{estimates phi0}), we have
\begin{align}
	\normxtjm{0}{B_q\cdot\nabla(\theta_{jm}\compt\phi^0-\theta_{jm})}&\lesssim \delta_{q+1}^{1/2}\mu\eta^{-1}\lambda_{q+1}, \label{aux fases tilde B 1} \\ \normxtjm{0}{B_q\cdot\nabla(\theta_{jm}\compt\phi^0)}&\lesssim \lambda_{q+1}. \label{aux fases tilde B 2}
\end{align}
These estimates, along with (\ref{estimates am}), lead to
\begin{align*}
	\normxt{0}{B^0_{q+1}-B_q}&\lesssim \normxt{1}{B_q}\normxt{0}{\psi^0-\pi}+  \frac{\delta_{q+1}^{1/2}\lambda_{q+1}}{\mu\eta}\normxt{1}{B_q}+\frac{\delta_{q+1}\mu\lambda_{q+1}}{\eta^2}\\&\hspace{12pt}+\frac{\delta_{q+1}^{1/2}\mu}{\eta}\normxt{0}{B^0_{q+1}}+\frac{\delta_{q+1}^{1/2}\mu^2}{\eta\lambda_{q+1}}\normxt{0}{B^0_{q+1}}+\frac{\delta_{q+1}^{1/2}\mu}{\eta}  \\ &\lesssim \frac{\delta_{q+1}^{1/2}\lambda_{q+1}}{\mu\eta}\normxt{1}{B_q}+\frac{\delta_{q+1}\mu\lambda_{q+1}}{\eta^2}+\frac{\delta_{q+1}^{1/2}\mu}{\eta}\normxt{0}{B^0_{q+1}}.
\end{align*}
Since the term multiplying the last norm can be made arbitrarily small by increasing $a>1$, by (\ref{hierarchy}), and since $\normxt{0}{B_q}\lesssim 1$, we have
\[\normxt{0}{B^0_{q+1}-B_q}\lesssim \frac{\delta_{q+1}^{1/2}\lambda_{q+1}}{\mu\eta}\normxt{1}{B_q}+\frac{\delta_{q+1}\mu\lambda_{q+1}}{\eta^2} \lesssim \delta_{q+2}\lambda_{q+1}^{-11\tau}\]
where we have used (\ref{inductive cota Bq C1}), (\ref{param 2}) and (\ref{param 5}). In particular, by (\ref{inductive tamaño Bq}) we have
\begin{equation}
\normxt{0}{B^0_{q+1}}\lesssim 1.
\label{tamaño tilde Bq+1}
\end{equation}

Next, we will focus on the $C^1$-norm. First of all, it follows from (\ref{inductive cota Bq C1}) and (\ref{estimates psi}) that the first term in (\ref{expresión tilde Bq+1}) can be estimated as
\begin{equation}
\normxt{1}{B_q\compt\psi^0}\lesssim \normxt{1}{B_q}\left(1+\normxt{1}{\psi^0}\right)\lesssim \delta_q^{1/2}\delta_{q+1}^{1/2}\lambda_q\eta^{-1}\lambda_{q+1}.
\label{cota BJ circ psi0}
\end{equation}
Concerning the second term, it follows from (\ref{inductive tamaño Bq}), (\ref{estimates psi}) and (\ref{diferencia aproximación B}) that
\begin{align*}
\normxt{1}{a_{jm}(B_q\compt \psi^0-B^q_{jm})}&\lesssim \normxt{1}{a_{jm}}\mu^{-1}\normxt{1}{B_q}+\normxt{0}{a_{jm}}\normxt{1}{B_q}\left(1+\normxt{1}{\psi^0}\right)\\&\lesssim \delta_{q+1}^{1/2}(\delta_q^{1/2}\lambda_q) +\delta_{q+1}^{1/2}(\delta_q^{1/2}\lambda_q)(\delta_{q+1}^{1/2}\eta^{-1}\lambda_{q+1}) \\&\lesssim \delta_q^{1/2}\delta_{q+1}\lambda_q\eta^{-1}\lambda_{q+1}.
\end{align*}
Therefore, up to a constant, the $C^1$-norm of the second term in (\ref{expresión tilde Bq+1}) can be bounded by
\begin{align*}
&\max_{j,m}\left(\frac{\lambda_{q+1}^2}{\eta}\normxt{0}{a_{jm}(B_q\compt \psi^0-B^q_{jm})}+\frac{\lambda_{q+1}}{\eta}\normxt{1}{a_{jm}(B_q\compt \psi^0-B^q_{jm})}\right)\\ &\hspace{60pt}\lesssim \delta_q^{1/2}\delta_{q+1}^{1/2}\lambda_q\mu^{-1}\eta^{-1}\lambda_{q+1}^2+\delta_q^{1/2}\delta_{q+1}\lambda_q\eta^{-2}\lambda_{q+1}^2 \\  &\hspace{60pt}\lesssim \delta_q^{1/2}\delta_{q+1}^{1/2}\lambda_q\mu^{-1}\eta^{-1}\lambda_{q+1}^2,
\end{align*}
where we have also used (\ref{inductive cota Bq C1}) and (\ref{diferencia aproximación B}). Regarding the third term in (\ref{expresión tilde Bq+1}), by (\ref{inductive tamaño Bq}), (\ref{inductive cota Bq C1}), (\ref{estimates psi}) and (\ref{estimates phi0}), we have
\begin{align}
	\normxtjm{1}{[B_q\cdot\nabla(\theta_{jm}\compt\phi^0-\theta_{jm})]\compt \psi^0}&\lesssim \normxt{1}{B_q}\normxtjm{1}{\theta_{jm}\compt\phi^0-\theta_{jm}}\left(1+\normxt{1}{\psi^0-\pi}\right) \nonumber\\&\hspace{12pt}+\normxt{0}{B_q}\normxtjm{2}{\theta_{jm}\compt\phi^0-\theta_{jm}}\left(1+\normxt{1}{\psi^0-\pi}\right) \label{aux BJ dtheta} \\&\lesssim\delta_{q+1}\mu\eta^{-2}\lambda_{q+1}^3, \nonumber
\end{align}
Hence, up to a constant, the $C^1$-norm of the third term in (\ref{expresión tilde Bq+1}) can be bounded by
\begin{align*}
&\max_{j,m}\bigg(\delta_{q+1}^{1/2}\eta^{-1}\lambda_{q+1}\normxtjm{0}{[B_q\cdot\nabla(\theta_{jm}\compt\phi^0-\theta_{jm})]\compt \psi^0}\\&\hspace{30pt}+\delta_{q+1}^{1/2}\eta^{-1}\normxtjm{1}{[B_q\cdot\nabla(\theta_{jm}\compt\phi^0-\theta_{jm})]\compt \psi^0}\bigg)\lesssim\delta_{q+1}^{3/2}\mu\eta^{-3}\lambda_{q+1}^3,
\end{align*}
where we have also used (\ref{def eta}) and (\ref{aux fases tilde B 1}). Next, since $\nabla\theta_{jm}$ is a constant of order $\lambda_{q+1}$, from (\ref{cota BJ circ psi0}) and (\ref{aux BJ dtheta}) we deduce
\begin{align*}
\normxtjm{1}{[B_q\cdot\nabla(\theta_{jm}\compt\phi^0)]\compt\psi^0}&\lesssim \delta_{q+1}\mu\eta^{-2}\lambda_{q+1}^3+\delta_q^{1/2}\delta_{q+1}^{1/2}\lambda_q\eta^{-1}\lambda_{q+1}^2\\&\lesssim \delta_{q+1}\mu\eta^{-2}\lambda_{q+1}^3.
\end{align*}
Using this estimate, along with (\ref{estimates am}) and (\ref{aux fases tilde B 2}), we see that the $C^1$-norm of the last term in (\ref{expresión tilde Bq+1}) is bounded by
\[\delta_{q+1}^{1/2}\mu\eta^{-1}\lambda_{q+1}+(\delta_{q+1}^{1/2}\mu\eta^{-1}\lambda_{q+1}^{-1})(\delta_{q+1}\mu\eta^{-2}\lambda_{q+1}^3)\stackrel{(\ref{param 1})}{\lesssim} \delta_{q+1}^{1/2}\mu\eta^{-1}\lambda_{q+1}.\]
Using (\ref{estimates am}) and (\ref{tamaño tilde Bq+1}) to estimate the rest of the terms in (\ref{expresión tilde Bq+1}) and simplifying the expressions using $\delta_{q+1}^{1/2}\lambda_{q+1}\gg \eta$ leads to
\[\normxt{1}{B^0_{q+1}}\lesssim \frac{\delta_q^{1/2}\delta_{q+1}^{1/2}\lambda_q\lambda_{q+1}^2}{\mu\eta}+\frac{\delta_{q+1}^{3/2}\mu\lambda_{q+1}^3}{\eta^3}+\frac{\delta_{q+1}^{1/2}\mu}{\eta}\normxt{1}{B^0_{q+1}}.\]
Since the term multiplying the last norm can be made arbitrarily small by taking $a>1$ sufficiently large, due to (\ref{hierarchy}), we conclude
\[\normxt{1}{B^0_{q+1}}\lesssim \frac{\delta_q^{1/2}\delta_{q+1}^{1/2}\lambda_q\lambda_{q+1}^2}{\mu\eta}+\frac{\delta_{q+1}^{3/2}\mu\lambda_{q+1}^3}{\eta^3}\lesssim \delta_{q+2}\lambda_{q+1}^{1-11\tau},\]
where we have used (\ref{param 2}) and (\ref{param 4}). 

Once we have suitable bounds for $B^0_{q+1}$, it is easy to estimate $b_{q+1}$. First of all, note that
\begin{equation}
	B_{q+1}=(\phi^c_t)_\ast\hspace{0.5pt}B^0_{q+1}=(D\phi^cB^0_{q+1})\compt\psi^c
	\label{aux Bq+1}
\end{equation}
because $\phi^{q+1}=\phi^c\compt\phi^0$. Hence,
\begin{align*}
	\normxt{0}{b_{q+1}}&=\normxt{0}{(D\phi^cB^0_{q+1})\compt\psi^c-B_q}\\&\leq \normxt{0}{(D\phi^cB^0_{q+1})\compt\psi^c-B^0_{q+1}\compt\psi^c}+\normxt{0}{B^0_{q+1}\compt\psi^c-B^0_{q+1}}+\normxt{0}{B^0_{q+1}-B_q} \\&\leq \normxt{1}{\phi^c-\pi}\normxt{0}{B^0_{q+1}}+\normxt{1}{B^0_{q+1}}\normxt{0}{\psi^c-\pi}+\normxt{0}{B^0_{q+1}-B_q} \\ &\lesssim \delta_{q+1}\mu^3\eta^{-2}+(\delta_{q+2}\lambda_{q+1}^{1-11\tau})(\delta_{q+1}\mu^2\eta^{-2})+\delta_{q+2}\lambda_{q+1}^{-11\tau}\lesssim \delta_{q+2}\lambda_{q+1}^{-11\tau}
\end{align*}
by (\ref{hierarchy}) and (\ref{param 5}). Regarding the $C^1$-norm, differentiating (\ref{aux Bq+1}) leads to
\begin{align*}
\normxt{1}{B_{q+1}}&\lesssim \normxt{1}{D\phi^c}\normxt{0}{B^0_{q+1}}\left(1+\normxt{1}{\psi^c-\pi}\right)+\normxt{0}{D\phi^c}\normxt{1}{B^0_{q+1}}\left(1+\normxt{1}{\psi^c-\pi}\right)\\&\lesssim\delta_{q+1}\mu^2\eta^{-2}\lambda_{q+1}^{1+\tau}+\delta_{q+2}\lambda_{q+1}^{1-11\tau}\lesssim \delta_{q+2}\lambda_{q+1}^{1-11\tau},
\end{align*}
where we have used (\ref{param 3}), (\ref{estimates phic}) and (\ref{estimates psic}). The claimed bound for the $C^1$-norm of $b_{q+1}$ follows from this and the inductive estimate (\ref{inductive cota Bq C1}) for $B_q$, due to (\ref{relación entre lambdas}).
\end{proof}

The previous lemma readily implies (\ref{cambio B q}) if $a>1$ is taken sufficiently large to compensate the numerical constants. In fact, we have obtained better behavior than what we claim in (\ref{cambio B q}), which would lead to better Hölder regularity of the final magnetic field. This is not very relevant, however,  because the improved Hölder exponent would still be very small, so we have decided to ignore it. On the other hand, we also have (\ref{inductive cota Bq C1}) for sufficiently large $a>1$. Regarding (\ref{inductive tamaño Bq}), we have
\[\delta_{q+2}\left(2^{-q}-2^{-q-1}\right)^{-1}\leq 2^{q+1}\delta_{q+2}=2^{q+1}a^{-2\beta b^{q+2}},\]
which can be made arbitrarily small by choosing $a>1$ sufficiently large (independent of $q$), so (\ref{inductive tamaño Bq}) holds, too. We conclude that the new magnetic field satisfies all of the required conditions.

We will now focus on the perturbation to the velocity field. We begin by estimating the main perturbation term, $w_0$, and the correction $w_c$, which were defined in (\ref{def w0}) and (\ref{def wc}), respectively. 
\begin{lemma}
	The correction $w_0$ satisfies
	\begin{equation}
		\normxt{N}{w_0}\leq  \frac{1}{2}\delta_{q+1}^{1/2}\lambda_{q+1}^N \qquad \forall N\leq 4. \label{cota w0 exacta}
	\end{equation}
	In addition, for any $N\geq0$ we have
	\begin{align}
		\normxt{N}{w_0}&\lesssim \delta_{q+1}^{1/2}\lambda_{q+1}^N,
		\label{cota w0 con constantes} \\
		\normxt{N}{w_c}&\lesssim \delta_{q+1}^{1/2}\mu\lambda_{q+1}^{N-1}. \label{cota wc}
	\end{align}
\end{lemma}
\begin{proof}
	By \cref{remark soporte ajm disjuntos}, any given point is in the support of at most 16 of the coefficients $a_{jm}$. Hence, it follows from (\ref{estimate am C0}) that $\norm{w_0}_0\leq 2^{-2}\delta_{q+1}^{1/2}$. Regarding the higher-order norms, it follows from (\ref{estimates am}) that the derivatives of $\cos(\theta_{jm})$ dominate because $\lambda_{q+1}\gg \mu$. Thus, we deduce (\ref{cota w0 con constantes}). An analogous argument applies to $w_c$, leading to (\ref{cota wc}).
	
	In addition, it follows from the definition of the parameters that $\mu^{-1}\lambda_{q+1}$ can be made arbitrarily small by taking $a>1$ sufficiently large. By (\ref{estimates am}), we thus see that in the expression of the $N$-th derivative of $w_0$, all of the other terms can be made arbitrarily small in comparison to the $N$-th derivative of the cosine. Combining this with \cref{remark soporte ajm disjuntos} and \eqref{estimate am C0}, we conclude (\ref{cota w0 exacta}).
\end{proof}

Note that the fields $w_0$ and $w_c$ are \textit{ad hoc} and it is not obvious a priori whether the remaining term $w_\phi=w_{q+1}-w_0-w_c$ is small, as we claim. To prove that this is the case, we need to further study the behavior of the map $\phi^{q+1}$. This is one of the most delicate parts of our construction. 

\begin{lemma} \label{lema wphi}
Given $N\in\{0,1,2,3,4\}$, the correction $w_\phi$ satisfies
\begin{equation}
	\label{estimates wphi}
	\normxt{N}{w_\phi}\lesssim \delta_{q+2}\lambda_{q+1}^{N-11\tau}.
\end{equation}
\end{lemma}
\begin{proof}
	Let us define
	\begin{align*}
		v^0_{q+1}&\coloneqq (\partial_t\phi^0+D\phi^0v_q)\compt\psi^0, \\
		w^0_\phi&\coloneqq v^0_{q+1}-v_q-w_0-w_c.
	\end{align*}
	Since $\phi^c$ and $\psi^c$ are very close to $\pi$ in $C^2$, we expect $v^0_{q+1}$ and $\widetilde{w}_\phi$ to be very close to $v_{q+1}$ and $w_\phi$, respectively. To simplify the expressions, we introduce the notation
	\[D_v\coloneqq \partial_t+v_q\cdot\nabla\,	.\]
	Note that, by (\ref{desarrollar producto nabla am}), we have
	\[D_v[(\nabla a_{jm}\compt\phi^0)\times (\zeta_j\times k_{jm})]=(\zeta_j\otimes k_{jm}-k_{jm}\otimes \zeta_j)\hspace{0.5pt}\left[(\nabla\partial_t a_{jm})\compt\phi^0+(D^2 a_{jm}\compt \phi^0)\hspace{0.5pt}D_v\phi^0\right].\]
	Hence, differentiating the expression in \cref{identidad inversa} leads to
	\begin{align*}
		D_v\phi^0=v_q&+\sum_{j,m}[D_v(\theta_{jm}\compt\phi^0)]\frac{a_{jm}\compt\phi^0}{\ell_m\eta}\,\zeta_j\cos(\theta_{jm}\compt\phi^0) \\&+\sum_{j,m} \frac{\partial_t a_{jm}\compt\phi^0+(\nabla a_{jm}\compt\phi^0)D_v\phi^0}{\ell_m\eta}\,\zeta_j\sin(\theta_{jm}\compt\phi^0) \\&+\sum_{j,m}[D_v(\theta_{jm}\compt\phi^0)]\frac{(\nabla a_{jm}\compt\phi^0)\times(\zeta_j\times k_{jm})}{\ell_m^2\eta\lambda_{q+1}}\sin(\theta_{jm}\compt\phi^0) \\&-\sum_{j,m}\frac{(\zeta_j\otimes k_{jm}-k_{jm}\otimes \zeta_j)\hspace{0.5pt}\left[(\nabla\partial_t a_{jm})\compt\phi^0+(D^2 a_{jm}\compt \phi^0)\hspace{0.5pt}D_v\phi^0\right]}{\ell_m^2\eta\lambda_{q+1}}\cos(\theta_{jm}\compt\phi^0).
	\end{align*}
	Note that due to the cancellation (\ref{derivada material thetam}), we have
	\[D_v(\theta_{jm}\compt\phi^0)=D_v(\theta_{jm}\compt\phi^0-\theta_{jm})+\ell_m\eta+\ell_m\lambda_{q+1}k_{jm}\cdot(v_q-v^q_{jm}).\]
	Therefore, we conclude
	\begin{align}
		v^0_{q+1}=v_q\compt\psi^0+w_0+w_c&+\sum_{j,m}\mu^{-1}[(D_v(\theta_{jm}\compt\phi^0-\theta_{jm}))\compt\psi^0]\,u_{jm} \nonumber \\
		\begin{split}
			&+\sum_{j,m}\ell_m\mu^{-1}\lambda_{q+1}k_{jm}\cdot(v_q\compt\psi^0-v^q_{jm})\,u_{jm} \\&+\sum_{j,m}\frac{\partial_t a_{jm}+v^0_{q+1}\cdot \nabla a_{jm}}{\ell_m\eta}\,\zeta_j\sin(\theta_{jm})
			\label{expresion tilde vq+1}
		\end{split}
		\\&-\sum_{j,m}\frac{(\zeta_j\otimes k_{jm}-k_{jm}\otimes \zeta_j)(\nabla\partial_t a_{jm}+D^2a_{jm}\hspace{0.5pt}v^0_{q+1})}{\ell_m^2\eta\lambda_{q+1}}\cos(\theta_{jm}), \nonumber
	\end{align}
	where
	\begin{equation}
		u_{jm}\coloneqq \mu\frac{ a_{jm}}{\ell_m\eta}\,\zeta_j\cos(\theta_{jm})+\mu\frac{\nabla a_{jm}\times(\zeta_j\times k_{jm})}{\ell_m^2\eta\lambda_{q+1}}\sin(\theta_{jm}).
		\label{def um}
	\end{equation}
	Note that the support of $u_{jm}$ is contained in the support of $a_{jm}$. By \eqref{estimates psi}, the image of this set is contained in $G_{jm}$ for sufficiently large $a>1$. Hence, it follows from (\ref{diferencia thetam con y sin phi0}) that
	\begin{align}
		[D_v(\theta_{jm}\compt\phi^0-\theta_{jm})]\compt\psi^0\big|_{\supp a_{jm}}&=\sum_{n\in\Lambda}[(D_v(\theta_{jn}\compt\phi^0-\theta_{jn}))\compt\psi^0]\frac{\ell_mk_{jm}\cdot[\nabla a_{jn}\times(\zeta_j\times k_{jn})]}{\ell_n^2\eta}\sin(\theta_{jn}) \nonumber \\ \begin{split} \label{aux Dv diferencia thetas}
			&\hspace{12pt}+\sum_{n\in\Lambda} \ell_n^{-1}\ell_mk_{jm}\cdot[\nabla a_{jn}\times(\zeta_j\times k_{jn})]\sin(\theta_{jn}) \\ &\hspace{12pt}+\sum_{n\in\Lambda}  \lambda_{q+1}k_{jn}\cdot(v_q\compt\psi^0-v^q_{jn})\,\frac{\ell_mk_{jm}\cdot[\nabla a_{jn}\times(\zeta_j\times k_{jn})]}{\ell_n\eta}\sin(\theta_{jn})
		\end{split}  \\ &\hspace{12pt}-\sum_{n\in\Lambda}\frac{\ell_m k_{jm}\cdot(\zeta_j\otimes k_{jn}-k_{jn}\otimes \zeta_j)(\nabla\partial_t a_{jn}+D^2a_{jn}v^0_{q+1})}{\ell_n^2\eta}\cos(\theta_{jn}) \nonumber.
	\end{align}
	Let 
	\[\Theta_{jm}\coloneqq \mu^{-1}[D_v(\theta_{jm}\compt\phi^0-\theta_{jm})]\compt\psi^0\big|_{\supp a_{jm}}\]
	and consider the coefficients
	\begin{align*}
		f_w&\coloneqq v_q\compt\psi^0-v_q+\sum_{j,m}\ell_m\mu^{-1}\lambda_{q+1}k_{jm}\cdot(v_q\compt\psi^0-v^q_{jm})\,u_{jm} \\&\hspace{12pt}+\sum_{j,m}\frac{D_v\hspace{0.5pt} a_{jm}+(w_0+w_c)\cdot\nabla a_{jm}}{\ell_m\eta}\,\zeta_j\sin(\theta_{jm})\\&\hspace{12pt}-\sum_{j,m}\frac{(\zeta_j\otimes k_{jm}-k_{jm}\times \zeta_j)[D_v\nabla a_{jm}+D^2a_{jm}(w_0+w_c)]}{\ell_m^2\eta\lambda_{q+1}}\cos(\theta_{jm}), \\[8pt]
		 M_w&\coloneqq \sum_{j,m}\left[\frac{\zeta_j\otimes \nabla a_{jm}}{\ell_m\eta}\sin(\theta_{jm})+\frac{(\zeta_j\otimes k_{jm}-k_{jm}\otimes \zeta_j)D^2a_{jm}}{\ell_m^2\eta\lambda_{q+1}}\cos(\theta_{jm})\right].
	\end{align*}
	Then, we can write
	\begin{equation}
	w^0_\phi=f_w+M_ww^0_\phi+\sum_{j,m}u_{jm}\Theta_{jm}.
	\label{fórmula compacta tildew}
	\end{equation}
	
    Similarly, if we define
	\begin{align*}
		f_{jm}&\coloneqq \sum_{n\in\Lambda} \ell_n^{-1}\ell_m\mu^{-1}k_{jm}\cdot[\nabla a_{jn}\times(\zeta_j\times k_{jn})]\sin(\theta_{jn}) \\ &\hspace{12pt}+\sum_{n\in\Lambda}  \lambda_{q+1}k_{jn}\cdot(v_q\compt\psi^0-v^q_{jn})\,\frac{\ell_mk_{jm}\cdot[\nabla a_{jn}\times(\zeta_j\times k_{jn})]}{\ell_n\mu\eta}\sin(\theta_{jn})
		\\ &\hspace{12pt}-\sum_{n\in\Lambda}\frac{\ell_m k_{jm}\cdot(\zeta_j\otimes k_{jn}-k_{jn}\otimes \zeta_j)[D_v\nabla a_{jn}+D^2a_{jn}(w_0+w_c)]}{\ell_n^2\mu\eta}\cos(\theta_{jn}) \\
		M_{jm}&\coloneqq -\sum_{n\in\Lambda}\frac{\ell_m k_{jm}\cdot(\zeta_j\otimes k_{jn}-k_{jn}\otimes \zeta_j)D^2a_{jn}}{\ell_n^2\mu\eta}\cos(\theta_{jn}), \\[6pt]
		c_{jmn}&\coloneqq \frac{\ell_mk_{jm}\cdot[\nabla a_{jn}\times(\zeta_j\times k_{jn})]}{\ell_n^2\eta}\sin(\theta_{jn}).
	\end{align*}
	we can then write
	\[\Theta_{jm}=f_{jm}+M_{jm}w^0_\phi+\sum_{n\in\Lambda}c_{jmn}\Theta_{jn}.\]
	If we define a vector $\Theta\coloneqq (\Theta_{jm})_{j,m}$ and we consider vector versions of the other parameters, we can write this as
	\begin{equation*}
		\begin{cases}
			w^0_\phi= f_w+M_ww^0_\phi+U\Theta, \\
			\hspace{4.5pt}\Theta=f_\Theta+M_\Theta w^0_\phi+C\Theta.
		\end{cases} 
	\end{equation*}
Furthermore, if we define
\[W\coloneqq\begin{pmatrix}
	w^0_{\phi} \\ \Theta
\end{pmatrix}, \qquad F\coloneqq\begin{pmatrix}
	f_w \\ f_\Theta
\end{pmatrix}, \qquad M\coloneqq\begin{pmatrix}
	M_w & U \\ M_\Theta & C
\end{pmatrix},\]
we can write the previous system as
\begin{equation}
W=F+M W.
\label{sistema W}
\end{equation}

Let us estimate each term. First of all, by (\ref{inductive cotas vq CN}) and (\ref{estimates psi}) we have
\[\normxt{0}{v_q\compt \psi^0-v_q}\lesssim \normxt{1}{v_q}\normxt{0}{\psi^0-\pi} \lesssim \delta_q^{1/2}\delta_{q+1}^{1/2}\lambda_q\eta^{-1}.\]
Meanwhile, arguing as in (\ref{diferencia aproximación B}) leads to
\[\normxtjm{0}{v_q\compt\psi^0-v^q_{jm}}\lesssim \delta_q^{1/2}\lambda_q\mu^{-1}.\]
Next, using (\ref{inductive cotas vq CN}) and (\ref{estimates psi}), it follows from (\ref{estimate composition}) that
\begin{align*}
\seminormxt{N}{v_q\compt\psi^0}&\lesssim \normxt{1}{v_q}\left(1+\normxt{N}{\psi^0-\pi}\right)+\normxt{N}{v_q}\left(1+\normxt{1}{\psi^0-\pi}\right)^N \\ &\lesssim \delta_q^{1/2}\delta_{q+1}^{1/2}\lambda_q\eta^{-1}\lambda_{q+1}^N+\delta_q^{1/2}\lambda_q^N(\delta_{q+1}^{1/2}\eta^{-1}\lambda_{q+1})^N\lesssim \delta_q^{1/2}\delta_{q+1}^{1/2}\lambda_q\eta^{-1}\lambda_{q+1}^N
\end{align*}
for $N=1,2,3,4$. Using (\ref{param 2}), we conclude that
\begin{align*}
\normxt{N}{v_q\compt\psi^0-v_q}&\lesssim \delta_q^{1/2}\delta_{q+1}^{1/2}\lambda_q\eta^{-1}\lambda_{q+1}^N\lesssim \delta_{q+2} \lambda_{q+1}^{N-12\tau}, \\[4pt] \normxt{N}{(v_q\compt\psi^0-v^q_{jm})u_{jm}}&\lesssim \normxtjm{0}{v_q\compt\psi^0-v^q_{jm}}\,\normxt{N}{u_{jm}}+\normxtjm{N}{v_q\compt\psi^0-v^q_{jm}}\normxt{0}{u_{jm}} \\&\lesssim \delta_q^{1/2}\delta_{q+1}^{1/2}\lambda_q\eta^{-1}\lambda_{q+1}^N \lesssim (\mu\lambda_{q+1}^{-1})\,\delta_{q+2}\lambda_{q+1}^{N-12\tau}
\end{align*}
for $N\leq 4$, $j\in\{1,\dots, 6\}$ and $m\in\Lambda$. Using (\ref{inductive cotas vq CN}), \eqref{param 9}, (\ref{estimates am}), (\ref{cota w0 con constantes}) and (\ref{cota wc}) to estimate the rest of the terms, we conclude
\begin{align}
	\normxt{N}{f_w}&\lesssim \delta_{q+2}\lambda_{q+1}^{N-12\tau} \qquad N=0, \dots, 4.
	\label{cota fw lema wphi}
\intertext{%
A similar argument leads to}
\normxt{N}{f_\Theta}&\lesssim \delta_{q+1}^{1/2}\lambda_{q+1}^{N} \hspace{32pt} N=0, \dots, 4.
	\label{cota fm lema wphi}
\intertext{%
In particular, we have}
\normxt{N}{F}&\lesssim \delta_{q+1}^{1/2}\lambda_{q+1}^{N} \hspace{32pt} N=0, \dots, 4. \nonumber
\intertext{%
On the other hand, it follows from (\ref{estimates am}) and (\ref{estimate product}) that}
\label{cotas matriz M sistema}
\normxt{N}{M}&\lesssim \frac{\delta_{q+1}^{1/2}\mu\lambda_{q+1}^N}{\eta} \hspace{24pt} N\geq 0.
\end{align}
Substituting these estimates into (\ref{sistema W}) yields
\[\normxt{0}{W}\lesssim \normxt{0}{F}+\normxt{0}{M}\normxt{0}{W}\lesssim \delta_{q+1}^{1/2}+\frac{\delta_{q+1}^{1/2}\mu}{\eta}\normxt{0}{W}.\]
Note that, although the number of pairs $(j,m)$ with $j\in\{1,\dots,6\}$ and $m\in\Lambda$ (and, hence, the size of these vectors) is proportional to $\mu^4$, each point has a neighborhood in which at most 16 of the amplitudes $a_{jm}$ is nonzero. Therefore, the number of entries of $M$ that are nonzero in a neighborhood of the point is independent of $\mu$. As a result, the implicit constants are, indeed, independent of the parameters of our construction. 

On the other hand, as argued several times, the coefficient multiplying the last norm can be made smaller that $1/2$ by choosing $a>1$ sufficiently large. Thus, we obtain
\[\normxt{0}{W}\lesssim \delta_{q+1}^{1/2}.\]
Regarding the $C^N$-norm for $N=1,2,3,4$, using (\ref{estimate product}) in (\ref{sistema W}) yields
\begin{align*}
\normxt{N}{W}&\lesssim\normxt{N}{F}+\normxt{N}{M}\normxt{0}{W}+\normxt{0}{M}\normxt{N}{W} \\ &\lesssim \delta_{q+1}^{1/2}\lambda_{q+1}^N+\delta_{q+1}^{1/2}\mu\eta^{-1}\normxt{N}{W}.
\end{align*}
Hence, for sufficiently large $a>1$, we have
\begin{equation}
\normxt{N}{W}\lesssim\delta_{q+1}^{1/2}\lambda_{q+1}^N \qquad\quad \forall N\leq 4. \label{cotas W lema wphi}
\end{equation}

We have estimated $W$, but the first component, $w^0_\phi$, satisfies better estimates. To derive them, we substitute these bounds for $W$ into (\ref{fórmula compacta tildew}), obtaining:
\begin{align}
\begin{split}
	\label{cotas tildew}
	\normxt{N}{w^0_\phi}&\lesssim\normxt{N}{f_w}+\normxt{N}{M}\normxt{0}{W}+\normxt{0}{M}\normxt{N}{W} \\ &\lesssim \delta_{q+2}\lambda_{q+1}^{N-11\tau}+\delta_{q+1}\mu\eta^{-1}\lambda_{q+1}^N\lesssim \delta_{q+2}\lambda_{q+1}^{N-11\tau} \hspace{30pt}\forall N\leq 4,
\end{split}
\end{align}
where we have used (\ref{param 3}). Combining this with (\ref{inductive tamaño vq}), (\ref{inductive cotas vq CN}), (\ref{cota w0 con constantes}) and (\ref{cota wc}), we also obtain
\begin{align}
\normxt{0}{v^0_{q+1}}&\lesssim 1, \label{cota tildev C0}\\
\normxt{N}{v^0_{q+1}}&\lesssim \delta_{q+1}^{1/2}\lambda_{q+1} \qquad N=1,2,3,4. \label{cota tildev CN}
\end{align}

Once we have studied $v^0_{q+1}$ and $w^0_\phi$, we can return to estimating $w_\phi$. Since the final diffeomorphism is $\phi^{q+1}=\phi^c\compt\phi^0$, by the chain rule, we have
\begin{align*}
	\partial_t\phi^{q+1}&=\partial_t\phi^c\compt\phi^0+(D\phi^c\compt\phi^0)\,\partial_t\phi^0, \\
	D\phi^{q+1}&=(D\phi^c\compt\phi^0)D\phi^0,
\intertext{%
so}
\partial_t\phi^{q+1}+D\phi^{q+1}v_q&=\partial_t\phi^c\compt\phi^0+(D\phi^c\compt\phi^0)(v^0_{q+1}\compt\phi^0).
\intertext{%
We conclude}
v_{q+1}&=\partial_t\phi^c\compt\psi^c+(D\phi^cv^0_{q+1})\compt \psi^c. \end{align*}
Therefore,
\begin{equation}
	w_\phi=w^0_\phi+\partial_t\phi^c\compt\psi^c+[(D\phi^c-\Id)v^0_{q+1}]\compt\psi^c+(v^0_{q+1}\compt\psi^c-v^0_{q+1}).
	\label{expresion wphi en función de tilde}
\end{equation}

We have already estimated the first term, so let us derive bounds for the $C^N$-norm of the rest. By (\ref{estimates phic}), (\ref{estimates psic}) and (\ref{estimate composition}), we have
\begin{align}
\normxt{0}{\partial_t\phi^c\compt\psi^c}&\lesssim \normxt{1}{\phi^c-\pi}\lesssim \delta_{q+1}\mu^3\eta^{-2}, \label{aux primer término wphi C0} \\ 
\normxt{N}{\partial_t\phi^c\compt\psi^c}&\lesssim \normxt{2}{\phi^c-\pi}\left(1+\normxt{N}{\psi^c-\pi}\right) \nonumber\\&\hspace{12pt}+\normxt{N+1}{\phi^c-\pi}\left(1+\normxt{1}{\psi^c-\pi}\right)^N \label{aux primer término wphi CN} \\&\lesssim \delta_{q+1}\mu^2\eta^{-2}\lambda_{q+1}^{N+\tau} \hspace{110pt} \forall N\geq 1. \nonumber
\end{align}
Concerning the next term, by (\ref{cota tildev C0}), we have
\begin{equation}
\label{aux segundo término wphi C0}
\normxt{0}{[(D\phi^c-\Id)v^0_{q+1}]\compt\psi^c}\lesssim \normxt{1}{\phi^c-\pi}\lesssim \delta_{q+1}\mu^3\eta^{-2}.
\end{equation}
Meanwhile, using (\ref{estimate composition}) with estimates (\ref{estimates phic}), (\ref{estimates psic}) and (\ref{cota tildev CN}) leads to the following bounds for $N=1,2,3,4$:
\begin{align}
\normxt{N}{[(D\phi^c-\Id)v^0_{q+1}]\compt\psi^c}&\lesssim \normxt{1}{\phi^c-\pi}\normxt{N}{v^0_{q+1}\compt\psi^c}+\normxt{N}{(D\phi^c-\Id)\compt\psi^c} \nonumber\\&\lesssim  \normxt{1}{\phi^c-\pi}\normxt{1}{v^0_{q+1}}\left(1+\normxt{N}{\psi^c-\pi}\right) \nonumber\\&\hspace{12pt}+ \normxt{1}{\phi^c-\pi}\normxt{N}{v^0_{q+1}}\left(1+\normxt{1}{\psi^c-\pi}\right)^N \nonumber \\
\begin{split}
\label{aux segundo término wphi CN}&\hspace{12pt}+ \normxt{2}{\phi^c-\pi}\left(1+\normxt{N}{\psi^c-\pi}\right) \\&\hspace{12pt}+\normxt{N+1}{\phi^c-\pi}\left(1+\normxt{1}{\psi^c-\pi}\right)^N 
\end{split} 
\\&\lesssim \delta_{q+1}^{5/2}\lambda_{q+1}\mu^5\eta^{-4}\max\{\mu^N,\lambda_{q+1}^{N-1+\tau}\} \nonumber \\&\hspace{12pt}+ \delta_{q+1}^{3/2}\mu^3\eta^{-2}\lambda_{q+1}^N+\delta_{q+1}\mu^2\eta^{-2}\lambda_{q+1}^{N+\tau} \nonumber  \\&\lesssim \delta_{q+1}^{3/2}\mu^3\eta^{-2}\lambda_{q+1}^{N}, \nonumber
\end{align}
where we have also used that $\phi^c$ and $\psi^c$ are close to $\pi$ in $C^2$ and we have simplified the expression using
\[\lambda_{q+1}^\tau\leq\delta_{q+1}^{1/2}\mu, \qquad \qquad \delta_{q+1}\mu^3\eta^{-2}\leq 1.\]

Regarding the last term in (\ref{expresion wphi en función de tilde}), we have
\begin{equation}
\label{aux tercer término wphi C0}
\normxt{0}{v^0_{q+1}\compt\psi^c-v^0_{q+1}}\lesssim \normxt{1}{v^0_{q+1}} \normxt{0}{\psi^c-\pi}\lesssim \frac{\delta_{q+1}^{3/2}\mu^2\lambda_{q+1}}{\eta^{2}}\lesssim \delta_{q+2}\lambda_{q+1}^{-11\tau},
\end{equation}
where we have used (\ref{param 5}). For $N=1,2,3,4$, it is better to isolate $w_0$ from the rest. First of all, by (\ref{estimate composition}), we have
\begin{align}
\normxt{N}{(v^0_{q+1}-w_0)\compt\psi^c-(v^0_{q+1}-w_0)} &\lesssim \normxt{N}{(v_q+w_c+w^0_\phi)\compt\psi^c}+\normxt{N}{v_q+w_c+w^0_\phi} \nonumber\\&\lesssim \normxt{1}{v_q+w_c+w^0_\phi}\left(1+\normxt{N}{\psi^c-\pi}\right) \nonumber \\&\hspace{12pt}+\normxt{N}{v_q+w_c+w^0_\phi}\left(1+\normxt{1}{\psi^c-\pi}\right)^N \label{aux tercer término wphi CN 1} \\&\lesssim \delta_{q+2}\lambda_{q+1}^{N-11\tau}+(\delta_{q+2}\lambda_{q+1}^{1-11\tau})\delta_{q+1}\mu^2\eta^{-2}\max\{\mu^N,\lambda_{q+1}^{N-1+\tau}\} \nonumber \\&\lesssim \delta_{q+2}\lambda_{q+1}^{N-11\tau}. \nonumber
\end{align}
Here we have substituted the estimates (\ref{inductive cotas vq CN}), (\ref{estimates psic}), (\ref{cota wc}) and (\ref{cotas tildew}) and we have used (\ref{def mu}) and (\ref{def eta}) to select the largest term. Meanwhile, by \cref{prop diferencia composición}, we have
\begin{align}
\normxt{N}{w_0\compt\psi^c-w_0}&\lesssim \normxt{1}{w_0}\normxt{N}{\psi^c-\pi}+\normxt{2}{w_0}\normxt{0}{\psi^c-\pi}\left(1+\normxt{N}{\psi^c-\pi}\right) \nonumber \\&\hspace{12pt}+\normxt{N+1}{w_0}\normxt{0}{\psi^c-\pi}\left(1+\normxt{1}{\psi^c-\pi}\right)^N\nonumber \\&\lesssim \delta_{q+1}^{3/2}\mu^2\eta^{-2}\lambda_{q+1}\max\{\mu^N,\lambda_{q+1}^{N-1+\tau}\}+\delta_{q+1}^{3/2}\mu^2\eta^{-2}\lambda_{q+1}^2 \label{aux tercer término wphi CN 2} \\&\hspace{12pt}+\delta_{q+1}^{5/2}\mu^4\eta^{-4}\lambda_{q+1}^2\max\{\mu^N,\lambda_{q+1}^{N-1+\tau}\}+\delta_{q+1}^{3/2}\mu^2\eta^{-2}\lambda_{q+1}^{N+1} \nonumber \\&\lesssim \delta_{q+2}\lambda_{q+1}^{N-11\tau}, \nonumber
\end{align}
where we have used (\ref{param 5}). 

Once we have estimated all of the terms in (\ref{expresion wphi en función de tilde}), to obtain the desired bounds (\ref{estimates wphi}) it suffices to insert (\ref{cotas tildew}) and (\ref{aux primer término wphi C0})-(\ref{aux tercer término wphi CN 2}) into (\ref{expresion wphi en función de tilde}) and use (\ref{param 5}).
\end{proof}

In particular, combining (\ref{cota w0 exacta}), (\ref{cota wc}) and (\ref{estimates wphi}) yields (\ref{inductive cotas vq CN}) and (\ref{cambio v q}) for sufficiently large $a>1$. Regarding (\ref{inductive tamaño vq}), we have
\[\delta_{q+1}^{1/2}\left(2^{-q}-2^{-(q+1)}\right)^{-1}= 2^{q+1}\delta_{q+1}^{1/2}=2^{q+1}a^{-\beta b^{q+1}},\]
which can be made arbitrarily small by taking $a>1$ large enough (depending on $\beta$ and $b$ but not on $q$), so (\ref{inductive tamaño vq}) holds, too. We conclude that the new velocity satisfies all of the required conditions.

Regarding the cross-helicity, we have
\[\int_{\RR^3}v_{q+1}\cdot B_{q+1}-\int_{\RR^3}v_q\cdot B_q=\int_{\RR^3}w_0\cdot B_q+\int_{\RR^3}(w_c\cdot B_q+w_\phi\cdot B_q+v_{q+1}\cdot b_{q+1}).\]
Using estimates (\ref{inductive tamaño vq}), (\ref{inductive cotas vq CN}) and (\ref{estimates am}), it follows from \cref{stationary phase lemma} with $m=1$ that
\[\abs{\int_{\RR^3}w_0\cdot B_q}\lesssim \delta_{q+1}^{1/2}(\mu+\delta_q^{1/2}\lambda_q)\lambda_{q+1}^{-1}\lesssim \delta_{q+2}\lambda_{q+1}^{-11\tau},\]
where we have also used (\ref{hierarchy}) and (\ref{def mu}). Recall that the support of all of the vector fields is contained in a common set $B(0,\bar{r})\times[0,T]$ independent of $q$, so its volume can be treated like a constant. Meanwhile, using (\ref{estimates b}), (\ref{cota wc}) and (\ref{estimates wphi}) to estimate the second integral in the change of the cross-helicity, we conclude
\[\abs{\int_{\RR^3}v_{q+1}\cdot B_{q+1}-\int_{\RR^3}v_q\cdot B_q}\lesssim \delta_{q+2}\lambda_{q+1}^{-11\tau},\]
from which (\ref{cambio helicidad q}) follows for sufficiently large $a>1$. Finally, let us study the change in the energy of the subsolution:
\begin{lemma}
\label{lema energía}
The new velocity $v_{q+1}$ and magnetic field $B_{q+1}$ satisfy
\[\abs{e(t)-\int_{\RR^3}\left(\abs{v_{q+1}}^2+\abs{B_{q+1}}^2\right)-\frac{1}{2}\delta_{q+2}\lambda_{q+1}^{-\tau}}\lesssim \delta_{q+2}\lambda_{q+1}^{-11\tau}.\]
\end{lemma}
\begin{proof}
By the definition of the perturbations, we have
\begin{align*}
\int_{\RR^3}\left(\abs{v_{q+1}}^2+\abs{B_{q+1}}^2\right)-\int_{\RR^3}\left(\abs{v_{q}}^2+\abs{B_{q}}^2\right)&=\int_{\RR^3}\left(\abs{w_0}^2+2w_0\cdot v_q\right)\\&\hspace{12pt}+\int_{\RR^3}\left[(v_q+v_{q+1})\cdot(w_c+w_\phi)+2b_{q+1}\cdot B_q+\abs{b_{q+1}}^2\right].
\end{align*}
It follows from (\ref{inductive tamaño vq}), (\ref{inductive tamaño Bq}), (\ref{estimates b}), (\ref{cota wc}) and (\ref{estimates wphi}) that 
\begin{align*}
\abs{\int_{\RR^3}\left[(v_q+v_{q+1})\cdot(w_c+w_\phi)+2b_{q+1}\cdot B_q+\abs{b_{q+1}}^2\right]}&\lesssim \delta_{q+1}^{1/2}\mu\lambda_{q+1}^{-1}+\delta_{q+2}\lambda_{q+1}^{-11\tau}\\&\lesssim \delta_{q+2}\lambda_{q+1}^{-11\tau},
\end{align*}
where we have also used (\ref{def mu}). Regarding the other term, since $\sum_m\chi_m^2=1$, it follows from \eqref{def am} that
\begin{equation}
\sum_{j,m}\frac{1}{2}a_{jm}^2\zeta_j\otimes \zeta_j=\sum_{j=1}^6\sigma_j(\mu \,\cdot\,)^2\hspace{1.5pt}\overline{\gamma}_{qj}^2\hspace{0.5pt}\zeta_j\otimes \zeta_j.
\label{media recupera matriz}
\end{equation}
Thus, using a well-known trigonometric identity yields
\begin{align*}
\int_{\RR^3}\left(\abs{w_0}^2+2w_0\cdot v_q\right)&=\int_{\RR^3}\sum_{j=1}^6\sigma_j(\mu \,\cdot\,)^2\hspace{1.5pt}\overline{\gamma}_{qj}^2\hspace{0.5pt}+\int_{\RR^3} \sum_{j,m} \frac{1}{2}a_{jm}^2\cos(2\theta_{jm})\\ &\hspace{12pt}+\int_{\RR^3}\sum_{j,\,m\neq n}\frac{1}{2}a_{jm}\hspace{1pt}a_{jn}[\cos(\theta_{jm}-\theta_{jn})+\cos(\theta_{jm}+\theta_{jn})] \\ &\hspace{12pt}+\int_{\RR^3}\sum_{j,m}2a_{jm}(\zeta_j\cdot v_q)\cos(\theta_{jm}).
\end{align*}
Remember that the product of two coefficients $a_{jm}$ with different values of $j$ vanishes, by \cref{remark soporte ajm disjuntos}. Using estimates (\ref{inductive cotas vq CN}) and (\ref{estimates am}) and \cref{stationary phase lemma} with $m=1$, we see that, up to a constant, the last three terms on the right-hand side are bounded by
\[\delta_{q+1}\mu\lambda_{q+1}^{-1}+\delta_{q+1}^{1/2}(\mu+\delta_q^{1/2}\lambda_q)\lambda_{q+1}^{-1}\lesssim \delta_{q+2}\lambda_{q+1}^{-11\tau},\]
where we have also used (\ref{def mu}) and (\ref{relación entre lambdas}). Using \eqref{diferencia gamatilde} to estimate the difference between $\overline{\gamma}_{qj}^2$ and $\gamma_{qj}^2$, we conclude that
\[\abs{\int_{\RR^3}\left(\abs{v_{q+1}}^2+\abs{B_{q+1}}^2\right)-\int_{\RR^3}\left(\abs{v_{q}}^2+\abs{B_{q}}^2\right)-\int_{\RR^3}\sum_{j=1}^6\sigma_j(\mu \,\cdot\,)^2\hspace{1.5pt}\gamma_{qj}^2}\lesssim \delta_{q+2}\lambda_{q+1}^{-11\tau}.\]
Here we have used that the support of $\overline{\gamma}_{qj}$ and $\gamma_{qj}$ is contained in $B(0,\overline{r})$, so its volume can be absorbed into the implicit constant. To study the third term, we decompose $\RR^3$ as a union of closed cubes $Q_i$ meeting only at their edges and having a side of length $2\pi\mu^{-1}$. Since we only have to cover $B(0,\overline{r})$, we may assume that the total volume of the cubes is a constant. We fix a point in each cube $x_i\in Q_i$ and we write:
\begin{equation}
\begin{split}
\int_{\RR^3}\sigma_j(\mu x)^2\hspace{1.5pt}\gamma_{qj}(x,t)^2\,dx&=\sum_{i}\int_{Q_i}\sigma_j(\mu x)^2\hspace{1.5pt}\gamma_{qj}(x_i,t)^2\,dx\\&\hspace{12pt}+\sum_{i}\int_{Q_i}\sigma_j(\mu x)^2\hspace{1.5pt}[\gamma_{qj}(x,t)^2-\gamma_{qj}(x_i,t)^2]\,dx.
\end{split}
\label{integral energía cubos}
\end{equation}
Since $\sigma_j$ is $2\pi$-periodic and since it does not matter the starting point when we integrate over a whole period, we have:
\begin{equation}
\int_{Q_i}\sigma_j(\mu x)^2\,dx=\int_{[0,2\pi]^3}\mu^{-3}\sigma_j(x)\,dx=(2\pi\mu^{-1})^3=\abs{Q_i},
\label{aux integral=volumen cubos}
\end{equation}
where we have used property (iii) in \cref{lemma mikado}. Thus, we can bound the second term on the right-hand side of \eqref{integral energía cubos} using \eqref{aux cotas gammaqj} and \eqref{relación entre lambdas}:
\begin{align*}
\sum_{i}\int_{Q_i}\sigma_j(\mu x)^2\hspace{1.5pt}\abs{\gamma_{qj}(x,t)^2-\gamma_{qj}(x_i,t)^2}\,dx&\lesssim \mu^{-1}\normxt{1}{\gamma_{qj}^2}\sum_{i}\int_{Q_i}\sigma_j(\mu x)^2\,dx=\mu^{-1}\normxt{1}{\gamma_{qj}^2}\sum_i \abs{Q_i}\\&\lesssim \mu^{-1}\normxt{1}{\gamma_{qj}^2}\lesssim \delta_{q+1}\lambda_q\mu^{-1}\lesssim \delta_{q+2}\lambda_{q+1}^{-11\tau}.
\end{align*}
Regarding the first term, we use \eqref{aux integral=volumen cubos} to write:
\[\sum_{i}\int_{Q_i}\sigma_j(\mu x)^2\hspace{1.5pt}\gamma_{qj}(x_i,t)^2\,dx=\sum_i \gamma_{qj}(x_i,t)^2 \abs{Q_i}=\int_{\RR^3} \gamma_{qj}^2\,dx+\int_{\RR^3} \left(\widetilde{\gamma}_{qj}^2-\gamma_{qj}^2\right)dx,\]
where
\[\widetilde{\gamma}_{qj}(x,t)\coloneqq \sum_i \gamma_{qj}(x_i,t) \hspace{1.5pt}\mathds{1}_{Q_i}\hspace{-1pt}(x) \qquad \Rightarrow \qquad \normxt{0}{\widetilde{\gamma}^2_{qj}-\gamma_{qj}^2}\lesssim \mu^{-1}\normxt{1}{\gamma_{qj}^2}\lesssim \delta_{q+2}\lambda_{q+1}^{-11\tau}.\]
In summary, so far we have proved:
\[\abs{\int_{\RR^3}\left(\abs{v_{q+1}}^2+\abs{B_{q+1}}^2\right)-\int_{\RR^3}\left(\abs{v_{q}}^2+\abs{B_{q}}^2\right)-\int_{\RR^3}\sum_{j=1}^6\gamma_{qj}^2}\lesssim \delta_{q+2}\lambda_{q+1}^{-11\tau}.\]
Due to \eqref{traza de Rq} and the fact that $\zeta_j$ is unitary, that taking the trace in \eqref{descomposición geométrica R_q} and integrating leads to
\[\int_{\RR^3}\sum_{j=1}^6\gamma_{qj}(x,t)^2\,dx=3\rho_q(t)\int_{\RR^3}\overline{\sigma}_q(x)^2\,dx\stackrel{\eqref{def rhoq}}{=}e(t)-\frac{1}{2}\delta_{q+2}\lambda_{q+1}^{-\tau}-\int_{\TT^3}\left(\abs{v_q}^2+\abs{B_q}^2\right)dx.\]
Substituting this in the previous estimate yields the result.
\end{proof}

We conclude that the new velocity $v_{q+1}$ and magnetic field $B_{q+1}$ satisfy all of the required properties: they are supported in $\Omega_{q,2}\times[0,T]$ and satisfy \eqref{inductive tamaño Bq}-\eqref{inductive cotas vq CN} and \eqref{inductive q energía}-\eqref{cambio helicidad q}. It remains to prove that the new Reynolds stress is sufficiently small.

\section{Estimates on the Reynolds stress}\label{section matrix}
The goal of this section is to check that the perturbations to the velocity and the magnetic field lead to a smaller Reynolds stress $R_{q+1}$. To derive the necessary estimates, it is convenient to consider the decomposition \eqref{def S1}-\eqref{def Rq+1}. Two of these terms can be easily estimated with the bounds that we already have. First of all, it follow from \eqref{aux cotas gammaqj}, \cref{lema aprox gammas} and the definition of $\mu\geq \lambda_q$ that the error term defined in (\ref{def S3}) satisfies:
\begin{equation}
	\normxt{0}{S_3}+\lambda_{q+1}^{-1}\normxt{1}{S_3}\lesssim\delta_{q+2}\lambda_{q+1}^{-11\tau}. \label{cotas S3}
\end{equation}
Similarly, it follows from (\ref{estimates b}), (\ref{cota wc}) and (\ref{estimates wphi}) that the error term defined in (\ref{def S6}) satisfies the following bounds:
\begin{equation}
\normxt{0}{S_6}+\lambda_{q+1}^{-1}\normxt{1}{S_6}\lesssim\delta_{q+2}\lambda_{q+1}^{-11\tau}. \label{cotas S6}
\end{equation}
The next subsection is devoted to studying $\mathcal{R}(\partial_t w_\phi)$, a term that is unique to our construction and that arises from defining the new velocity $v_{q+1}$ via diffeomorphisms rather than a closed expression. In the last subsection we will derive bounds for the remaining error terms, which are more standard.

\subsection{The new error term}
\label{SS.Rpdtwphi}
As we shall see later, $\mathcal{R}(\partial_t w_c)$ is easy to estimate due to its simple explicit formula. In contrast, $\mathcal{R}(\partial_t w_\phi)$ poses a greater challenge since it lacks an explicit formula and is instead given by a complex implicit expression.  

Beyond these technical difficulties, a crucial factor in obtaining acceptable bounds is that $\phi^c$ and $\psi^c$ are very close to $\pi$ in $C^2$. This is ensured by selecting a diffeomorphism $\psi^0$ that is nearly volume-preserving. If we had not required to be so close to~$\pi$, we could have omitted the third term in (\ref{def psi}), which would significantly simplify \cref{section estimates diffeo} and \cref{section estimates vector fields}. In conclusion, estimating $\mathcal{R}(\partial_t w_\phi)$ is the most delicate part of our construction, introducing additional complexity that would otherwise be unnecessary.  

Before we begin with more complex estimates, we establish the following elementary fact. We refer to Appendix~\ref{appendix a} for the definition of the ``triple'' norm we use in the following. Recall that $\tau>0$ is a very small parameter given by \cref{lema relaciones entre parámetros}.
\begin{lemma}
	\label{lema fácil Rdt}
A map $g\in C^2(\RR^3\times[0,T], \RR^3)$ satisfies:
\[\normatres{\mathcal{R}(\partial_tg)}_{\tau}\lesssim \normxt{1}{g}, \qquad \normatres{\mathcal{R}(\partial_tg)}_{1+\tau}\lesssim \normxt{2}{g}.\]
\end{lemma}
\begin{proof}
Since $\mathcal{R}$ is a bounded operator from $B^{s}_{\infty,\infty}(\RR^3)$ to $B^{s+1}_{\infty,\infty}(\RR^3)$, for a fixed time $t\in[0,T]$ we have
\[\normx{\tau}{\mathcal{R}[\partial_tg(\cdot,t)]}\lesssim \left\|\partial_tg(\cdot,t)\right\|_{B^{-1+\tau}_{\infty,\infty}}\lesssim \normx{0}{\partial_tg(\cdot,t)}\lesssim \normxt{1}{g}.\]
Similarly, 
\begin{align*}
	\normx{\tau}{\partial_t\mathcal{R}[\partial_tg(\cdot,t)]}&=\normx{\tau}{\mathcal{R}[\partial_t^2g(\cdot,t)]}\lesssim \left\|\partial_t^2g(\cdot,t)\right\|_{B^{-1+\tau}_{\infty,\infty}}\lesssim\normxt{2}{g}, \\
	 \normx{1+\tau}{\mathcal{R}[\partial_tg(\cdot,t)]}&\lesssim \normx{\tau}{\partial_tg(\cdot,t)}\lesssim\normxt{2}{g},
\end{align*}
where we have used that $\mathcal{R}$ commutes with $\partial_t$. The required estimates follow.
\end{proof}
 Although the previous estimates may appear crude, one cannot expect a much better behavior from a general function, since the smoothing effect of $\mathcal{R}$ applies only to spatial oscillations. Consequently, only functions whose temporal oscillations are not significantly faster than their spatial ones can benefit from improved estimates due to the action of $\mathcal{R}$. Fortunately, many functions in our construction exhibit this behavior, which we make more precise in the following definition:

\begin{definition}
	\label{def F}
	Given $N^\ast,r^\ast\in\NN$, we define $\mathcal{F}(N^\ast,r^\ast)$ to be the class of smooth of functions that can be written as
	\[f=g_0+\sum_{j,l} g_{jl}\,\sin(\tilde{\theta}_{jl})\]
	and such that
	\begin{enumerate}
		\item there exists $J^\ast (f)\in\NN$ such that every point has a neighborhood in which at most $J^\ast (f)$ of the amplitudes $g_{jl}$ do not vanish, 
		\item the functions $g_0,g_{jl}$ are smooth and there exists constants $C_0(f)$, $C_1(f)$ such that:
		\begin{align}
		 \normxt{N}{g_0}&\lesssim C_0(f)\hspace{0.5pt}\mu^N, \qquad N=0,1,\dots, N^\ast,\label{cotas coef g0}\\ \normxt{N}{g_{jl}}&\lesssim C_1(f)\hspace{0.5pt}\mu^N, \hspace{21pt} N=0,1,\dots, N^\ast, \label{cotas coef gi}
		\end{align}
		where the implicit constants depend on $N$ but not on $\mu$ or other parameters of the iterative process, 
		\item The functions $\tilde{\theta}_{jl}$ are non-vanishing sums of the form
		\[\theta_{jm_1}\pm \cdots \pm \theta_{jm_r}, \quad\qquad m_s\in\Lambda,\; r\leq r^\ast,\]
		plus possibly a constant angle $\pm\pi/2$, with the phase functions defined in~\eqref{def thetam}.
		\item If $\theta_{jm}$ is present in the sum $\tilde{\theta}_{jl}$, then the support of $g_{jl}$ is contained in the support of $a_{jm}$.
	\end{enumerate}
	To exclude artificial terms, we also assume that the amplitudes $g_{jl}$ are not identically 0 and that the phases $\tilde{\theta}_{jl}$ are different from one another. 
\end{definition}
We establish some basic properties of these functions in the following lemma:
\begin{lemma}
	\label{lema propiedades mathcal F}
	Let $N^\ast,r^\ast_1,r^\ast_2 \in \NN$. Given $f\in \mathcal{F}(N^\ast,r^\ast_1)$, the phases $\tilde{\theta}_{jl}$ satisfy
	\begin{equation}
		|\nabla\tilde{\theta}_{jl}|\sim \lambda_{q+1}, \qquad |\partial_t\tilde{\theta}_{jl}|\lesssim \lambda_{q+1}, \label{derivada thetai}
	\end{equation}
	where the implicit constants depend only on $r^\ast_1$. In addition, given $f_1\in \mathcal{F}(N^\ast,r^\ast_1)$ and $f_2\in \mathcal{F}(N^\ast,r^\ast_2)$, the product $f_1f_2$ is in $\mathcal{F}(N^\ast,r^\ast_1+r^\ast_2)$ and we have
	\begin{align}
	J^\ast(f_1f_2)&\leq J^\ast(f_1)+J^\ast(f_2)+  2J^\ast(f_1)J^\ast(f_2),
	\label{number of terms}\\
	C_0(f_1f_2)&\leq C_0(f_1)\,C_0(f_2)+ J_\ast(f_1)\,J_\ast(f_2)\,C_1(f_1)\,C_1(f_2), \label{constante C0 producto} \\
	\begin{split}C_1(f_1f_2)&\leq J_\ast(f_2)\,C_0(f_1)\,C_1(f_2)+J_\ast(f_1)\,C_1(f_1)\,C_0(f_2)\\&\hspace{12pt}+2J_\ast(f_1)\,J_\ast(f_2)\,C_1(f_1)\,C_1(f_2). \label{constante C1 producto}\end{split}
	\end{align}
\end{lemma}
\begin{proof}
	It follows from the definition of $\theta_{jm}$ that \[\abs{\nabla\theta_{jm}}=\ell_m\lambda_{q+1}, \qquad \abs{\partial_t\theta_{jm}}\leq \ell_m\lambda_{q+1}.\] Therefore, to prove (\ref{derivada thetai}), it suffices to show that the corresponding sum of the coefficients $\ell_m$ is nonzero, and to find a lower bound on its absolute value depending only on $r^\ast$.
	
	It follows from the definition of $\mathcal{L}$ (cf. Equation~\eqref{eq.L}) that the only linear combinations of elements $l_j\in\mathcal{L}$ with coefficients $\pm 1$ that can be equal to 0 are the trivial ones, i.e.,
	\[l_{1}-l_{1}+l_{2}-l_{2}+\cdots+l_{n}-l_{n}.\]
	From this we immediately deduce that if $r$ is odd, any sum $\sum_{j=1}^r\pm \ell_{m_j}$ with $\ell_{m_j}\in \mathcal{L}$ must be nonzero. It is clear that the minimum of all such sums only depends on $r$, hence, on $r^\ast$. 
	
	On the other hand, suppose that $r$ is even and that we have a sum $\sum_{j=1}^r\pm \ell_{m_j}=0$ with $\ell_{m_i}\in \mathcal{L}$. Then, after possibly reordering the terms, we must have 
	\[\ell_{m_1}=\ell_{m_2}, \quad \ell_{m_3}=\ell_{m_4}, \quad \dots, \quad \ell_{m_{r-1}}=\ell_{m_r}\]
	and the sum must be alternated. Suppose that one of the phases $\tilde{\theta}_{jl}$ in $f$ is associated to this sum, that is, \[\tilde{\theta}_{jl}=\theta_{jm_1}-\theta_{jm_2}+\theta_{jm_3}-\theta_{jm_4}+\cdots+\theta_{jm_{r-1}}-\theta_{jm_r}.\]
	By property $(\text{iv})$, the support of the non-vanishing amplitude $g_{jl}$ must be contained in the intersection of the supports of $\chi_{m_1}$ and $\chi_{m_2}$. Since $\ell_{m_1}=\ell_{m_2}$, by construction of the coefficients $\ell_m$, this intersection is empty unless $m_1=m_2$. Repeating this argument for the other indices, we conclude that they form pairs of equal indices but associated to opposite signs. This means that $\tilde{\theta}_{jl}=0$. Since we are excluding such terms, this means that no phase $\tilde{\theta}_{jl}$ is associated to a sum $\sum_{j=1}^r\pm \ell_{m_j}$ that vanishes. Again, the minimum value that such sums can take depends only on $r$ (hence, on $r^\ast$) and our choice of $\mathcal{L}$. 
	
	Concerning the second claim in the lemma, it follows from property (iv) in the definition and from \cref{remark soporte ajm disjuntos} that we only have to consider the product of terms with the same value of $j$. By a well-known product-to-sum trigonometric identity,
	\begin{align}
	\left[g_0^1+\sum_{j,l_1} g^1_{jl_1}\sin(\tilde{\theta}^1_{jl_1})\right]\left[g_0^2+\sum_{j,l_2} g^2_{jl_2}\sin(\tilde{\theta}^2_{jl_2})\right]&=g_0^1g_0^2+\sum_{j,l_1} g_0^2g^1_{jl_1}\sin(\tilde{\theta}^1_{jl_1})+\sum_{j,l_2} g_0^1g^2_{jl_2}\sin(\tilde{\theta}^2_{jl_2}) \nonumber \\&\hspace{12pt}+\sum_{j,l_1,l_2}\left(\frac{1}{2}g_{jl_1}^1g_{jl_2}^2\right)\sin\left(\tilde{\theta}^1_{jl_1}-\tilde{\theta}^2_{jl_2}+\frac{\pi}{2}\right) \label{producto F} \\&\hspace{12pt}+\sum_{j,l_1,l_2}\left(-\frac{1}{2}g_{jl_1}^1g_{jl_2}^2\right)\sin\left(\tilde{\theta}^1_{jl_1}+\tilde{\theta}^2_{jl_2}+\frac{\pi}{2}\right). \nonumber
	\end{align}
	The new phases are a sum of at most $r_1^\ast+r_2^\ast$ phases $\theta_{jm}$ with $m\in\Lambda$ and the corresponding amplitudes satisfy the support property $(\text{iv})$ because so did $g^1_{jl_1}$ and $g_{jl_2}^2$.
	
	Regarding the number of terms, it follows from (\ref{producto F}) that each high-frequency term in one of the functions yields a high-frequency term when multiplied by the low-frequency term of the other function. Thus, each point has a neighborhood in which at most $J^\ast(f_1)+J^\ast(f_2)$ of such terms do not vanish. Meanwhile, we observe in (\ref{producto F}) that the product of a high-frequency term coming from each function leads to (at most) two high frequency terms. Hence, each point has a neighborhood in which at most $2J^\ast(f_1)J^\ast(f_2)$ of such terms do not vanish. Note that some of the sums $\tilde{\theta}^1_{jl_1}\pm \tilde{\theta}^2_{jl_2}$ could be constant. Those terms should be included in the low-frequency term.
	
	Finally, let us derive bounds for the constants $C_0(f_1f_2)$ and $C_1(f_1f_2)$. To estimate the former, it suffices to count the number of terms in (\ref{producto F}) that could contribute to the low frequency term and apply (\ref{estimate product}). Regarding the crude bound (\ref{constante C1 producto}), we have just assumed the worst case scenario, i.e., that all of the new phases are equal. Then, it suffices to count the number of high-frequency terms and apply (\ref{estimate product}), as in (\ref{constante C0 producto}).
\end{proof}
Therefore, these function classes are stable under products (in an appropriate sense). Moreover, if the number of terms that are nonzero in a neighborhood of a point is independent of the parameters $\delta_q$ and $\lambda_q$, then the same holds for the product. We will now study the action of $\mathcal{R}$ on these spaces:
\begin{lemma}
\label{R en mathcal F}
Let $N^\ast, r^\ast\geq 1$ be integers and let $f\in C^\infty(\RR^3\times[0,T],\RR^3)$ be a smooth vector field such that $f\in \mathcal{F}(N^\ast,r^\ast)$. With the notation of \cref{def F}, we have
\begin{align}
\normatres{\mathcal{R}f}_{\tau}&\lesssim C_0(f)\hspace{1pt}+C_1(f)\left(\lambda_{q+1}^{-1+\tau}+\mu^{N^\ast}\lambda_{q+1}^{-N^\ast}\right), \label{cota tres barras tau mathcal F} \\
\normatres{\mathcal{R}f}_{1+\tau}&\lesssim C_0(f)\hspace{1pt}\mu+C_1(f)\hspace{1pt}\lambda_{q+1}\left(\lambda_{q+1}^{-1+\tau}+\mu^{N^\ast}\lambda_{q+1}^{-N^\ast}\right), \label{cota tres barras tau mathcal F 2}
\end{align}
where the implicit constants depend on $\tau$ and on $N^\ast$, $r^\ast$ and $J^\ast(f)$.
\end{lemma}
\begin{proof}
Let us first focus on (\ref{cota tres barras tau mathcal F}). Taking into account that $\mathcal{R}$ is a bounded operator from $B^{-1+\tau}_{\infty,\infty}(\RR^3)$ to $C^\tau(\RR^3)$, we have
\[\normx{\tau}{\mathcal{R}f(\cdot,t)}\lesssim \norm{f(\cdot,t)}_{B^{-1+\tau}_{\infty,\infty}}\]
for any fixed $t\in[0,T]$. We estimate the Besov norm of the low-frequency part of $f$, that is, $g_0$, by its $C^0$-norm. To estimate the Besov norm of the remaining terms we use \cref{stationary phase lemma}. We conclude
\[\normx{\tau}{\mathcal{R}f(\cdot,t)}\lesssim C_0(f)\hspace{1pt}+C_1(f)\left(\lambda_{q+1}^{-1+\tau}+\mu^{N^\ast}\lambda_{q+1}^{-N^\ast}\right).\]
Taking the supremum in $t\in[0,T]$ yields (\ref{cota tres barras tau mathcal F}). Note that so far we have only needed $N^\ast\geq 0$; the condition $N^\ast\geq 1$ is only required in (\ref{cota tres barras tau mathcal F 2}). Regarding the spatial derivatives, we have:
\[ \max_{t\in[0,T]}\normx{1+\tau}{\mathcal{R}f(\cdot,t)}\lesssim \max_{t\in[0,T]}\normx{\tau}{f(\cdot,t)}\lesssim C_0\hspace{1pt}(f)\mu+C_1(f)\hspace{1pt}\lambda_{q+1}^\tau,\]
where we have used again that $\mathcal{R}$ is an operator of order $-1$. Next, with the notation of \cref{def F}, we write the temporal derivative as
\[\partial_t f=\left[\partial_t g_0+\sum_{j,l} \partial_t g_{jl}\,\sin(\tilde{\theta}_{jl})\right]+\sum_{j,l} g_{jl}\,\partial_t\tilde{\theta}_{jl}\cos(\tilde{\theta}_{jl})\equiv \tilde{f}_1+\tilde{f}_2.\]
By \cref{def F}, we see that $\tilde{f}_1\in\mathcal{F}(N^\ast-1,r^\ast)$ with $C_0(\tilde{f}_1)=C_0(f)\mu$ and $C_1(\tilde{f}_1)=C_1(f)\mu$. Hence, we can apply (\ref{cota tres barras tau mathcal F}), obtaining:
\[\max_{t\in[0,T]}\normx{\tau}{\mathcal{R}\tilde{f}_1(\cdot,t)}=\normatres{\mathcal{R}\tilde{f}_1}_{\tau}\lesssim C_0(f)\hspace{1pt}\mu+C_1(f)\hspace{1pt}\mu\left[\lambda_{q+1}^{-1+\tau}+\mu^{N^\ast-1}\lambda_{q+1}^{-(N^\ast-1)}\right].\]
As mentioned earlier, (\ref{cota tres barras tau mathcal F}) holds for any $N^\ast \geq 0$, so the previous estimate is justified. Meanwhile, by (\ref{inductive tamaño vq}), (\ref{def thetam}) and property (iii) in \cref{def F}, we see that $\partial_t\tilde{\theta}_{jl}$ is constant and that $|\partial_t\tilde{\theta}_{jl}|\lesssim \lambda_{q+1}$. Therefore, $\tilde{f}_2\in\mathcal{F}(N^\ast,r^\ast)$ with $C_0(\tilde{f}_1)=0$ and $C_1(\tilde{f}_2)=C_1(f)\lambda_{q+1}$. Again, applying (\ref{cota tres barras tau mathcal F}) yields
\[\max_{t\in[0,T]}\normx{\tau}{\mathcal{R}\tilde{f}_2(\cdot,t)}=\normatres{\mathcal{R}\tilde{f}_2}_{\tau}\lesssim C_1(f)\hspace{1pt}\lambda_{q+1}\left(\lambda_{q+1}^{-1+\tau}+\mu^{N^\ast}\lambda_{q+1}^{-N^\ast}\right).\]
Taking into account that $\partial_t$ and $\mathcal{R}$ commute, we conclude
\begin{align*}
\normatres{\mathcal{R}f}_{1+\tau}&=\max_{t\in[0,T]}\left(\normx{\tau}{\mathcal{R}f(\cdot,t)}+\normx{1+\tau}{\mathcal{R}f(\cdot,t)}+\normx{\tau}{\partial_t\mathcal{R}f(\cdot,t)}\right) \\&\lesssim \max_{t\in[0,T]}\left(\normx{\tau}{\mathcal{R}f(\cdot,t)}+\normx{1+\tau}{\mathcal{R}f(\cdot,t)}+\normx{\tau}{\mathcal{R}\tilde{f}_1(\cdot,t)}+\normx{\tau}{\mathcal{R}\tilde{f}_2(\cdot,t)}\right) \\ &\lesssim C_0(f)\hspace{1pt}\mu+C_1(f)\hspace{1pt}\lambda_{q+1}\left(\lambda_{q+1}^{-1+\tau}+\mu^{N^\ast}\lambda_{q+1}^{-N^\ast}\right),
\end{align*}
as desired.
\end{proof}
We are now ready to determine the effect of $\mathcal{R}\partial_t$ on $\mathcal{F}(N^\ast,r^\ast)$. We do not expect the bound for the low frequency term $g_0$ to be better than \cref{lema fácil Rdt}. Nevertheless, this estimate is acceptable in this case because $g_0$ oscillates with frequency $\mu\ll \lambda_{q+1}$. Meanwhile, we expect $\mathcal{R}\partial_t$ to preserve the size of the high-frequency terms because the frequency of the temporal oscillations is not larger than the frequency of the spatial oscillations. In practice we get a factor $\lambda_{q+1}^\tau$ because we have to use the $C^\tau_x$ norm instead of the $C^0$ norm.
\begin{corollary}
\label{Rdt en mathcal F}
Let $N^\ast,r^\ast\in \NN$ with $N^\ast\geq 3$ and let $f\in \mathcal{F}(N^\ast,r^\ast)$. With the notation of \cref{def F}, we have
\begin{align}
\normatres{\mathcal{R}(\partial_tf)}_\tau&\lesssim C_0(f)\hspace{0.5pt}\mu+C_1(f)\hspace{0.5pt}\lambda_{q+1}^\tau, \label{estimates Rpartialt f C0} \\ \normatres{\mathcal{R}(\partial_tf)}_{1+\tau}&\lesssim C_0(f)\hspace{0.5pt}\mu^2+C_1(f)\hspace{0.5pt}\lambda_{q+1}^{1+\tau}, \label{estimates Rpartialt f C1}
\end{align}
where the implicit constants depend on $r^\ast$ and $J^\ast(f)$ but not on $\mu$ or $\lambda_{q+1}$.
\end{corollary}
\begin{proof}
With the notation of \cref{def F}, we write the temporal derivative as
\[\partial_t f=\left[\partial_t g_0+\sum_{j,l} \partial_t g_{jl}\,\sin(\tilde{\theta}_{jl})\right]+\sum_{j,l} g_{jl}\,\partial_t\tilde{\theta}_{jl}\cos(\tilde{\theta}_{jl})\equiv \tilde{f}_1+\tilde{f}_2.\]
We see that $\tilde{f}_1\in\mathcal{F}(2,r^\ast)$ with $C_0(\tilde{f}_1)=C_0(f)\hspace{1pt}\mu$ and $C_1(\tilde{f}_1)=C_1(f)\hspace{1pt}\mu$. Hence, by \cref{R en mathcal F}, we have
\begin{align*}
\normatres{\mathcal{R}\tilde{f}_1}_\tau&\lesssim C_0(f)\hspace{1pt}\mu+ C_1(f)\hspace{1pt}\mu\hspace{0.5pt}\left(\lambda_{q+1}^{-1+\tau}+\mu^2\lambda_{q+1}^{-2}\right), \\
\normatres{\mathcal{R}\tilde{f}_1}_{1+\tau}&\lesssim C_0(f)\hspace{1pt}\mu^2+ C_1(f)\hspace{1pt}\mu\hspace{1pt}\lambda_{q+1}\left(\lambda_{q+1}^{-1+\tau}+\mu^2\lambda_{q+1}^{-2}\right).
\end{align*}
Meanwhile, as argued in \cref{R en mathcal F}, we have $\tilde{f}_2\in\mathcal{F}(3,r^\ast)$ with $C_0(\tilde{f}_2)=0$ and $C_1(\tilde{f}_2)=C_1(f)\hspace{1pt}\lambda_{q+1}$. Applying \cref{R en mathcal F} to $\tilde{f}_2$ thus yields
\begin{align*}
	\normatres{\mathcal{R}\tilde{f}_2}_\tau&\lesssim  C_1(f)\hspace{1pt}\lambda_{q+1}\hspace{0.5pt}\left(\lambda_{q+1}^{-1+\tau}+\mu^3\lambda_{q+1}^{-3}\right), \\
	\normatres{\mathcal{R}\tilde{f}_2}_{1+\tau}&\lesssim  C_1(f)\hspace{1pt}\hspace{1pt}\lambda_{q+1}^2\left(\lambda_{q+1}^{-1+\tau}+\mu^3\lambda_{q+1}^{-3}\right).
\end{align*}
Combining the bounds for both terms and using
\begin{equation}
	\mu^3\leq \lambda_{q+1}^{2+\tau}
	\label{aux relación parámetros g0gi}
\end{equation}
leads to the desired estimates for $\mathcal{R}(\partial_tf)$. This identity follows from (\ref{def mu}).
\end{proof}
Our strategy to obtain suitable bounds for $\mathcal{R}(\partial_t w_\phi)$ will then be to write $w_\phi$ as a sum of a function in $\mathcal{F}(3,r^\ast)$ for some $r^\ast\in \NN$ independent of $a>1$ and a remainder that is very small in $C^1$. The former can be estimated using the previous lemma, whereas the latter can be estimated using \cref{lema fácil Rdt}.

Nevertheless, there is an obstacle when trying to put this strategy into practice: the last term in (\ref{expresion wphi en función de tilde}) is not in $\mathcal{F}(3,r^\ast)$ and neither is it small in $C^1$. Indeed, since $\psi^c$ is $C^2$-close to the identity, one would expect to obtain an appropriate bound by simply applying \cref{prop diferencia composición}. However, that is not the case, because $\lambda_{q+1}^2\normxt{0}{\psi^c-\Id}$ is not small. In order to derive the required estimates, it is necessary to exploit the special structure of the functions in $\mathcal{F}(3,r^\ast)$:
\begin{lemma}
\label{lema Rdt diferencia F}
Let $N^\ast,r^\ast\in \NN$ with $N^\ast\geq 3$ and let $f\in \mathcal{F}(N^\ast,r^\ast)$. With the notation of \cref{def F}, we have
\begin{align}
	\normatres{\mathcal{R}\partial_t(f\compt\psi^c-f)}_\tau&\lesssim \frac{\delta_{q+1}\mu^2}{\eta^{2}}\left[C_0(f)\hspace{0.5pt}\mu^2+C_1(f)\hspace{0.5pt}\lambda_{q+1}^{1+\tau}\right] \label{estimates Rpartialt f diferencia composición C0} \\ \normatres{\mathcal{R}\partial_t(f\compt\psi^c-f)}_{1+\tau}&\lesssim \frac{\delta_{q+1}\mu^2}{\eta^{2}}\left[C_0(f)\hspace{0.5pt}\mu^2+C_1(f)\hspace{0.5pt}\lambda_{q+1}^{1+\tau}\right]\lambda_{q+1}, \label{estimates Rpartialt f diferencia composición C1}
\end{align}
where the implicit constants depend on $\tau$, $r^\ast$ and $J^\ast(f)$ but not on $\mu$ or $\lambda_{q+1}$.
\end{lemma}
\begin{proof}
We will use the notation of \cref{def F}. First of all, applying \cref{prop diferencia composición} with estimates (\ref{estimates psic}) and (\ref{cotas coef gi}) yields
\begin{align*}
\normxt{0}{g_{jl}\compt\psi^c-g_{jl}}&\lesssim \normxt{1}{g_{jl}}\normxt{0}{\psi^c-\pi}\lesssim C_1(f)\hspace{0.5pt}\delta_{q+1}\mu^3\eta^{-2}, \\ 
\normxt{N}{g_{jl}\compt\psi^c-g_{jl}} &\lesssim \normxt{1}{g_{jl}}\normxt{N}{\psi^c-\pi}+ \normxt{0}{\psi^c-\pi}\normxt{2}{g_{jl}}\left(1+\normxt{N}{\psi^c-\pi}\right)\\&\hspace{12pt}+\normxt{0}{\psi^c-\pi}\normxt{N+1}{g_{jl}}\left(1+\normxt{1}{\psi^c-\pi}\right)^N \\ &\lesssim C_1(f)\hspace{0.5pt}\delta_{q+1}\mu^3\eta^{-2}\max\{\mu^N,\lambda_{q+1}^{N-1+\tau}\},
\end{align*}
for $N=1,2$. We have simplified the expressions using that $\psi^c$ is close to $\pi$ in $C^2$. In conclusion,
\begin{equation}
\normxt{N}{g_{jl}\compt\psi^c-g_{jl}}\lesssim C_1(f)\hspace{0.5pt}\delta_{q+1}\mu^3\eta^{-2}\max\{\mu^N,\lambda_{q+1}^{N-1+\tau}\} \qquad N=0,1,2.
\label{cotas diferencia gj con psic}
\end{equation}
It is clear that $g_0$ satisfies analogous bounds, with $C_1(f)$ replaced by $C_0(f)$. Hence, by \cref{lema fácil Rdt}, we have
\begin{align}
\normatres{\mathcal{R}\partial_t(g_0\compt\psi^c-g_0)}_\tau &\lesssim C_0(f)\hspace{0.5pt}\delta_{q+1}\mu^4\eta^{-2}, \label{cota Rdt diferencia g0 C0}\\
\normatres{\mathcal{R}\partial_t(g_0\compt\psi^c-g_0)}_{1+\tau} &\lesssim C_0(f)\hspace{0.5pt}\delta_{q+1}\mu^3\eta^{-2}\lambda_{q+1}^{1+\tau}\lesssim C_0(f)\hspace{0.5pt}\delta_{q+1}\mu^4\eta^{-2}\lambda_{q+1}. \label{cota Rdt diferencia g0 C1}
\end{align}
	
Next, we study how the operator $\mathcal{R}\partial_t$ acts on the remaining terms. Since the supports of at most $J_\ast(f)$ of the functions $g_{jl}$ have nonempty intersection, it follows from \cref{stationary phase lemma} that we can focus on a single term of the sum in \cref{def F}. Note that the implicit constants in the estimates will then depend on $J_\ast(f)$. Furthermore, since $\mathcal{R}\partial_t$ commutes with taking the imaginary part, it suffices to study
\[z_{jl}\coloneqq g_{jl}\hspace{0.5pt}e^{i\tilde{\theta}_{jl}}.\]
We see that we can write
\begin{align*}
z_{jl}\compt\psi^c-z_{jl}&=(g_{jl}\compt\psi^c)\hspace{0.5pt}e^{\tilde{\theta}_{jl}\compt\psi^c}-g_{jl}\hspace{0.5pt}e^{i\tilde{\theta}_{jl}}\\&=(g_{jl}\compt\psi^c-g_{jl})\hspace{0.5pt}e^{\tilde{\theta}_{jl}\compt\psi^c}+g_{jl}\left(e^{\tilde{\theta}_{jl}\compt\psi^c}-e^{i\tilde{\theta}_{jl}}\right) \\&=\left[(g_{jl}\compt\psi^c-g_{jl})\hspace{0.5pt}e^{\tilde{\theta}_{jl}\compt\psi^c-\tilde{\theta}_{jl}}+g_{jl}\left(e^{\tilde{\theta}_{jl}\compt\psi^c-\tilde{\theta}_{jl}}-1\right)\right]\hspace{0.5pt}e^{i\tilde{\theta}_{jl}} \equiv \tilde{g}_{jl}\hspace{0.5pt}e^{i\tilde{\theta}_{jl}}.
\end{align*}
Note that
\[\normxt{N}{\exp(\tilde{\theta}_{jl}\compt\psi^c-\tilde{\theta}_{jl})-1}\lesssim \normxt{N}{\tilde{\theta}_{jl}\compt\psi^c-\tilde{\theta}_{jl}}\lesssim \lambda_{q+1}\normxt{N}{\psi^c-\pi}. \]
Combining this with (\ref{cotas diferencia gj con psic}) and using (\ref{estimate product}) yields
\begin{equation}
\normxt{N}{\tilde{g}_{jl}}\lesssim C_1(f)\hspace{0.5pt}\delta_{q+1}\mu^2\eta^{-2}\lambda_{q+1}\max\{\mu^N,\lambda_{q+1}^{N-1+\tau}\} \qquad N=0,1,2.
\end{equation}
We have
\[\partial_t\left(\tilde{g}_{jl}\hspace{0.5pt}e^{i\tilde{\theta}_{jl}}\right)=(\partial_t\tilde{g}_{jl}+i\partial_t\tilde{\theta}_{jl}\tilde{g}_{jl})\hspace{0.5pt}e^{i\tilde{\theta}_{jl}},\]
where $\partial_t\tilde{\theta}_{jl}$ is just a constant of order $\lambda_{q+1}$. Therefore, applying \cref{stationary phase lemma} to the first term with $m=1$ and to the second term with $m=2$ yields
\begin{align*}
\norm{\partial_t\left(\tilde{g}_{jl}\hspace{0.5pt}e^{i\tilde{\theta}_{jl}}\right)(\cdot,t)}_{B^{-1+\tau}_{\infty,\infty}}&\lesssim \lambda_{q+1}^{-1+\tau}\normxt{1}{\tilde{g}_{jl}}+\lambda_{q+1}^{-1}\normxt{2}{\tilde{g}_{jl}}+\lambda_{q+1}^{\tau}\normxt{0}{\tilde{g}_{jl}}\\&\lesssim C_1(f)\hspace{0.5pt}\delta_{q+1}\mu^2\eta^{-2}\lambda_{q+1}^{1+\tau}.
\end{align*}
for a fixed time $t\in[0,T]$. Thus,
\begin{align}
\normatres{\mathcal{R}\partial_t\left(\tilde{g}_{jl}\hspace{0.5pt} e^{i\tilde{\theta}_{jl}}\right)}_{\tau}&=\max_{t\in[0,T]}\normx{\tau}{\mathcal{R}\partial_t\left(\tilde{g}_{jl}\hspace{0.5pt}e^{i\tilde{\theta}_{jl}}\right)(\cdot,t)}\nonumber\\&\lesssim \max_{t\in[0,T]}\norm{\partial_t\left(\tilde{g}_{jl}\hspace{0.5pt}e^{i\tilde{\theta}_{jl}}\right)(\cdot,t)}_{B^{-1+\tau}_{\infty,\infty}}\lesssim C_1(f)\hspace{0.5pt}\delta_{q+1}\mu^2\eta^{-2}\lambda_{q+1}^{1+\tau}.
\label{Rdt tildegj}
\end{align}
On the other hand,
\[\partial_t^2\left(\tilde{g}_{jl}\hspace{0.5pt}e^{i\tilde{\theta}_{jl}}\right)=[\partial_t^2\tilde{g}_{jl}+2i\partial_t\tilde{g}_{jl}\partial_t\tilde{\theta}_{jl}-\tilde{g}_{jl}(\partial_t\tilde{\theta}_{jl})^2]\hspace{0.5pt}e^{i\tilde{\theta}_{jl}}.\]
Hence, for a fixed $t\in[0,T]$ we have
\begin{align*}
\norm{\partial_t^2\left(\tilde{g}_{jl}\hspace{0.5pt}e^{i\tilde{\theta}_{jl}}\right)(\cdot,t)}_{B^{-1+\tau}_{\infty,\infty}}&\lesssim \normxt{2}{\tilde{g}_{jl}}+\lambda_{q+1}^\tau\normxt{1}{\tilde{g}_{jl}}+\lambda_{q+1}^{1+\tau}\normxt{0}{\tilde{g}_{jl}}\\&\lesssim C_1(f)\hspace{0.5pt}\delta_{q+1}\mu^2\eta^{-2}\lambda_{q+1}^{2+\tau},
\end{align*}
where we have estimated the Besov norm of the first term by its $C^0$-norm and we have applied \cref{stationary phase lemma} with $m=1$ to the second term and with $m=2$ to the third term. Then,
\begin{align*}
\max_{t\in[0,T]}\normx{\tau}{\partial_t \mathcal{R}\partial_t\left(\tilde{g}_{jl}\hspace{0.5pt}e^{i\tilde{\theta}_{jl}}\right)(\cdot,t)}&= \max_{t\in [0,T]}\normx{\tau}{ \mathcal{R}\left[\partial_t^2\left(\tilde{g}_{jl}\hspace{0.5pt}e^{i\tilde{\theta}_{jl}}\right)\right](\cdot,t)}\\&\lesssim \max_{t\in[0,T]}\norm{\partial_t^2\left(\tilde{g}_{jl}\hspace{0.5pt}e^{i\tilde{\theta}_{jl}}\right)(\cdot,t)}_{B^{-1+\tau}_{\infty,\infty}}\lesssim C_1(f)\hspace{0.5pt}\delta_{q+1}\mu^2\eta^{-2}\lambda_{q+1}^{2+\tau}.
\end{align*}
Meanwhile,
\begin{align*}
	\max_{t\in [0,T]}\normx{1+\tau}{ \mathcal{R}\left[\partial_t\left(\tilde{g}_{jl}\hspace{0.5pt}e^{i\tilde{\theta}_{jl}}\right)\right](\cdot,t)}&\leq \max_{t\in[0,T]}\normx{\tau}{\partial_t\left(\tilde{g}_{jl}\hspace{0.5pt}e^{i\tilde{\theta}_{jl}}\right)(\cdot,t)}\lesssim \normxt{1+\tau}{\tilde{g}_{jl}\hspace{0.5pt}e^{i\tilde{\theta}_{jl}}}\\&\lesssim \normxt{1+\tau}{\tilde{g}_{jl}}+\normxt{0}{\tilde{g}_{jl}}\normxt{1+\tau}{\exp(i\tilde{\theta}_{jl})}\\&\lesssim \normxt{2}{\tilde{g}_{jl}}+\normxt{0}{\tilde{g}_{jl}}\lambda_{q+1}^{1+\tau}\lesssim C_1(f)\hspace{0.5pt}\delta_{q+1}\mu^2\eta^{-2}\lambda_{q+1}^{2+\tau},
\end{align*}
where we have used that $\mathcal{R}$ is an operator of order $-1$. We conclude
\[\normatres{\mathcal{R}\partial_t\left(\tilde{g}_{jl}\hspace{0.5pt}e^{i\tilde{\theta}_{jl}}\right)}_{1+\tau}\lesssim C_1(f)\hspace{0.5pt}\delta_{q+1}\mu^2\eta^{-2}\lambda_{q+1}^{2+\tau}.\]
Combining this with (\ref{cota Rdt diferencia g0 C0}), (\ref{cota Rdt diferencia g0 C1}) and (\ref{Rdt tildegj}) yields the desired bounds.
\end{proof}

We are now ready to estimate $\mathcal{R}(\partial_t w_\phi)$:
\begin{lemma}
The error matrix $S_2=\mathcal{R}(\partial_t w_\phi)$ satisfies
\begin{equation}
\normatres{S_2}_\tau+\lambda_{q+1}^{-1}\normatres{S_2}_{1+\tau}\lesssim \delta_{q+2}\lambda_{q+1}^{-11\tau}.
\label{cotas S2}
\end{equation}
\end{lemma}
\begin{proof}
To write $w_\phi$ in a suitable form, we will use a system of equations similar to (\ref{sistema W}). The difference is that here we must pay more attention to which terms are in $\mathcal{F}(3,r^\ast)$ for some $r^\ast\in \NN$. For instance, the difference $v_q\compt\psi^0-v_q$ does not have the desired form. Since it is not small in $C^1$, we cannot simply dismiss this term, it requires further study. By Taylor's theorem with remainder in integral form, we have 
\[v_q\compt\psi^0-v_q=Dv_q(\psi^0-\pi)+Z,\]
where
\begin{equation}
	Z\coloneqq \sum_{i,j=1}^3\int_0^1\frac{1-s}{2-\delta_{ij}}\left[\partial_{ij} v_q\compt\left((1-s)\pi+s\psi^0\right)\right](\psi^0-\pi)_i(\psi^0-\pi)_j\;ds.
	\label{def R error taylor}
\end{equation}
It follows from (\ref{inductive cotas vq CN}), (\ref{hierarchy}), (\ref{def psi}) and (\ref{estimates am}) that $Dv_q(\psi^0-\pi)$ is in $\mathcal{F}(3,1)$, but we cannot say anything about $Z$. Therefore, we will try to write a system similar to (\ref{sistema W}) but isolating the terms that contain $Z$. Let us define the coefficients
\begin{align*}
	h_w&\coloneqq Z+\sum_{j,m}\ell_m\mu^{-1}\lambda_{q+1}(k_{jm}\cdot Z)\,u_{jm},\\[6pt]
	\tilde{f}_w&\coloneqq Dv_q(\psi^0-\Id)+\sum_{j,m} \ell_m\mu^{-1}\lambda_{q+1}k_{jm}\cdot[(v_q-v^q_{jm})+Dv_q(\psi^0-\pi)]\,u_{jm} \\ &\hspace{12pt}+\sum_{j,m}\frac{D_v a_{jm}+(w_0+w_c)\cdot\nabla a_{jm}}{\ell_m\eta}\,\zeta_j\sin(\theta_{jm}) \\ &\hspace{12pt}+\sum_{j,m}\frac{(\zeta_j\otimes k_{jm}-k_{jm}\otimes \zeta_j)[D_v\nabla a_{jm}+D^2 a_{jm}(w_0+w_c)]}{\ell_m^2\eta\lambda_{q+1}}\cos(\theta_{jm}), \\[6pt] 
	h_{jm}&\coloneqq \sum_{n\in\Lambda}\lambda_{q+1}(k_{jn}\cdot Z)\frac{\ell_mk_{jm}\cdot[\nabla a_{jn}\times(\zeta_j\times k_{jn})]}{\ell_n\mu\eta}\sin(\theta_{jn}),
	\\[6pt]
	\tilde{f}_{jm}&\coloneqq \sum_{n\in\Lambda}\frac{\ell_mk_{jm}\cdot[\nabla a_{jn}\times(\zeta_j\times k_{jn})]}{\ell_n\mu}\sin(\theta_{jn})\\&\hspace{12pt}+\sum_{n\in\Lambda}\lambda_{q+1}k_{jn}\cdot [(v_q-v_{q,n})+Dv_q(\psi^0-\Id)]\,\frac{\ell_mk_{jm}\cdot[\nabla a_{jn}\times(\zeta_j\times k_{jn})]}{\ell_n\mu\eta}\sin(\theta_{jn})\\&\hspace{12pt}-\sum_{n\in\Lambda}\frac{\ell_mk_{jm}\cdot(\zeta_j\otimes k_{jm}-k_{jm}\otimes \zeta_j)[D_v\nabla a_{jn}+D^2a_{jn}(w_0+w_c)]}{\ell_n^2\mu\eta}\cos(\theta_{jn}),
\end{align*}
Furthermore, let us construct a vector $H$ containing $h_w$ as the first entry and $h_{jm}$, with $j\in\{1,\dots,6\}$ and $m\in\Lambda$, as the rest of the coordinates. We do the same with $\tilde{f}_w$ and $\tilde{f}_{jm}$, obtaining a vector $\widetilde{F}$. With the notation of \cref{lema wphi}, we see that 
\begin{equation}
F=H+\widetilde{F}.
\label{F suma H y tildeF}
\end{equation}
Therefore, we can write (\ref{sistema W}) as
\[W=H+\widetilde{F}+MW.\]
It follows from (\ref{inductive cotas vq CN}), (\ref{hierarchy}), (\ref{def psi}) and (\ref{estimates am}) that $\widetilde{F}, M\in \mathcal{F}(3,2)$. Taking into account that we only control up to the $C^4$-norm of $v_q$ and that $Dv_q$ is present in the expression of $\widetilde{F}$, we see that we can only control up to the $C^3$-norm of the amplitudes.

Substituting the expression for $W$ into the right-hand side of the system, we see that we can write
\[W=[(\Id+M)H+M^2W]+(\Id+M)\widetilde{F}.\]
By (\ref{number of terms}), we see that the last term is in $\mathcal{F}(3,4)$. Let us denote by $P_1(\cdot)$ the projection onto the first coordinate, so that $P_1(W)=w^0_\phi$, $P_1(H)=h_w$ and so on. Then, we may write
\[w^0_\phi=\hat{h}+\hat{f},\]
where
\begin{align}
	\hat{h}&\coloneqq P_1[(\Id+M)H+M^2W], \\ \hat{f}&\coloneqq \tilde{f}_w+P_1(M\widetilde{F}). \label{def f hat}
\end{align}
In conclusion, we have written $w^0_\phi$ as a term $\hat{f}\in\mathcal{F}(3,4)$ and a remainder $\hat{h}$, which will turn out to be small in $C^1$. 

Let us estimate each term. First of all, by (\ref{inductive cotas vq CN}) and (\ref{estimates psi}), we have
\[\normxt{0}{Z}\lesssim \normxt{2}{v_q}\normxt{0}{\psi^0-\pi}^2\lesssim \frac{\delta_q^{1/2}\delta_{q+1}\lambda_q^2}{\eta^2}.\]
Next, we fix $N\in\{1,2\}$. It follows from (\ref{estimate product}) that
\begin{align*}
\normxt{N}{Z}&\lesssim \max_{s\in[0,1]}\normxt{N}{D^2v_q\compt((1-s)\pi+s\psi^0)}\normxt{0}{\psi^0-\pi}^2\\&\hspace{12pt}+\normxt{2}{v_q}\normxt{0}{\psi^0-\pi}\normxt{N}{\psi^0-\pi}.
\end{align*}
Regarding the first term, by (\ref{estimate composition}), we have
\begin{align*}
\normxt{N}{D^2v_q\compt((1-s)\pi+s\psi^0)}&\lesssim \normxt{3}{v_q}\left(1+\normxt{N}{\psi^0-\pi}\right)\\&\hspace{12pt}+\normxt{N+2}{v_q}\left(1+\normxt{1}{\psi^0-\pi}\right)^N \\&\lesssim \delta_q^{1/2}\lambda_q^3\delta_{q+1}^{1/2}\eta^{-1}\lambda_{q+1}^N+\delta_q^{1/2}\lambda_{q}^{N+2}\left(\delta_{q+1}^{1/2}\eta^{-1}\lambda_{q+1}\right)^N \\ &\lesssim \delta_{q}^{1/2}\delta_{q+1}^{1/2}\lambda_q^3\eta^{-1}\lambda_{q+1}^N,
\end{align*}
where we have used that $\delta_{q+1}^{1/2}\lambda_q\eta^{-1}\leq 1$. Substituting this into our previous expression and using (\ref{estimates psi}) leads to
\[\normxt{N}{Z}\lesssim \delta_{q}^{1/2}\delta_{q+1}^{3/2}\lambda_q^3\eta^{-3}\lambda_{q+1}^N+\delta_q^{1/2}\delta_{q+1}\lambda_q^2\eta^{-2}\lambda_{q+1}^N \qquad N\in\{1,2\}.\]
Using again $\delta_{q+1}^{1/2}\lambda_q\eta^{-1}\leq 1$, we conclude
\[\normxt{N}{Z}\lesssim \delta_q^{1/2}\delta_{q+1}\lambda_q^2\eta^{-2}\lambda_{q+1}^N \hspace{30pt} N=0,1,2.\]
Hence, applying (\ref{estimate product}) to the definition of $H$ and using (\ref{estimates am}) leads to
\begin{align}
\begin{split}
\label{cotas H error chungo}
\normxt{N}{H}&\lesssim \frac{\delta_{q+1}^{1/2}\lambda_{q+1}}{\eta}\normxt{N}{Z}+\frac{\delta_{q+1}^{1/2}\lambda_{q+1}^{N+1}}{\eta}\normxt{0}{Z}\\
&\lesssim \frac{\delta_q^{1/2}\delta_{q+1}^{3/2}\lambda_{q}^2\lambda_{q+1}^{N+1}}{\eta^3}\lesssim \delta_{q+2}\lambda_{q+1}^{N-1-11\tau}, \hspace{40pt} N=0,1,2,
\end{split}
\end{align}
where we have used (\ref{param 6}). Meanwhile, $M$ was estimated in (\ref{cotas matriz M sistema}):
\[\normxt{N}{M}\lesssim \frac{\delta_{q+1}^{1/2}\mu\lambda_{q+1}^N}{\eta}\qquad N\geq 0.\]
By (\ref{estimate product}), we conclude that
\[\normxt{N}{P_1[(\Id+M)H]}\lesssim \delta_{q+2}\lambda_{q+1}^{N-1-11\tau} \qquad N=0,1,2. \]
Regarding the other term in $\hat{h}$, we have
\begin{align*}
\normxt{N}{P_1(M^2W)}&\lesssim \normxt{0}{M}^2\normxt{N}{W}+\normxt{0}{M}\normxt{N}{M}\normxt{0}{W}\\&\lesssim\left(\delta_{q+1}^{1/2}\mu\eta^{-1}\right)^2\delta_{q+1}^{1/2}\lambda_{q+1}^N\lesssim \delta_{q+2}\lambda_{q+1}^{N-1-11\tau}, \hspace{30pt} N=0,1,2.
\end{align*}
where we have used (\ref{param 5}) and (\ref{cotas W lema wphi}). We conclude
\begin{equation}
\normxt{N}{\hat{h}}\lesssim \delta_{q+2}\lambda_{q+1}^{N-1-11\tau} \qquad N=0,1,2.
\label{cotas h hat}	
\end{equation}

Let us study $\hat{f}$. First of all, using (\ref{cota fw lema wphi}), (\ref{cota fm lema wphi}) and (\ref{cotas H error chungo}), it follows from (\ref{F suma H y tildeF}) that
\begin{align*}
	\normxt{0}{\tilde{f}_w}&\leq \normxt{0}{f_w}+\normxt{0}{h_w}\lesssim \delta_{q+2}\lambda_{q+1}^{-12\tau}, \\
	\normxt{0}{\tilde{f}_m}&\leq \normxt{0}{f_m}+\normxt{0}{h_m}\lesssim \delta_{q+1}^{1/2} \hspace{30pt} m\in \Lambda,
\end{align*}
which yields 
\[\normxt{0}{\widetilde{F}}\lesssim \delta_{q+1}^{1/2}.\]
substituting these estimates into (\ref{def f hat}) leads to
\[\normxt{0}{\hat{f}}\lesssim \normxt{0}{\tilde{f}_w}+\normxt{0}{M}\normxt{0}{\tilde{F}}\lesssim \delta_{q+2}\lambda_{q+1}^{-12\tau}+\delta_{q+1}\mu\eta^{-1}\lesssim \delta_{q+2}\lambda_{q+1}^{-12\tau},\]
by (\ref{param 9}). Note that in these estimates we have always taken the supremum on the size of the various terms and we have never relied on cancellations between different terms. Therefore, this bound on the size of $\hat{f}$ applies to each separate term. Hence, we have
\begin{equation}
C_0(\hat{f}\hspace{1pt})+C_1(\hat{f}\hspace{1pt})\lesssim \delta_{q+2}\lambda_{q+1}^{-12\tau}. \label{cota C1 f hat}
\end{equation}
We can find a better bound for $C_0(\hat{f}\hspace{1pt})$ by studying $C_0(\tilde{f}_w)$. Substituting (\ref{def psi}), (\ref{def w0}), (\ref{def wc}) and (\ref{def um}) into the definition of $\tilde{f}_w$, we see that the only terms that can yield low frequencies are
\begin{itemize}
	\setlength\itemsep{5pt}
	\item $\sum_{j,m} \ell_m\mu^{-1}\lambda_{q+1}k_{jm}\cdot Dv_q(\psi^0-\pi)u_{jm}$, 
	\item $\sum_{j,m}\ell_m^{-1}\eta^{-1}\nabla a_{jm}\cdot w_c\,\zeta_j \sin(\theta_{jm})$, 
	\item $\sum_{j,m}\ell_m^{-2}\eta^{-1}\lambda_{q+1}^{-1}(\zeta_j\otimes k_{jm}-k_{jm}\otimes\zeta_j)D^2a_{jm}w_0\cos(\theta_{jm})$.
\end{itemize}
Taking into account that $\theta_{jm}\pm \theta_{jn}\neq 0$ for any two distinct indices $m,n\in\Lambda$ such that the support of $\chi_m$ and $\chi_n$ have nonempty intersection, we see that the low frequency term in $\tilde{f}_w$ is
\begin{align*}
g_0(\tilde{f}_w)&=\sum_{j,m}\frac{a_{jm}}{2\ell_m^2\eta^2}k_{jm}\cdot Dv_q\left[(\nabla a_{jm}\times(\zeta_j\times k_{jm}))\otimes\zeta_j-\zeta_j\otimes(\nabla a_{jm}\times(\zeta_j\times k_{jm}))\right] \\ &\hspace{12pt}+\sum_{j,m}\frac{a_{jm}}{2\ell_m^2\eta\lambda_{q+1}}(\zeta_j\otimes k_{jm}-k_{jm}\otimes \zeta_j)D^2a_{jm}\hspace{0.5pt}\zeta_j,
\end{align*}
where we have just substituted (\ref{def psi}), (\ref{def w0}), (\ref{def wc}) and (\ref{def um}) and used the cancellation
\[\nabla a_{jm}\cdot[\nabla a_{jm}\times(\zeta_j\times k_{jm})]=0.\]
Using (\ref{inductive cotas vq CN}) and (\ref{estimates am}), we conclude that 
\[C_0(\tilde{f}_w)\lesssim \frac{\delta_q^{1/2}\delta_{q+1}\lambda_q\mu}{\eta^2}+\frac{\delta_{q+1}\mu^2}{\eta\lambda_{q+1}}\lesssim \delta_{q+2}\mu^{-1}\lambda_{q+1}^{-11\tau},\]
by (\ref{param 7}) and (\ref{param 3}). Using these bounds and the relationship (\ref{param 3}), it follows from (\ref{def f hat}) that
\begin{align}
\begin{split}
	C_0(\hat{f}\hspace{1pt})\lesssim C_0(\tilde{f}_w)+\normxt{0}{M}\normxt{0}{\widetilde{F}}&\lesssim \delta_{q+2}\mu^{-1}\lambda_{q+1}^{-11\tau}+\delta_{q+1}\mu\eta^{-1}\\&\lesssim \delta_{q+2}\mu^{-1}\lambda_{q+1}^{-11\tau}. \label{cota C0 f hat}
\end{split}
\end{align}

Let us summarize what we have achieved so far: we have written $w^0_\phi$ as a the sum of a term $\hat{h}$, which is small in $C^1$, and a term $\hat{f}\in \mathcal{F}(3,4)$. In addition, we have obtained bounds (\ref{cotas h hat}) for $\hat{h}$ and bounds (\ref{cota C1 f hat}), (\ref{cota C0 f hat}) for $\hat{f}$. Once we have derived these estimates, along with the previous lemmas, we are finally ready to study $\mathcal{R}\partial_t(w_\phi)$. First, we use (\ref{expresion wphi en función de tilde}) to write
\begin{align}
\begin{split}
\label{expresión wphi error}
w_\phi&=\left[\hat{h}+\partial_t\phi^c\compt\psi^c+[(D\phi^c-\Id)v^0_{q+1}]\compt\psi^c\right]+\hat{f}+(v^0_{q+1}\compt\psi^c-v^0_{q+1}) \\ &\equiv w_I+\hat{f}+w_{III}.
\end{split}
\end{align}
It follows from \cref{lema fácil Rdt} with estimates (\ref{aux primer término wphi CN}), (\ref{aux segundo término wphi CN}) and (\ref{cotas h hat}) that
\begin{align*}
\normatres{\mathcal{R}\partial_t(w_I)}_\tau+\lambda_{q+1}^{-1}\normatres{\mathcal{R}\partial_t(w_I)}_{1+\tau}&\lesssim \delta_{q+2}\lambda_{q+1}^{-11\tau}+\delta_{q+1}\mu^2\eta^{-2}\lambda_{q+1}^{1+\tau}+\delta_{q+1}^{3/2}\mu^3\eta^{-2}\lambda_{q+1}\\&\lesssim \delta_{q+2}\lambda_{q+1}^{-11\tau}
\end{align*}
by (\ref{param 1}). Next, it follows from \cref{Rdt en mathcal F} with estimates (\ref{cota C1 f hat}) and (\ref{cota C0 f hat}) that
\[\normatres{\mathcal{R}\partial_t(\hat{f})}_\tau+\lambda_{q+1}^{-1}\normatres{\mathcal{R}\partial_t(\hat{f})}_{1+\tau}\lesssim \delta_{q+2}\lambda_{q+1}^{-11\tau}.\] 
Finally, it follows from (\ref{inductive tamaño vq}) and (\ref{cota w0 con constantes})-(\ref{estimates wphi}) that
\[C_0(v^0_{q+1})\lesssim 1, \qquad C_1(v^0_{q+1})\lesssim \delta_{q+1}^{1/2}.\]
Then, by \cref{lema Rdt diferencia F}, we have
\[\normatres{\mathcal{R}\partial_t(w_{III})}_\tau+\lambda_{q+1}^{-1}\normatres{\mathcal{R}\partial_t(w_{III})}_{1+\tau}\lesssim \frac{\delta_{q+1}\mu^4}{\eta^2}+\frac{\delta_{q+1}^{3/2}\mu^2\lambda_{q+1}^{1+\tau}}{\eta^2}\lesssim \delta_{q+2}\lambda_{q+1}^{-11\tau},\]
by (\ref{param 3}) and (\ref{param 5}). Combining the estimates for the three terms in (\ref{expresión wphi error}) yields the desired bound.
\end{proof}

\subsection{Remaining error terms}
We will now study the error matrices $S_1=\mathcal{R}(z_1)$,  $S_4=\mathcal{R}(z_2)$ and $S_5=\mathcal{R}(z_3)$, which were defined in \eqref{def S1}, \eqref{def S4} and \eqref{def S5}, respectively. Taking into account that $\nabla\theta_{jm}$ is orthogonal to $\zeta_j$ for any $m\in\Lambda$, it follows from (\ref{def w0}) and (\ref{def wc}) that
\begin{align*}
z_1&=\sum_{j,m}-a_{jm}(\partial_t\theta_{jm}+v_q\cdot\nabla\theta_{jm})\zeta_j\sin(\theta_{jm}) \\&\hspace{12pt}+\sum_{j,m}\left[(\partial_t a_{jm}+v_q\cdot\nabla a_{jm})\zeta_j+(\zeta_j\cdot\nabla a_{jm})v_q+a_{jm}(\zeta_j\cdot\nabla v_q)\right]\cos(\theta_{jm})\\&\hspace{12pt}+\sum_{j,m}\left[\frac{(\nabla\partial_ta_{jm})\times(\zeta_j\times k_{jm})}{\ell_m\lambda_{q+1}}\sin(\theta_{jm})+\frac{\nabla a_{jm}\times(\zeta_j\times k_{jm})}{\ell_m\lambda_{q+1}}(\partial_t\theta_{jm})\cos(\theta_{jm})\right].
\end{align*}
If follows from (\ref{def thetam}) that we may write
\[\partial_t\theta_{jm}+v_q\cdot\nabla\theta_{jm}=(v_q-v^q_{jm})\cdot\nabla \theta_{jm}+\ell_m\eta,\]
where $v^q_{jm}$ was defined in (\ref{def vJm}). Since $a_{jm}$ is supported in a cube centered at $\mu^{-1}m$ and whose sidelength is comparable to $\mu^{-1}$, we see that
\[\sup_{(x,t)\in\supp a_{jm}}\abs{v_q(x,t)-v^q_{jm}}\lesssim \normxt{1}{v_q}\mu^{-1}\stackrel{(\ref{inductive cotas vq CN})}{\lesssim} \delta_q^{1/2}\lambda_q\mu^{-1}.\]
Taking into account that $\lambda_q\ll\mu$ and
\[\delta_q^{1/2}\lambda_q\mu^{-1}\lambda_{q+1}\leq \eta\]
by (\ref{param 2}), it follows from (\ref{inductive cotas vq CN}) and (\ref{estimates am}) that
\[\normxt{N}{\partial_t\theta_{jm}+v_q\cdot\nabla\theta_{jm}}\lesssim \eta\mu^N\]
for $N=0,1,2,3$. Using (\ref{inductive tamaño vq}), (\ref{inductive cotas vq CN}) and (\ref{estimates am}) to estimate the rest of the amplitudes in $z_1$ and taking into account that $\mu^2\lambda_{q+1}^{-1}\leq \eta$, we conclude that $z_1\in \mathcal{F}(3,1)$ with $C_0(z_1)=0$ and $C_1(z_1)=\delta_{q+1}^{1/2}\eta$. \cref{R en mathcal F} thus leads to
\begin{align}
\begin{split}
\normatres{S_1}_\tau+\lambda_{q+1}^{-1}\normatres{S_1}_{1+\tau}\lesssim \delta_{q+1}^{1/2}\eta(\lambda_{q+1}^{-1+\tau}+\mu^3\lambda_{q+1}^{-3})&\stackrel{(\ref{aux relación parámetros g0gi})}{\lesssim}\delta_{q+1}^{1/2}\eta\lambda_{q+1}^{-1+\tau}\\&\hspace{2pt}\stackrel{(\ref{def eta})}{=} \delta_{q+2}\lambda_{q+1}^{-11\tau}.
\end{split}
\label{cotas S1}
\end{align}

Regarding $S_4$, using \eqref{media recupera matriz} and a well-known trigonometric identity, we can write:
\[z_2= \Div\hspace{-1.5pt}\left(\sum_{j,\,m\neq n}\frac{1}{2}a_{jm}a_{jn}\zeta_j\otimes\zeta_j[\cos(\theta_{jm}-\theta_{jn})+\cos(\theta_{jm}+\theta_{jn})] +\sum_{j,m} \frac{1}{2}a_{jm}^2\zeta_j\otimes\zeta_j\cos(2\theta_{jm})\right).\]
Remember that the product of two coefficients $a_{jm}$ with different values of $j$ vanishes by property (ii) in \cref{lemma mikado}. Since $\nabla\theta_{jm}$ is orthogonal to $\zeta_j$ for any $m\in\Lambda$, we obtain
\[z_2= \sum_{j,\, m\neq n}a_{jn}(\zeta_j\cdot\nabla a_{jm})\,\zeta_j\,[\cos(\theta_{jm}-\theta_{jn})+\cos(\theta_{jm}+\theta_{jn})] +\sum_{j,m} a_{jm}(\zeta_j\cdot\nabla a_{jm})\,\zeta_j\cos(2\theta_{jm}).\]
By (\ref{estimates am}), we see that $z_2\in \mathcal{F}(3,2)$ with $C_0(z_2)=0$ and $C_1(z_2)=\delta_{q+1}\mu$. \cref{R en mathcal F} thus leads to
\begin{align}
\begin{split}
\normatres{S_4}_\tau+\lambda_{q+1}^{-1}\normatres{S_4}_{1+\tau}\lesssim \delta_{q+1}\mu(\lambda_{q+1}^{-1+\tau}+\mu^3\lambda_{q+1}^{-3})&\stackrel{(\ref{aux relación parámetros g0gi})}{\lesssim}\delta_{q+1}\mu\lambda_{q+1}^{-1+\tau}\\&\hspace{2.5pt}\stackrel{(\ref{def mu})}{\lesssim} \delta_{q+2}\lambda_{q+1}^{-11\tau}.
\end{split}
\label{cotas S4}
\end{align}

Concerning $S_5$, it follows from property (iii) in \cref{lemma mikado} that $1-\sigma_j^2$ is mean-free. Therefore, we can write
\[1-\sigma_j(x)^2=\sum_{l\in\ZZ^3\backslash 0}c_{jl}e^{il\cdot x}\]
for some coefficients $c_{jl}\in \mathbb{C}$ which, on account of the smoothness of $\sigma_j$, can be assumed to satisfy
\[\abs{c_{jl}}\leq \frac{C}{\abs{l}^4},\]
for some constant $C>0$ independent of $j$ and $l$. Hence, $z_3$ can be written as:
\[z_3=\sum_{l\in \ZZ^3\backslash 0}\sum_{j=1}^6c_{jl}\Div(\overline{\gamma}_{qj}^2\zeta_j\otimes \zeta_j)e^{i\mu l\cdot x}.\]
Applying \cref{stationary phase lemma} and summing over $l\in \ZZ^3\backslash 0$, we obtain
\[\norm{z_3(\cdot,t)}_{B^{-1+\tau}_{\infty,\infty}}\lesssim \frac{\delta_{q+1}}{\mu^{1-\tau}}\left(\frac{\delta_{q+1}\lambda_q\lambda_{q+1}^{11\tau}}{\delta_{q+2}}\right)+\frac{\delta_{q+1}}{\mu^{N}}\left(\frac{\delta_{q+1}\lambda_q\lambda_{q+1}^{11\tau}}{\delta_{q+2}}\right)^{\hspace{-2pt}N+1},\]
where we have used \eqref{cotas gamatilde} and we are free to choose $N\geq 0$. Using \eqref{relación entre lambdas}, we write
\begin{align*}
\norm{z_3(\cdot,t)}_{B^{-1+\tau}_{\infty,\infty}}&\lesssim \delta_{q+2}\lambda_{q+1}^{-11\tau}\left[\frac{\delta_{q+1}^2\lambda_q\lambda_{q+1}^{20\tau}}{\delta_{q+2}^2\mu^{1-\tau}}+\frac{\delta_{q+1}^2\lambda_q\lambda_{q+1}^{20\tau}}{\delta_{q+2}^2\mu^{N\tau}}\left(\frac{\delta_{q+1}^2\lambda_q\lambda_{q+1}^{20\tau}}{\delta_{q+2}^2\mu^{1-\tau}}\right)^{\hspace{-3pt}N\;}\right] \\
&\lesssim \delta_{q+2}\lambda_{q+1}^{-11\tau}\left(1+\delta_{q+1}^2\delta_{q+2}^{-2}\lambda_q\lambda_{q+1}^{20\tau}\mu^{-N\tau}\right).
\end{align*}
It follows from \eqref{relación entre lambdas} that the second term in parenthesis can be made of order 1 by taking $N$ sufficiently large (depending on $\tau$). The implicit constant will then depend on $\tau$ through $N$, but having fixed the rest of the parameters, we will take $a>1$ sufficiently large to get rid of the constants. Therefore, since $\mathcal{R}$ is an operator of order -1, we get the following bound on $S_5=\mathcal{R}(z_3)$:
\[\normxt{0}{S_5}\leq\max_{t\in [0,T]}\norm{S_5(\cdot,t)}_{C^\tau}\lesssim \max_{t\in [0,T]} \norm{z_3(\cdot,t)}_{B^{-1+\tau}_{\infty,\infty}}\lesssim \delta_{q+1}\lambda_{q+1}^{-11\tau}.\]
Regarding the $C^1$ norm, we have
\[\normxt{1}{S_5}\leq\max_{t\in[0,T]}\norm{S_5(\cdot,t)}_{C^{1+\tau}}+\max_{t\in[0,T]}\norm{\partial_tS_5(\cdot,t)}_{C^{\tau}}.\]
For the first term, we use again the fact that $\mathcal{R}$ is an operator of order -1:
\[\norm{S_5(\cdot,t)}_{C^{1+\tau}}\lesssim \norm{z_3(\cdot,t)}_{C^{\tau}}\lesssim \delta_{q+1}\lambda_q\mu^\tau\lesssim \delta_{q+2}\lambda_{q+1}^{1-11\tau}.\]
Concerning the time derivative, we write
\[\partial_tz_3=\sum_{l\in \ZZ^3\backslash 0}\sum_{j=1}^6c_{jl}[\partial_t\Div(\overline{\gamma}_{qj}^2\zeta_j\otimes \zeta_j)]e^{i\mu l\cdot x}.\]
Thus, the same argument as with $z_3$ applies, albeit with an extra factor $\mu$ because of the extra derivative on the coefficients. We conclude
\begin{equation}
\normxt{0}{S_5}+\lambda_{q+1}^{-1}\normxt{1}{S_5}\lesssim \delta_{q+2}\lambda_{q+1}^{-11\tau}.
\label{cotas S5}
\end{equation}

Collecting the estimates for the error terms $S_1, \dots S_6$ that we have obtained in this section and applying \eqref{desigualdad norma tres}, we finally obtain
\begin{equation}
	\normxt{0}{R_{q+1}}+\lambda_{q+1}^{-1}\normxt{1}{R_{q+1}}\lesssim \delta_{q+2}\lambda_{q+1}^{-11\tau}.
	\label{cotas primer error}
\end{equation}
Therefore, the new Reynolds stress is small enough, since \eqref{inductive cota Rq C0} and \eqref{inductive cota Rq C1} hold for $q+1$. We are not done, however, because $R_{q+1}$ is not compactly supported. This is fixed in the following section.

\section{Final correction and conclusion}\label{sec correction}
In the previous sections we have checked that the new subsolution $(v_{q+1},B_{q+1},p_{q+1},R_{q+1})$ satisfies the inductive hypotheses (\ref{inductive tamaño Bq})-\eqref{cambio helicidad q} with $q$ replaced by $q+1$. In addition, the support of $(v_{q+1},B_{q+1},p_{q+1})$ is contained in $\overline{\Omega}_{q,2}\times[0,T]$.
Nevertheless, the new Reynolds stress is not compactly supported, so we have to perform a further correction, as described in \cref{subsec final correction}. 

The additional perturbation $w_L$ is chosen so that the field $\tilde{v}_{q+1}\coloneqq v_{q+1}+w_L$ satisfies the compatibility condition (\ref{momento angular cero}). Let us estimate this new term. The first step is to study the angular momentum of the perturbation $w_0+w_c$: 
\begin{lemma}
Given $N\geq 0$, the function $L\in C^\infty([0,T])$ defined in (\ref{def L}) satisfies 
\begin{equation}
\max_{t\in[0,T]}\abs{\partial_t^NL(t)}\lesssim \delta_{q+2}\lambda_{q+1}^{N-1-11\tau}.
\label{cotas L}
\end{equation}
\end{lemma}
\begin{proof}
We use (\ref{fórmula w0+wc}) and (\ref{identidad momento angular}) to write $L$ as
\[L=\sum_{j,m}\int_{\RR^3}\frac{2a_{jm}}{\ell_m\lambda_{q+1}}(\zeta_j\times k_{jm})\sin(\theta_{jm}).\]
Using the fact that $\partial_t\theta_{jm}$ is constant for any $m\in\Lambda$, we can use Leibniz's rule to write
\[\partial_t^NL=\sum_{j,m}\sum_{M=0}^N\int_{\RR^3}\begin{pmatrix}N\\M\end{pmatrix}\frac{2\,\partial_t^{N-M}a_{jm}\,(\partial_t\theta_{jm})^M}{\ell_m\lambda_{q+1}}\,\sin\left(\theta_{jm}+M\frac{\pi}{2}\right).\]
Recall that at most 16 of the coefficients $a_{jm}$ are nonzero in a small neighborhood of any given point and that their support is contained in $B(0,\bar{r})\times[0,T]$. In addition, we have $\abs{\partial_t\theta_{jm}}\lesssim\lambda_{q+1}$. Thus, applying \cref{stationary phase lemma} leads to
\begin{align*}
\abs{\partial_t^NL}&\lesssim \sum_{M=0}^N\abs{\int_{\RR^3}\sum_{j,m}\frac{\partial_t^{N-M}a_{jm}\,(\partial_t\theta_{jm})^M}{\ell_m\lambda_{q+1}}\,\sin\left(\theta_{jm}+M\frac{\pi}{2}\right)}
\\&\lesssim \sum_{M=0}^N\max_{j,m}\normxt{M+1}{\partial_t^{N-M}a_{jm}}\lambda_{q+1}^{(-1+M)-(M+1)}\lesssim \delta_{q+1}^{1/2}\mu^{N+1}\lambda_{q+1}^{-2}.
\end{align*}
Substituting (\ref{def mu}) into the previous bound yields (\ref{cotas L}).
\end{proof}
Substituting these bounds, along with (\ref{estimates wphi}), into the definition of $w_L$, we obtain
\begin{equation}
\normxt{N}{w_L}\lesssim \delta_{q+2}\lambda_{q+1}^{N-3\tau} \qquad N=0,1,2,3,4.
\label{estimates wL}
\end{equation}
This readily implies (\ref{inductive tamaño vq}), (\ref{inductive cotas vq CN}) and (\ref{cambio v q}) for the modified field $\tilde{v}_{q+1}$, since they hold for $v_{q+1}$ and the supports of $w_L$ and $v_{q+1}$ are disjoint. Furthermore, since the supports of the perturbation $w_L$ and the magnetic field $B_{q+1}$ are disjoint, the cross-helicity is unchanged and (\ref{cambio helicidad q}) also holds for $\tilde{v}_{q+1}$. Regarding the energy, we have
\[\abs{\int_{\RR^3}\abs{\tilde{v}_{q+1}}^2-\int_{\RR^3}\abs{v_{q+1}}^2}\lesssim \normx{0}{v_{q+1}}\normx{0}{w_L}+\normx{0}{w_L}^2\lesssim \delta_{q+2}\lambda_{q+1}^{-3\tau}.\]
Combining this with \cref{lema energía} yields \eqref{inductive q energía} for sufficiently large $a>1$.

The only remaining issue is the Reynolds stress. It follows from (\ref{inductive tamaño vq}), (\ref{inductive cotas vq CN}) and (\ref{estimates wL}) that the matrix $S_8$ defined in (\ref{def S8}) satisfies
\begin{equation}
\normxt{0}{S_8}+\lambda_{q+1}^{-1}\normxt{1}{S_8}\lesssim \delta_{q+2}\lambda_{q+1}^{-3\tau}.
\label{cotas S8}
\end{equation}
Concerning the matrix $S_7$ defined in (\ref{def S7}), we exploit the scaling to write
\[S_7(x,t)=-\frac{1}{2}\mathcal{R}\curl(\chi_LL')(x,t)+\lambda_{q+1}^{3\tau}S_2(\lambda_{q+1}^{\tau}(x-x_q),t),\]
where $S_2$ was defined in (\ref{def S2}). Since $\mathcal{R}\curl$ is bounded on spatial Besov spaces, it follows from (\ref{cotas L}) and the definition of $\chi_L$ that
\begin{align*}
\normatres{\mathcal{R}\curl(\chi_LL')}_{\tau}+\lambda_{q+1}^{-1}\normatres{\mathcal{R}\curl(\chi_LL')}_{1+\tau}&\lesssim \normatres{\chi_LL'}_{\tau}+\lambda_{q+1}^{-1}\normatres{\chi_LL'}_{1+\tau} \\&\lesssim \delta_{q+2}\lambda_{q+1}^{-3\tau}.
\end{align*}
Using (\ref{cotas S2}) to estimate the other term, we conclude
\begin{equation}
\normatres{S_7}_{\tau}+\lambda_{q+1}^{-1}\normatres{S_7}_{1+\tau}\lesssim \delta_{q+2}\lambda_{q+1}^{-3\tau}.
\label{cotas S7}
\end{equation}
Since $w_L$ was constructed so that (\ref{momento angular suma matrices}) holds, by \cref{invertir divergencia matrices} and \cref{remark modificar matriz}, there exists a matrix $S_9\in C^\infty(\RR^3\times[0,T],\RR^{3\times 3}_\text{sym})$ such that $\Div S_9=0$ and
\[S_1+S_2+S_4+S_5+S_7+S_9\]
is supported in $\Omega_{q+1}\times[0,T]$. This matrix satisfies the estimates
\begin{align*}
\normatres{S_9}_{\tau}+\lambda_{q+1}^{-1}\normatres{S_9}_{1+\tau}&\lesssim \normatres{S_1+S_2+S_4+S_5+S_7}_{\tau}+\lambda_{q+1}^{-1}\normatres{S_1+S_2+S_4+S_5+S_7}_{1+\tau}\\&\lesssim \delta_{q+2}\lambda_{q+1}^{-3\tau},
\end{align*}
where we have used (\ref{cotas S2}), (\ref{cotas S1}), (\ref{cotas S4}), \eqref{cotas S5} and (\ref{cotas S7}). Collecting the estimates for all the error terms $S_1, \dots, S_9$ and using (\ref{desigualdad norma tres}), we conclude that for sufficiently large $a>1$ the final Reynolds stress satisfies
\[\normxt{0}{\widetilde{R}_{q+1}}+\lambda_{q+1}^{-1}\normxt{1}{\widetilde{R}_{q+1}}\leq \delta_{q+2}\lambda_{q+1}^{-2\tau}.\]

As discussed in \cref{subsec final correction}, $(\tilde{v}_{q+1},B_{q+1},p_{q+1},\widetilde{R}_{q+1})$ is a subsolution satisfying (\ref{inductive support q}). In this section we have checked that it satisfies \eqref{inductive tamaño Bq}-\eqref{cambio helicidad q}. This completes the proof of \cref{prop steps}.

\section{Proof of \cref{main theorem}}\label{section th main}
First of all, it is straightforward to check that the MHD equations (\ref{MHD}) are invariant under the following transformation: 
\begin{align*}
	v(x,t)&\mapsto \lambda v(\lambda x,\lambda^2t), \\ 
	B(x,t)&\mapsto \lambda B(\lambda x,\lambda^2t), \\ 
	p(x,t)&\mapsto \lambda^2 p(\lambda x,\lambda^2t).
\end{align*}
Therefore, we may assume that the given solenoidal field $\bar{B}_0$ satisfies
\begin{equation}
\normxt{1}{\bar{B}_0}\leq \frac{1}{4},
\label{assumption size}
\end{equation}
We simply have to scale with a suitable choice of $\lambda$ above, construct the desired weak solution and then undo the scaling. Hence, it suffices to prove the result assuming the previous bound.

Next, we let $\beta>\alpha$ be the parameter given by \cref{prop steps} and then choose $\tau>0$ sufficiently small and $a>1$ sufficiently large so that \cref{prop steps} holds. Let $\varepsilon>0$ be as in the statement of \cref{main theorem} and let $\bar{r}>0$ such that the support of $\bar{B}_0$ is contained in the ball $B(0,\bar{r}/2)$. Without loss of generality, we may assume that $\varepsilon<1<\bar{r}$. We denote
\[a_\varepsilon\coloneqq \frac{\varepsilon}{1+\bar{r}}, \qquad \bar{r}_\varepsilon\coloneqq a_{\varepsilon}^{\hspace{1pt}-1/4}\,\bar{r}\]
and we define the solenoidal fields
\begin{align}
B_0(x,t)&\coloneqq \bar{B}_0(x)-a_\varepsilon\;\bar{B}_0\left(\frac{\bar{r}}{\bar{r}_\varepsilon}x\right), \label{def B0} \\
v_0(x,t)&\coloneqq \left(1-s\delta_1\lambda_0^{-2\tau}\frac{t}{T}\right)B_0(x,t).
\end{align}
where $s\in(0,1)$ will be chosen later, while $\lambda_q, \delta_q$ for $q\geq 0$ are the parameters defined in (\ref{def lambdaq}), (\ref{def deltaq}), respectively. We see that the support of $v_0$ and $B_0$ is contained in $\Omega_0\times[0,T]$, where we have set $\Omega_0\coloneqq B(0,\bar{r}_\varepsilon/2)$.

The definition of $a_\varepsilon$ and $\bar{r}_\varepsilon$ is tailored to yield an appropriate difference between $B_0$ and $\bar{B}_0$. Indeed, by (\ref{assumption size}), we have 
\[\normx{0}{B_0-\bar{B}_0}\leq \frac{1}{2}a_\varepsilon, \qquad \left[B_0-\bar{B}_0\right]_{C^1_x}\leq \frac{1}{2}a_\varepsilon^{5/4}.\]
It follows from the definition of the seminorm that $[f]_{C^\alpha_x}\leq [f]_{C^1_x}\text{diam}(\supp f)^{1-\alpha}$. Therefore, by our choice of $a_\varepsilon$ and $\bar{r}_\varepsilon$, we have
\begin{align}
\begin{split}
\normx{\alpha}{B_0-\bar{B}_0}&\leq \normx{0}{B_0-\bar{B}_0}+\left[B_0-\bar{B}_0\right]_{C^1_x}\bar{r}_\varepsilon^{\hspace{1pt}1-\alpha}\leq\frac{1}{2}a_\varepsilon+\frac{1}{2}a_\varepsilon^{5/4} \hspace{1pt}\bar{r}_\varepsilon^{\hspace{1pt}1-\alpha}\\&\leq \frac{1}{2}a_\varepsilon+\frac{1}{2}a_\varepsilon^{5/4} \bar{r}_\varepsilon=\frac{1}{2}(1+\bar{r})a_\varepsilon=\frac{\varepsilon}{2}.
\label{comparar con barraB0}
\end{split}
\end{align}

Meanwhile, the purpose of the additional term in $B_0$ was to ensure that
\begin{equation}
\int_{\Omega_0}B_0(x,t)\cdot \xi(x)\;dx=0 \qquad \forall \xi\in \ker\nabla_{\text{sym}},\;\forall t\in[0,T].
\label{producto Killing B0}
\end{equation}
Indeed, since $B_0$ is a compactly supported solenoidal field, it is easy to see that the previous integral vanishes for constant fields $\xi(x)\equiv b\in \RR^3$. One only has to apply the divergence theorem to the identity $\Div[(b\cdot x)B_0]=b\cdot B_0$. Therefore, it suffices to evaluate the integral for the non-constant vectors in the basis (\ref{definition base Killing}). We see that for these fields we have\[a_\varepsilon^{1/4}\xi_i(x)=\xi_i\left(a_\varepsilon^{1/4}x\right).\]
Hence,
\begin{align*}
\int_{\RR^3}a_\varepsilon\,\bar{B}_0\left(\frac{\bar{r}}{\bar{r}_\varepsilon}x\right)\cdot \xi_i(x)\;dx&=\int_{\RR^3}a_\varepsilon^{3/4}\bar{B}_0\left(a_\varepsilon^{1/4}x\right)\cdot \xi_i\left(a_\varepsilon^{1/4}x\right)\;dx\\&=\int_{\RR^3}\bar{B}_0(y)\cdot \xi_i(y)\;dy
\end{align*}
for any vector $\xi_i$ in the basis (\ref{definition base Killing}). We conclude (\ref{producto Killing B0}). Hence, by \cref{invertir divergencia matrices}, there exists a smooth symmetric matrix $S\in C^\infty_c(\Omega_0\times[0,T],\RR^{3\times 3}_\text{sym})$ such that
\[\Div S=\partial_t v_0=-sT^{-1}\delta_1\lambda_0^{-2\tau}B_0.\]	
Defining
\[R_0\coloneqq S-\delta_1\lambda_0^{-2\tau}\left(2s\frac{t}{T}-s^2\delta_1\lambda_0^{-2\tau}\frac{t^2}{T^2}\right)B_0\otimes B_0,\]
we see that
\[\partial_t v_0+\Div(v_0\otimes v_0-B_0\otimes B_0)=\Div R_0.\]
In addition, since the vector fields $\partial_t+v_0$ and $B_0$ commute in $\RR^3\times [0,T]$, we infer (from the definition of Lie derivative) that
\[
\frac{\partial }{\partial s}[(\Phi^0_{s})_\ast\hspace{1pt}B_0]=0\hspace{1pt},
\]
where $\Phi^0_s$ is the flow defined by $\partial_t+v_0$. Setting $t=s$ (which solves the $t$-component of the ODE defining $\Phi^0_s$), this flow is precisely the non-autonomous flow $X^0_t$ of $v_0$, thus implying that
\[
(X^0_{t})_\ast\,B_0=B_0\,,
\]
which is Equation~\eqref{subs c}. Defining $p_0=0$, we conclude that $(v_0,B_0,p_0,R_0)$ is a subsolution satisfying (\ref{inductive support q}).

Let us now check the required estimates. It follows from (\ref{assumption size}) that (\ref{inductive tamaño Bq})-(\ref{inductive cotas vq CN}) hold because $\delta_0^{1/2}\lambda_0=a^{1-\beta}$ can be made arbitrarily large by further increasing $a>1$, if necessary. Concerning the Reynolds stress, it follows from \cref{invertir divergencia matrices} that
\[\normx{1+\tau}{S(\cdot,t)}\lesssim sT^{-1}\delta_1\lambda_0^{-2\tau}\normx{\tau}{B_0(\cdot,t)}\lesssim sT^{-1}\delta_1\lambda_0^{-2\tau}.\]
Since $S$ is does not depend on time, because neither does $B_0$, we deduce
\[\normxt{1}{S}\lesssim sT^{-1}\delta_1\lambda_0^{-2\tau}.\]
Taking $s>0$ sufficiently small so as to compensate the numerical constants, we conclude (\ref{inductive cota Rq C0}) and (\ref{inductive cota Rq C1}). Finally, since the energy profile
\[\int_{\RR^3}\left(\abs{v_0(x,t)}^2+\abs{B_0(x,t)}^2\right)dx=\left[1+\left(1-s\delta_1\lambda_0^{-2\tau}\frac{t}{T}\right)^2\right]\left(1+a_\varepsilon^{1/4}\right)\norm{\bar{B}_0}_{L^2(\RR^3)}^2\]
is strictly decreasing for $t\in[0,T]$, we can find a strictly decreasing smooth function $e\in C^\infty([0,T])$ satisfying (\ref{inductive q energía}).

We conclude that $(v_0,B_0,p_0,R_0)$ is a subsolution satisfying the inductive hypotheses (\ref{inductive support q})-(\ref{inductive q energía}). Applying \cref{prop steps} iteratively, we obtain a sequence of subsolutions $\{(v_q,B_q,p_q,R_q)\}_{q=0}^\infty$ satisfying these hypotheses and also (\ref{cambio v q}), (\ref{cambio B q}) and \eqref{cambio helicidad q}. From (\ref{cambio v q}) and (\ref{cambio B q}) we deduce that $v_q$ and $B_q$ converge uniformly to some continuous weakly divergence-free fields $v$ and $B$, respectively. Since $R_q$ converges uniformly to 0, we see that $(v,B)$ is a weak solution of MHD. Recall that (\ref{MHD b}) and (\ref{MHD c}) hold throughout the whole process, due to \cref{equivalencia ec B}.

Furthermore, it follows from (\ref{cambio v q}) that 
\[\normxt{\alpha}{v_{q+1}-v_q}\lesssim \normxt{0}{v_{q+1}-v_q}^{1-\alpha}\normxt{1}{v_{q+1}-v_q}^{\alpha}\lesssim \delta_{q+1}^{\frac{1-\alpha'}{2}}\left(\delta_{q+1}^{1/2}\lambda_{q+1}\right)^{\alpha}=\lambda_{q+1}^{\alpha-\beta}.\]
Since $\beta>\alpha$, we see that the right-most term defines a convergent series. Thus, $v_q$ converges to $v$ in $C^{\alpha}(\RR^3\times[0,T])$. A completely analogous argument holds for the magnetic field:
\[\normxt{\alpha}{B_{q+1}-B_q}\lesssim \lambda_{q+1}^{\alpha-\beta} \qquad \Rightarrow\qquad \normxt{\alpha}{B-B_0}\lesssim \sum_{q=0}^\infty\lambda_{q+1}^{\alpha-\beta}\lesssim a^{-b(\beta-\alpha)}<\frac{\varepsilon}{2}, \]
where we have assumed that $a>1$ is sufficiently large. Combining this with (\ref{comparar con barraB0}) yields property $(iii)$ in \cref{main theorem}. Meanwhile, it follows from (\ref{inductive support q}) that the support of $(v,B)$ is contained in $B(0,\bar{r})$. We conclude that \[(v,B)\in C^{\alpha}(\RR^3\times[0,T])\cap C^{\hspace{0.5pt}0}([0,T],L^2(\RR^3)).\]
In fact, taking the limit in (\ref{inductive q energía}) yields
\[\int_{\RR^3}\left(\abs{v(x,t)}^2+\abs{B(x,t)}^2\right)dx=e(t).\]
Since the energy profile $e(t)$ was chosen to be strictly decreasing, we see that property $(i)$ in \cref{main theorem} holds.

Regarding the cross-helicity, we have
\begin{align}
\int_{\RR^3}v(x,0)\cdot B(x,0)\,dx-\int_{\RR^3}v(x,T)\cdot B(x,T)\,dx &\geq \int_{\RR^3}v_0(x,0)\cdot B_0(x,0)\,dx-\int_{\RR^3}v_0(x,T)\cdot B_0(x,T)\,dx  \nonumber\\&\hspace{12pt}-2\max_{t\in[0,T]}\abs{\int_{\RR^3}v\cdot B-\int_{\RR^3}v_0\cdot B_0}. \label{comparar helicidad}
\end{align}
By our choice of $v_0$ and $B_0$, we have
\begin{align*}
\int_{\RR^3}v_0(x,0)\cdot B_0(x,0)\,dx-\int_{\RR^3}v_0(x,T)\cdot B_0(x,T)\,dx&=s\delta_1\lambda_0^{-2\tau}\norm{B_0}_{L^2(\RR^3)}^2\\&=s\left[1+\left(\frac{\varepsilon}{2}\right)^{1/4}\right]\norm{\bar{B}_0}_{L^2(\RR^3)}^2\delta_1\lambda_0^{-2\tau}.
\end{align*}
Meanwhile, it follows from (\ref{cambio helicidad q}) that
\begin{align*}
\abs{\int_{\RR^3}v\cdot B-\int_{\RR^3}v_0\cdot B_0}&\leq \sum_{q=0}^\infty\abs{\int_{\RR^3}v_{q+1}\cdot B_{q+1}-\int_{\RR^3}v_q\cdot B_q}\\&\leq \sum_{q=0}^\infty \delta_{q+2}=\delta_2\sum_{q=0}^\infty a^{2\alpha \left(b^{q+2}-b^2\right)}\lesssim \delta_2.
\end{align*}
We see that the third term on the right-hand side of (\ref{comparar helicidad}) can be made smaller than the difference between the other two terms by taking $\tau>0$ sufficiently small and $a>1$ sufficiently large. Therefore, the cross-helicity will be strictly larger at the initial time than at time $t=T$, so property $(ii)$ in \cref{main theorem} holds, too. This completes the proof.

\section*{Acknowledgements}

This work has received funding from the European Research Council (ERC) under the European Union's Horizon 2020 research and innovation programme through the grant agreement~862342 (A.E.). It is partially supported by the grants CEX2023-001347-S, RED2022-134301-T and PID2022-136795NB-I00 (A.E. and D.P.-S.) funded by MCIN/AEI, and Ayudas Fundaci\'on BBVA a Proyectos de Investigaci\'on Cient\'ifica 2021 (D.P.-S.). q.P.-T. was partially supported by an FPI grant CEX2019-000904-S-21-4 funded by MICIU/AEI/10.13039/501100011033 and by FSE+.

\appendix 
\section{Hölder and Besov spaces}\label{appendix a}

Given $f\in C^\infty_c(\RR^3)$, we denote the supremum norm as \[\norm{f}_0\coloneqq \sup_{x\in\RR^3}\abs{f(x)},\]
where $\abs{\cdot}$ denotes the absolute value when working with real-valued maps, the modulus when working with vector-valued maps and the operator norm when working with matrix-valued maps, that is,
\[\abs{M}\coloneqq \sup_{\abs{v}\leq 1}\abs{Mv}.\]
This is the norm that we will use whenever we are working with matrices. It is not very important, however, since all norms on a finite-dimensional vector space are equivalent.

Given an integer $N\geq 0$ and a bounded domain $\Omega\subset \RR^n$ with Lipschitz boundary, we consider the following seminorm in the space $C^N(\Omega)$:
\[[f]_{C^N(\Omega)}\coloneqq \max_{\abs{\beta}=N}\norm{\partial^\beta f}_0.\]
Here $\beta$ is a multi-index. The following expression defines a norm on $C^N(\Omega)$:
\[\norm{f}_{C^N(\Omega)}\coloneqq \sum_{j\leq N}[f]_{C^j(\Omega)}.\]

Next, given $\alpha \in (0,1)$, we define the Hölder space $C^{N+\alpha}(\Omega)$ to be the subset of functions in $C^N(U)$ for which the following quantity, which can be shown to define a seminorm, is finite:
\[[f]_{C^{N+\alpha}(\Omega)}\coloneqq \max_{\abs{\beta}=N}\sup_{x\neq y}\frac{\abs{\partial^\beta f(x)-\partial^\beta f(y)}}{\abs{x-y}^\alpha}.\]
The space $C^{N+\alpha}(\Omega)$ becomes a Banach space when equipped with the  norm
\[\norm{f}_{C^{N+\alpha}(\Omega)}\coloneqq \norm{f}_{C^N(\Omega)}+[f]_{C^{N+\alpha}(\Omega)}.\] 
To simplify the notation, we will omit the domain if it can be inferred from the context. In fact, most of the time we work with compactly supported smooth functions whose support is contained in $U\times [0,T]$ for a fixed bounded domain with smooth boundary $U\subset \RR^3$ and some $T>0$. To make expressions even more compact, we denote the space-time Hölder norm by $\normxt{N+\alpha}{\cdot}$, since we use this norm very often. To distinguish it from the spatial Hölder norm in each time slice, which is used more sparingly, we add a subscript $x$, i.e.,
\begin{align*}
\normxt{N+\alpha}{f}&\coloneqq \norm{f}_{C^{N+\alpha}(U\times[0,T])}, \\
\normx{N+\alpha}{f(\cdot,t)}&\coloneqq \norm{f(\cdot,t)}_{C^{N+\alpha}(U)}, \qquad t\in [0,T].
\end{align*}

Unfortunately, some elements of our construction are not well-suited to work with space-time norms, so we need to introduce the following norms:
\[\normatres{f}_{N+\tau}\coloneqq \sup_{t\in[0,T]}\sum_{M_1+M_2\leq N}\normx{M_2+\tau}{\partial_t^{M_1}f(\cdot,t)},\]
where $N, M_1, M_2\geq 0$ are non-negative integers and $\tau\in(0,1)$. Note that we have 
\begin{equation}
\normxt{N}{f}\lesssim \normatres{f}_{N+\tau},
\label{desigualdad norma tres}
\end{equation}
so these norms are slightly stronger than space-time $C^N$ norms. The reason for which we introduce these norms is that, when estimating some time-dependent solutions to certain equations, it will be useful to keep a Hölder exponent in the spatial estimates.

We now summarize some properties of Hölder spaces. For a proof, see~\cite{CDK}.
\begin{proposition}
Let $\Omega\subset \RR^n$ be a bounded open Lipschitz set. Let $N_1, N_2, N_3\geq 0$ be integers and $0\leq \alpha_1,\alpha_2,\alpha_3< 1$ such that
\[N_1+\alpha_1\leq N_2+\alpha_2\leq N_3+\alpha_3.\]
Let $\theta\in [0,1]$ such that
\[N_2+\alpha_2=\theta(N_1+\alpha_1)+(1-\theta)(N_3+\alpha_3).\]
Then, there exists a constant $C=C(\Omega,N_3)>0$ such that
\[\normxt{N_2+\alpha_2}{f}\leq C\normxt{N_1+\alpha_1}{f}^{\theta}\normxt{N_3+\alpha_3}{f}^{1-\theta}.\]
\end{proposition}
\begin{proposition}
Let $\Omega\subset \RR^n$ be a bounded open Lipschitz set, $N\geq 0$ be an integer and $0\leq \alpha< 1$. Let $f,g\in C^{N+\alpha}(\Omega\times[0,T])$. Then, there exists a constant $C_1=C_1(\Omega, N)>0$ such that 
\begin{equation}
	\label{estimate product}
	\normxt{N+\alpha}{fg}\leq C_1(\normxt{0}{f}\normxt{N+\alpha}{g}+\normxt{N+\alpha}{f}\normxt{0}{g}).
\end{equation}
In addition, if $c^{-1}\leq f\leq c$ on $\Omega$ for some constant $c>0$, there exists a constant $C_2=C_2(\Omega,c,N)$ such that
\begin{equation}
	\label{estimate quotient}
	\normxt{N+\alpha}{\frac{1}{f}}\leq C_2 \normxt{N+\alpha}{f}.
\end{equation}
\end{proposition}

\begin{proposition}
Let $\Omega\subset \RR^n$ be a bounded open Lipschitz set, $N\geq 0$ be an integer and $0\leq \alpha< 1$. Let $\Phi\in C^{N+\alpha}$ be a homeomorphism of $\Omega$. Suppose that
\[\normxt{1}{\Phi}+\normxt{1}{\Phi^{-1}}\leq c\]
for some $c>0$. Then, there exists a constant $C=C(c,N,\Omega)$ such that
\begin{equation}
\normxt{N+\alpha}{\Phi^{-1}}\leq C\normxt{N+\alpha}{\Phi}.
\label{estimate inverse}
\end{equation}
\end{proposition}

\begin{proposition}
Let $f:\Omega\to \RR$ and $\Psi:\RR^n\to \Omega$ be smooth functions with $\Omega\subset \RR^m$. Then, for every $N\in \NN\backslash \{0\}$ there is a constant $C=C(m,n,N)$ such that
\begin{equation}
\label{estimate composition}
[f\circ \Psi]_{C^N}\leq C\left([f]_{C^1}\norm{D\Psi}_{C^{N-1}}+\norm{D f}_{C^{N-1}}[\Psi]_{C^1}^N\right).
\end{equation}
\end{proposition}

The following proposition is a simplified version of the results found in~\cite{CDK}. Since this result could be considered to be less standard than the others, and we have adapted the statement to suit our needs, we include the proof for convenience of the reader. Check \cref{section preliminaries} for the notation $\compt$ and $\pi$.
\begin{proposition}
\label{prop diferencia composición}
Let $U\subset \RR^3$ be a bounded domain with smooth boundary and let $T>0$. Let $N\geq 1$ be an integer. Consider a function $f\in C^N(U\times[0,T])$ and a map $\psi\in C^\infty(U\times[0,T],U)$. Then, there exists a constant $C=C(N,U,T)$ such that
\begin{align*}
\normxt{0}{f\compt\psi-f}&\leq \normxt{1}{f}\normxt{0}{\psi-\pi}, \\
\normxt{N}{f\compt\psi-f}&\leq C\normxt{1}{f}\normxt{N}{\psi-\pi}+C \normxt{0}{\psi-\pi}\normxt{2}{f}\left(1+\normxt{N}{\psi-\pi}\right)\\&\hspace{12pt}+C\normxt{0}{\psi-\pi}\normxt{N+1}{f}\left(1+\normxt{1}{\psi-\pi}\right)^N.
\end{align*}
\end{proposition}
\begin{proof}
The bound for the $C^0$-norm is straightforward, so let us focus on the other one. Note that we can write
\begin{align*}
(f\compt\psi)(x,t)-f(x,t)&=\int_0^1\frac{d}{ds}\left[f(x+s(\psi(x,t)-x),t)\right]\,ds \\ &=\int_0^1Df(x+s(\psi(x,t)-x),t)\,(\psi(x,t)-x)\;ds.
\end{align*}
By (\ref{estimate product}), we then have
\[\normxt{N}{f\compt\psi-f}\lesssim\normxt{1}{f}\normxt{N}{\psi-\pi}+\normxt{0}{\psi-\pi}\max_{s\in[0,1]}\normxt{N}{Df\compt(\pi+s(\psi-\pi))}.\]
Using (\ref{estimate composition}) to estimate the last term yields the desired bound.
\end{proof}

To estimate the $C^0$-norm of the solution to certain differential equations, we need to introduce an additional functional space. Let us consider a Littlewood--Payley decomposition, e.g.\ as in~\cite[Section 2.2]{BCD}. For this, we take smooth radial functions $\chi,\varphi:\RR^3\to[0,1]$, whose supports are contained in the ball $B(0,\frac43)$ and in the annulus $\{\frac34 <|\xi|<\frac83\}$ respectively, with the property that
\[
\chi(\xi)+\sum_{N=0}^\infty\varphi(2^{-N}\xi)=1
\]
for all~$\xi\in\RR^3$. Given $s\in\RR$, the Besov norm $B^s_{\infty,\infty}$ can be written in terms of the Fourier multipliers $P_{<}:=\chi(-i\nabla)$ and $P_N:=\varphi(-2^{-N}i\nabla)$ as 
\begin{equation}\label{E.BesovF}
	\|h\|_{B^{s}_{\infty,\infty}}:=\|P_{<}h\|_{L^\infty}+\sup_{N\geq0}2^{Ns}\|P_Nh\|_{L^\infty}.
\end{equation}
It can be proved that this norm is equivalent to the H\"older norm $C^s$ if $s\in\RR^+\backslash\NN$, and strictly weaker if $s\in\NN$.

The following lemma shows that we can estimate the Besov norms of a type of function that is widely used throughout the paper: 
\begin{lemma}
	\label{stationary phase lemma}
	Let $s>-1$, $J\in \NN$ and $\lambda\geq 1$. For $j=1, \dots J$, let $c_j\in C^\infty_c(\RR^3)$ and $k_j\in \RR^3$ such that $C_0^{-1}\leq \abs{k_j}\leq C_0$ for some $C_0>0$. Suppose that the supports of at most $J_*$ of the functions $\{c_j\}_{j=1}^q$ have a nonempty intersection, so that the product of any $J_*+1$ of the functions~$c_j$ is identically~$0$. Let
	\[f(x)\coloneqq \sum_{j=1}^J c_j(x)e^{i\lambda k_j\cdot x}.\]
	Then, for any nonnegative integer $m$ we have
	\begin{align*}
	\abs{\int_{\RR^3}f}&\leq C_0\hspace{1pt}J_*\frac{\mathcal{C}_m}{\lambda^m}\abs{\supp f},\\[3pt]
	\norm{f}_{B^{s}_{\infty,\infty}}&\leq C_*\left(\lambda^s\mathcal C_0+\lambda^{-m}\mathcal C_{m+s_+}\right),
	\end{align*}
	where $s_+:=\max\{s,0\}$ and $\mathcal C_N:=\max_{j\leq q}\|c_j\|_N$. The constant~$C_*$ depends on $s,m, J_*$ and $C_0$ but not on~$J$.
\end{lemma}
\begin{proof}
The above bound for the Besov norm is proved in \cite{steady_Euler}. Regarding the bound for the integral, it is immediate to adapt the argument from \cite{Continuous} to our setting. We define
\begin{align*}
A_n&\coloneqq\sum_{j=1}^J -i\left[\frac{k_j}{\abs{k_j}^2}\left(i\frac{k_j}{\abs{k_j}^2}\cdot\nabla\right)^nc_j\right]e^{i\lambda k_j\cdot x}, \\
F_n&\coloneqq\sum_{j=1}^J \left[\left(i\frac{k_j}{\abs{k_j}^2}\cdot\nabla\right)^nc_j\right]e^{i\lambda k_j\cdot x}.
\end{align*}
A straightforward calculation shows that
\[F_n=\frac{1}{\lambda}\Div A_n+\frac{1}{\lambda}F_{n+1}.\]
In particular, for any $m\geq 0$:
\[f=F_0=\sum_{n=0}^{m-1}\frac{1}{\lambda^{n+1}}\Div A_n+\frac{1}{\lambda^m}F_m.\]
Integrating this expression over the support of $f$ yields the desired bound.
\end{proof}

\section{The divergence equation} \label{appendix divergence}
We recall some tools concerning the divergence equation for vector fields and matrix-valued fields. The main issue is whether we can obtain compactly supported solutions when the source term is compactly supported. This is possible if and only if the source satisfies certain compatibility conditions. 

Let us start with the case of matrices, which is more complicated. We denote by $\RR^{3\times 3}_\text{sym}$ the set of symmetric $3\times 3$ matrices. We introduce the following operators:
\begin{align*}
	\Div: C^\infty(\RR^3,\RR^{3\times 3}_\text{sym})&\to C^\infty(\RR^3,\RR^3), \hspace{20.5pt} S\mapsto v_i\coloneqq\partial_j S_{ij}, \\
	\nabla_{\text{sym}}: C^\infty(\RR^3,\RR^3)&\to C^\infty(\RR^3,\RR^{3\times 3}_\text{sym}), \hspace{12pt} v\mapsto S_{ij}\coloneqq\frac{1}{2}(\partial_i v_j+\partial_j v_i).
\end{align*}
We use the convention that partial derivatives with latin subscripts denote partial derivatives in the spatial coordinates, whereas temporal partial derivatives are always denoted by $\partial_t$. Using the divergence theorem, it is easy to see that the following identity holds in any bounded domain with smooth boundary $\Omega\subset\RR^3$:
\begin{equation}
	\int_{\Omega}\xi\cdot \Div S+\int_{\Omega} \nabla_{\text{sym}}\xi:S=\int_{\partial \Omega}\xi^t S\,\nu.
	\label{identidad de Green matrices}
\end{equation}

The kernel of the operator $\nabla_{\text{sym}}$ are the so-called Killing vector fields. It is a finite-dimensional vector space that plays an important role in Riemannian geometry. It is well known (see \cite[page 52]{Petersen}) that in $\mathbb{R}^3$ a basis of this vector space is given by
\begin{equation}
	\label{definition base Killing}
	\mathcal{B}\coloneqq \{e_1, e_2, e_3, \xi_{1}, \xi_{2}, \xi_{3}\},
\end{equation}
where $\{e_1,e_2,e_3\}$ is the Cartesian basis and where
\begin{equation}
	\label{definition Killing momento angular}
	\xi_{i}(x)\coloneqq e_i\times x \hspace{40pt} i=1,2,3.
\end{equation}

Next, we introduce a differential operator that is standard in convex integration:
\begin{definition}
	Let $f\in C^\infty_c(\RR^3,\RR^3)$ be a smooth vector field. We define a matrix valued function $\mathcal{R}f$ as
	\begin{equation}
		\label{def mathcal R}
		(\mathcal{R}f)_{ij}\coloneqq \Delta^{-1}(\partial_i f_j+\partial_j f_i)-\delta_{ij}\Delta^{-1}\Div f.
	\end{equation}
	Here $\delta_{ij}$ is the Kronecker delta and $\Delta^{-1}$ refers to the potential-theoretic solution of the Poisson equation, that is, the (spatial) convolution of the source term  with the fundamental solution of the Laplace equation in $\mathbb R^3$. 
\end{definition}
\begin{lemma}
	For any $f\in C^\infty_c(\RR^3\times[0,T],\RR^3)$, we have
	\begin{enumerate}
		\setlength\itemsep{5pt}
		\item $\mathcal{R}f(x,t)$ is a symmetric matrix for all $(x,t)\in \RR^3\times[0,T]$,
		\item $\Div (\mathcal{R}f)=f$,
		\item $\mathcal{R}f\in C^\infty(\RR^3\times[0,T],\RR^{3\times3}_{\text{sym}})$ and for any $t\in[0,T]$ we have
		\begin{equation}
			\normx{N+\tau}{\partial_t^M\mathcal{R}f(\cdot,t)}\leq C(N,\tau)\norm{\partial_t^Mf(\cdot,t)}_{B^{N-1+\tau}_{\infty,\infty}}
			\label{R acotado en Besov}
		\end{equation}
		for any integers $N,M\geq 0$ and $\tau\in(0,1)$.
	\end{enumerate}
\end{lemma}
\begin{proof}
	It is clear from the definition that $\mathcal{R}f$ is symmetric. Meanwhile, $(ii)$ is easily obtained by taking the divergence in the definition of $\mathcal{R}f$ and using that $\Delta(\Delta^{-1}\rho)=\rho$. Regarding $(iii)$, since $\mathcal{R}$ is an operator of order $-1$, we have
	\[\normx{N+\tau}{\mathcal{R}f(\cdot,t)}\leq C(N,\tau)\norm{f(\cdot,t)}_{B^{N-1+\tau}_{\infty,\infty}}.\]
	On the other hand, since $\Delta^{-1}\rho$ is the unique solution of $\Delta u=\rho$ in $L^2(\RR^3)$, we deduce that $\Delta^{-1}$ and $\partial_t$ commute. As a result, $\partial_t$ and $\mathcal{R}$ commute, too. Hence, applying the previous estimate to $\partial_t^M f$ yields (\ref{R acotado en Besov}) and, in particular, smoothness in time of $\mathcal{R}v$.
\end{proof}

It follows from the previous lemma that the operator $\mathcal{R}$ can be used to solve the divergence equation $\Div S=f$. However, $\mathcal{R}f$ will not be compactly supported, in general. Nevertheless, it can be modified to obtain a compactly supported solution, provided that $f$ satisfies the appropriate conditions:
\begin{lemma}
	\label{invertir divergencia matrices}
	Let $\Omega\subset \RR^3$ be a bounded domain with smooth boundary and let $f\in C^\infty_c(\Omega,\RR^3)$. There exists $S\in C^\infty_c(\Omega,\RR^{3\times3}_{\text{sym}})$ such that $\Div S=f$ if and only if
	\begin{equation}
	\int_{\Omega} f\cdot\xi\,dx=0 \qquad \forall \xi\in \ker \nabla_{\text{sym}},\; \forall t\in[0,T].
	\label{compatibilidad divergencia matrices}
	\end{equation}
	In that case, $S$ may be chosen so that for any integers $N,M\geq 0$ and $\tau\in(0,1)$ we have
	\begin{equation}
	\normx{N+\tau}{\partial^M_tS(\cdot,t)}\lesssim\norm{\partial_t^Mf(\cdot,t)}_{B^{N-1+\tau}_{\infty,\infty}},
	\label{cotas invertir div matrices}
	\end{equation}
	where the implicit constants depend on $N$, on $\alpha$ and on $\Omega$.
\end{lemma}
\begin{proof}
In \cite[Lemma 2.9]{Extension} it is shown that condition (\ref{compatibilidad divergencia matrices}) is necessary and, provided that it holds, it explains how to construct a matrix $S\in C^\infty_c(\Omega,\RR^{3\times3}_{\text{sym}})$ such that $\Div S=f$ and satisfying the estimates
\[\normx{N+\tau}{S(\cdot,t)}\leq C(N,\tau,\Omega)\,\norm{f(\cdot,t)}_{B^{N-1+\tau}_{\infty,\infty}} \qquad t\in[0,T]\]
for any integer $N\geq 0$ and $\tau\in (0,1)$. Thus, we only need to estimate the time derivatives. Further inspection of the proof of \cite[Lemma 2.9]{Extension} shows that the construction commutes with taking the partial derivative. Hence, applying the result to $\partial^M_t f$ yields (\ref{cotas invertir div matrices}).
\end{proof}
\begin{remark}\label{remark modificar matriz}
Finding a solution $S\in C^\infty_c(\Omega\times[0,T],\RR^{3\times 3}_\text{sym})$ to $\Div S=f$ is equivalent to finding a divergence-free symmetric matrix $M$ such that the support of $\mathcal{R}f+M$ is contained in $\Omega\times[0,T]$. The previous lemma proves the necessary and sufficient conditions for this to be possible. In addition, further inspection of the proof of \cite[Lemma 2.9]{Extension} shows that one has the following estimates:
\[\normx{N+\tau}{\partial^M_tM(\cdot,t)}\lesssim \normx{N+\tau}{\partial^M_t\mathcal{R}f(\cdot,t)}\]
for any integers $N,M\geq 0$ and $\tau\in (0,1)$. 
\end{remark}
\begin{remark} \label{remark dependencia dominio}
In the previous results the implicit constants depend on the domain $\Omega$. However, since we only work on balls whose radii are in the interval $(\bar{r}/2, \bar{r})$ for some $\bar{r}>0$, we can choose common constants. Hence, we can essentially ignore the dependence on the domain. See~\cite[Lemma A.6]{steady_Euler} for the explicit dependence on the radius. 
\end{remark}

An analogous result holds for the vector-valued divergence equation, which is more standard:
\begin{lemma}\label{invertir div vectores}
	Let $T>0$ and let $\Omega\subset \RR^3$ be a bounded domain with smooth boundary. Given a density $f\in C^\infty_c(\Omega\times[0,T])$, there exists a vector field $u\in C^\infty_c(\Omega\times[0,T],\RR^3)$ such that $\Div u=f$ if and only if
	\[\int_\Omega f(x,t)\,dx=0 \qquad \forall t\in[0,T].\]
	In that case, $u$ can be chosen so for that for any integers $N,M\geq 0$ and $\tau\in(0,1)$ we have
	\begin{equation}
		\normx{N+\tau}{\partial^M_tu(\cdot,t)}\lesssim\norm{\partial_t^Mf(\cdot,t)}_{B^{N-1+\tau}_{\infty,\infty}},
		\label{cotas invertir div vectores}
	\end{equation}
	where the implicit constants depend on $N$, on $\alpha$ and on $\Omega$.
\end{lemma}
For the same reason as in \cref{remark dependencia dominio}, we can essentially ignore the dependence of the constants on the domain.

\section{Other auxiliary results}\label{app.C}
The following lemma, which is standard, shows the existence of cut-off functions with good bounds for their derivatives. For a proof, see e.g.~\cite[Lemma B.1]{Extension}.
\begin{lemma}
	\label{lemma cutoff}
	Let $A\subset \mathbb{R}^3$ be a measurable set and let $r>0$. There exists a cutoff function $\chi_r\in C^\infty(\mathbb{R}^3,[0,1])$ whose support is contained in $A+B(0,r)$ and such that $\chi_r\equiv 1$ in a neighborhood of $A$. Furthermore, for any $N\geq 0$ we have \[\left\|\chi_r\right\|_N\leq C(N)\,r^{-N}.\]
	for some universal constants $C(N)$ depending on $N$ but not on $A$ or $r$. 
\end{lemma}

For convenience of the reader, we give the proof of a standard result concerning the decomposition of the Reynolds stress into simpler matrices:
\begin{proof}[Proof of \cref{geometric lemma}]
	We define
	\begin{align*}
		\zeta_1 &\coloneqq \frac{1}{\sqrt{2}}(e_1+e_2), \qquad \qquad \zeta_2\coloneqq  \frac{1}{\sqrt{2}}(e_2+e_3), \qquad \qquad \zeta_3\coloneqq \frac{1}{\sqrt{2}}(e_3+e_1), \\  \zeta_4&\coloneqq \frac{1}{\sqrt{2}}(e_1-e_2), \qquad \qquad
		\zeta_5\coloneqq\frac{1}{\sqrt{2}}(e_2-e_3), \qquad \qquad \zeta_6\coloneqq \frac{1}{\sqrt{2}}(e_3-e_1).
	\end{align*}
	One can check that $\{\zeta_j\otimes\zeta_j:j=1,\dots, 6\}$ is a basis of the vector space $\RR^{3\times 3}_{\text{sym}}$. Therefore, the map  $\RR^6\to\RR^{3\times 3}_{\text{sym}},\;\; (x_1,\dots, x_6)\mapsto \sum_{j=1}^6x_j\zeta_j\otimes\zeta_j$ is a linear isomorphism between both vector spaces. Since
	\[\Id=\sum_{j=1}^6\frac{1}{2}\zeta_j\otimes\zeta_j,\]
	the coordinates of the inverse map must be positive in a sufficiently small neighborhood of $\Id$, so we can take the square root and the claim follows.
\end{proof}

\bibliographystyle{amsplain}

\end{document}